\documentclass[final,3p,sort&compress]{elsarticle}

\usepackage[T1]{fontenc}
\usepackage[bitstream-charter]{mathdesign}

\usepackage[english]{babel}
\usepackage[colorlinks=true]{hyperref}

\usepackage{caption}
\usepackage{subcaption}

\usepackage{amsmath,amsthm}
\usepackage{mathtools}
\usepackage{commath}
\usepackage{enumerate}

\usepackage{graphicx}
\usepackage{booktabs}
\usepackage[nameinlink]{cleveref}
\usepackage{comment}

\usepackage{bm}
\usepackage{stmaryrd}
\SetSymbolFont{stmry}{bold}{U}{stmry}{m}{n}
\DeclareMathAlphabet{\pazocal}{OMS}{zplm}{m}{n}
\DeclareMathAlphabet{\pazocalbf}{OMS}{cmsy}{b}{n}
\usepackage{etoolbox}

\preto\subequations{\ifhmode\unskip\fi}
\theoremstyle{plain}

\theoremstyle{definition}

\newtheorem{assumption}{Assumption}

\theoremstyle{remark} 
\newtheorem{remark}{Remark}

\newcommand{\Ohn}{\Omega_h^n}
\newcommand{\Ohndt}{\widetilde{\Omega}_{h,\delta_h}^n}
\newcommand{\Ohnt}{\widetilde{\Omega}_{h,0}^n}
\newcommand{\Ghn}{\Gamma_h^n}

\newcommand{\Tht}{\widetilde{\mathcal{T}}_h}
\newcommand{\Thn}{\mathcal{T}_{h,0}^n}
\newcommand{\Thnd}{\mathcal{T}_{h,\delta_h}^n}
\newcommand{\ThnS}{\mathcal{T}_{h,\pazocal{S}}^n}
\newcommand{\ThnG}{\mathcal{T}_{h,\Gamma}^n}

\newcommand{\Fht}{\widetilde{\mathcal{F}}_{h}}
\newcommand{\Fhnd}{\mathcal{F}^n_{h,\delta_h}}
\newcommand{\FhndAl}{\mathcal{F}^{n,\text{Al}}_{h,\delta_h}}

\newcommand{\FhnG}{\mathcal{F}^n_{h,\Gamma}}
\newcommand{\FhnGAl}{\mathcal{F}^{n,\text{Al}}_{h,\Gamma}}

\newcommand{\ub}{\bm{u}}
\newcommand{\bu}{\bm{u}}
\newcommand{\vb}{\bm{v}}
\newcommand{\bv}{\bm{v}}
\newcommand{\wb}{\bm{w}}
\newcommand{\xb}{\bm{x}}
\newcommand{\nb}{\bm{n}}
\newcommand{\blf}{\bm{f}}
\newcommand{\bd}{\bm{d}}

\newcommand{\uhn}{\ub_h^n}
\newcommand{\vh}{\vb_h}

\newcommand{\uninf}{u_\Gamma^\infty}

\newcommand{\dn}{\partial_{\nb}}
\DeclareMathOperator{\diver}{div}
\DeclareMathOperator{\dist}{dist}
\DeclareMathOperator{\meas}{meas}

\newcommand{\jump}[1]{\llbracket#1\rrbracket}
\newcommand{\restr}[2]{{\left.\kern-\nulldelimiterspace#1\right|_{#2}}}

\newcommand{\gpv}{\gamma_{\text{gp}}^{v}}
\newcommand{\gpp}{\gamma_{\text{gp}}^{p}}
\newcommand{\ggd}{\gamma_{\text{gd}}}

\journal{Journal Name}

\begin{document}

\begin{frontmatter}
\title{Stability of instantaneous pressures in an Eulerian finite element method for moving boundary flow problems}

\author[1]{Maxim Olshanskii}
\ead{molshan@math.uh.edu}
\author[2]{Henry von Wahl}
\ead{henry.von.wahl@uni-jena.de}

\affiliation[1]{
  organization={Department of Mathematics, University of Houston},
  addressline={651 PGH},
  city={Houston},
  postcode={77204},
  state={TX},
  country={USA}}

\affiliation[2]{
  organization={Friedrich-Schiller-Universität, Fakultät für Mathematik und Informatik},
  addressline={Ernst-Abber-Platz 2},
  city={Jena},
  postcode={07743},
  country={Germany}}

\begin{abstract}
This paper focuses on identifying the cause and proposing a remedy for the problem of spurious pressure oscillations in a sharp-interface immersed boundary finite element method for incompressible flow problems in moving domains.
The numerical method belongs to the class of Eulerian unfitted finite element methods. It employs a cutFEM discretization in space and a standard BDF time-stepping scheme, enabled by a discrete extension of the solution from the physical domain into the ambient space using ghost-penalty stabilization.
To investigate the origin of spurious temporal pressure oscillations, we revisit a finite element stability analysis for the steady domain case and extend it to derive a stability estimate for the pressure in the  $L^\infty(L^2)$-norm that is uniform with respect to discretization parameters. By identifying where the arguments fail in the context of a moving domain, we propose a variant of the method that ensures unconditional stability of the instantaneous pressure. As a result, the modified method eliminates spurious pressure oscillations.
We also present extensive numerical studies aimed at illustrating our findings and exploring the effects of fluid viscosity, geometry approximation, mass conservation, discretization and stabilization parameters, and the choice of finite element spaces on the occurrence and magnitude of spurious temporal pressure oscillations.
The results of the experiments demonstrate a significant improvement in the robustness and accuracy of the proposed method compared to existing approaches.
\end{abstract}


\end{frontmatter}

\section{Introduction}
\label{sec.intro}

Fluid equations defined in domains with moving boundaries and interfaces are central to mathematical models across diverse applications, such as cardiovascular modeling and aerospace engineering~\cite{formaggia2010cardiovascular,bazilevs2013computational}. Developing computational methods for solving such problems numerically has been recently a focus of extensive research. Established techniques include immersed boundary methods, fictitious domain methods, and approaches based on Lagrangian and arbitrary Lagrangian-Eulerian formulations; see, e.g., \cite{hirt1974arbitrary,peskin1977numerical,tezduyar1992new,masud1997space,glowinski1999distributed,formaggia1999stability,duarte2004arbitrary,gross2011numerical}.

In this paper, we are interested in an Eulerian finite element method that employs a fixed triangulation of $\mathbb{R}^d$ to solve a system of governing equations within a domain $\Omega(t)\subset\mathbb{R}^d$, where the boundary $\Gamma(t)$ moves smoothly through the static background mesh. This renders the method an immersed boundary method (IBM), an important class of techniques that simplifies handling complex deformations and topological transitions in $\Omega(t)$ by avoiding the need to fit the mesh to the evolving geometry.

While using background meshes decoupled from boundary motion offers clear algorithmic benefits and implementation conveniences, this approach also presents challenges. When applied to numerically solve the incompressible Navier-Stokes equations, immersed boundary methods are known to produce spurious \emph{temporal} pressure oscillations, often resulting in non-physical instantaneous forces (such as drag and lift) acting on immersed moving bodies; see, e.g.,~\cite{uhlmann2005immersed,seo2011sharp}.  Moreover, the amplitude of spurious oscillations increases when the specific CFL number, 
\begin{equation}\label{eq:CFL}
C_{\text{CFL}} = (V_\Gamma \Delta t) / h
\end{equation}
--- the ratio between the boundary speed, time step, and local mesh size --- becomes \emph{small}~\cite{lee2011sources,XCM24}.

The phenomenon of spurious pressure oscillations in immersed boundary methods has been studied extensively in the literature from various perspectives, leading to several interpretations and amelioration techniques~\cite{uhlmann2005immersed,mittal2008versatile,luo20103d,liao2010simulating,yang2009smoothing,shirgaonkar2009new,lee2011sources,seo2011sharp,lee2013implicit,schneiders2013accurate,ruberg2014fixed,goza2016accurate,kim2019immersed,XCM24}. For example, \cite{uhlmann2005immersed} used a diffusive interface IBM and an explicit formulation of the fluid-solid interaction forces to reduce spurious forces. In \cite{yang2009smoothing}, small but non-negligible spurious force oscillations were observed for an IBM with direct volume forcing. These oscillations were mitigated by smoothing the discrete delta function used in the IBM. Similarly, \cite{shirgaonkar2009new} applied an IBM in the context of a constraint-based formulation for the problem of self-propulsion, where force oscillations were also evident.

On the other hand, \cite{mittal2008versatile} developed a sharp-interface IBM, which, while effective, suffers from spurious pressure oscillations over time~\cite{luo20103d,liao2010simulating}. The lack of geometric (mass) conservation was identified as a source of these oscillations~\cite{seo2011sharp}, and an IBM variant strictly enforcing conservation of area was proposed, reducing oscillations by an order of magnitude. In \cite{lee2011sources}, velocity and pressure discontinuities, which occur when a grid point transitions into the fluid domain, were identified as another source of spurious oscillations. To address this, the authors proposed adding numerical mass sources and sinks. Additionally, \cite{lee2013implicit} introduced a fully implicit backward differentiation approach to handle the effects of multiple layers of new cells entering the fluid region.

In \cite{schneiders2013accurate}, abrupt changes in discretization operators caused by moving boundaries were identified as a key source of spurious force oscillations. To counteract this, the authors developed a cut-cell method and a moving boundary formulation that responds smoothly to these changes, avoiding discontinuities arising from the transition of cells into or out of the fluid domain. This formulation conserves mass, momentum, and energy through a flux-redistribution technique. In \cite{goza2016accurate}, the authors identified the ill-posedness of surface stress equations as a source of inaccurate stresses and forces in IBMs with smeared delta functions. They demonstrated that smoothing the delta function, combined with a filtering technique, significantly improves the representation of physical stresses. A detailed review of IBMs for fluid-structure interaction (FSI) problems, along with various approaches to mitigate spurious force oscillations, is presented in \cite{kim2019immersed}.

The problem of spurious pressure oscillations has been studied less extensively in the context of unfitted finite element methods for fluid problems with moving boundaries or interfaces. In \cite{ruberg2014fixed}, a fixed-grid B-spline finite element method was employed to discretize the fluid-structure interaction problem. The robustness of this scheme was enhanced through a Dirichlet-Robin coupling algorithm and a conservative operator to transfer data between the fluid and the solid. More recently, \cite{XCM24} introduced a weighted shifted boundary method for Stokes flow, achieving higher accuracy in fluid volume preservation compared to the standard shifted boundary method. This approach resulted in significantly reduced pressure drag oscillations.

Despite this understanding, developing immersed boundary methods that are robust with respect to 
$C_{\text{CFL}}$  and employ a sharp treatment of the evolving geometry and boundary/interface forces remains a largely open problem. The use of unfitted finite elements, also known as cutFEM~\cite{BCH14}, is a promising approach, as many of the aforementioned issues can be addressed accurately. The cutFEM approach uses an unfitted background mesh to define finite element spaces. Dirichlet boundary conditions are enforced either as constraints in a variational setting, leading to Lagrange multipliers~\cite{BH10}, or using Nitsche's method~\cite{BH12}. Stability and conditioning problems in the resulting linear systems, due to elements with a small intersection with the physical domain, are handled using ghost-penalty stabilization~\cite{Bur10}. Integration over the resulting unfitted domains is addressed with custom quadrature approaches. Consequently, in the present setting of a fluid problem in a moving domain, the conservation of area (fluid volume) can be approximated to higher order in cutFEM, depending on the quadrature used, or even achieved exactly if the geometry allows for exact quadrature rules. The conservation of mass of the fluid can also be enforced to varying degrees in this context. For example, grad-div stabilization is a well-established tool for improving mass conservation in finite element methods that are not exactly divergence-free~\cite{olshanskii2002low}. Alternatively, finite element spaces can be constructed such that divergence condition is fulfilled pointwise, although, we note that limited work has been done in this direction in the unfitted setting; see \cite{LNO23,FNZ23}.

For time-dependent problems, two common approaches to employing cutFEM are using either space-time unfitted finite elements~\cite{hansbo2016cut} or a hybrid method that applies finite difference discretizations to the time derivatives~\cite{LO19}. In this paper, we focus on the latter approach, an extension of the method of lines for evolving domains. This method leverages ghost-penalty stabilization in cutFEM to create a discrete extension of the velocity, ensuring that finite differences in time
\begin{equation*}
  \partial_t u \approx \frac{1}{\Delta t}(u^n - u^{n-1})
\end{equation*}
are well-defined, even when $u^n$ and $u^{n-1}$ are defined on different domains, $\Omega(t^n)$ and $\Omega(t^{n-1})$, respectively.

This fully Eulerian time-stepping scheme was introduced and analyzed in~\cite{LO19} for a convection-diffusion problem, applied to two-phase flow problems in~\cite{claus2019cutfem}, and used for flows with immersed moving bodies in~\cite{von2021falling}. A similar unfitted approach using space-time elements for flow problems was used in \cite{AB22}. Extensions for the Stokes and linearized Navier-Stokes equations in deforming domains were analyzed in~\cite{BFM22, vWRL21, NO23}, with optimal convergence in the energy norm for velocity and a scaled $L^2(H^1)$-norm for pressure established for unfitted Taylor--Hood finite elements of any degree.

Despite this solid mathematical foundation, we find here that the known variants of the method  are not free from spurious temporal oscillations in hydrodynamic forces. These oscillations, though small and rapidly decreasing with mesh refinement, persist even when volume and mass conservation are exact --- contrary to the common assumption that failures in conservation principles are the primary cause. Moreover, consistent with earlier studies, we observe that oscillation amplitude increases as the time step decreases, although the time-averaged amplitude remains bounded. This behavior suggests that the method's numerical stability bounds vary across different norms.

We argue that controlling spurious pressure oscillations, resulting in spurious hydrodynamic forces acting on immersed boundary, requires a discretization method that is (unconditionally) stable in the $L^\infty(H^1)$-norm for pressure. Standard finite element analyses typically assess pressure stability in the $L^2(L^2)$-norm, which only ensures average stability over time. Under plausible assumptions, $L^\infty(H^1)$-pressure stability can be reduced to $L^\infty(L^2)$-pressure stability. Both mean \emph{the stability of instantaneous pressure}.
Our numerical experiments show that all available unfitted finite element method variants exhibit unconditional stability in the $L^2(L^2)$-norm, yet the pressure blows up in the $L^\infty(L^2)$-norm as the time step approaches zero for a fixed spatial mesh. We revisit the analysis of finite element methods on steady domains to identify key arguments needed for $L^\infty(L^2)$-pressure stability and explain where these arguments fail in moving boundary cases.

Identifying where the arguments fail allows us to propose a variant of the method that ensures unconditional stability of the instantaneous pressure, which \emph{de facto} suppresses spurious pressure oscillations and eliminates any dependence on the specific CFL number. A key modification to the method is redefining the ghost-penalty term globally, rather than restricting it to a narrow band around the boundary, as previously suggested.

In addition to this finding, the paper examines the effects of fluid viscosity, geometry approximation, mass conservation, discretization and stabilization parameters, and the choice of finite element spaces on the occurrence and magnitude of spurious temporal pressure oscillations in the narrow-band immersed boundary FEM.
Our results indicate that some of these parameters have little to no impact on pressure temporal stability, whereas varying others may reduce spurious oscillations in the ‘traditional’ variant to levels below physical significance.

The remainder of this paper is organized as follows. In \Cref{sec.problem}, we briefly present the equations under consideration. Next, in \Cref{sec.method}, we introduce the unfitted Eulerian time-stepping method used to solve the problem numerically and summarize the key available results regarding this approach. Several variants of the stabilization are presented, including the new one with a global ghost-penalty term. 
\Cref{sec.math-persective} looks into the problem of spurious temporal oscillations of boundary forces from a mathematical perspective. Numerical results for our method are presented in \Cref{sec.numerics}, where we conduct a stability and convergence study of the method in various relevant norms. 
Additionally, we examine the method's performance in two challenging example scenarios from the literature, extended to the Navier-Stokes setting considered here. In a separate \Cref{sec.numerics:subsec.other-methods} we discuss numerical results obtained for several related unfitted finite element methods. 
A few concluding remarks are offered in \Cref{sec.conclusion}.

\section{Problem Description}
\label{sec.problem}

We consider the Navier-Stokes equations for incompressible viscous fluid
\begin{align}\label{NSE1}
  \partial_t \ub -\nu\Delta \ub + \ub\cdot\nabla\ub + \nabla p &= 0\\ \label{NSE2}
  \nabla\cdot\ub &= 0
\end{align}
in a time-dependent domain $\Omega(t)\subset\mathbb{R}^d$ which we assume to be embedded in a fixed and bounded domain $\widetilde\Omega$ for all times $t\in[0, T]$. In this study, we assume that the evolution of $\Omega(t)$ is given and smooth in the sense that there is a one-to-one continuous mapping $\Phi(t):\Omega(0)\to\Omega(t)$ such that $\Phi\in C^2(\overline{\Omega(0)}\times[0,T]; \mathbb{R}^d)$.

On the moving part of the boundary $\Gamma(t)$, we consider no-slip boundary conditions
\begin{equation*}
  \ub = \ub_\Gamma\coloneqq\partial_t \Phi.
\end{equation*}
For simplicity, on the fixed part of the boundary, we consider no-slip and free outflow boundary conditions.

In this set-up, the Reynolds number associated with the moving boundary is
\begin{equation*}
  Re = \frac{u_0 D}{\nu},
\end{equation*}
where $u_0$ is the maximal speed of the boundary and $D$ is characteristic size of $\Omega(t)$.

To quantify temporal spurious oscillations of the pressure, we will consider  various pressure norms and  the \emph{pressure drag coefficient} in the direction of motion of the moving boundary. This is defined as
\begin{equation}\label{eqn.pressure-drag}
  C_p(p) \coloneqq \frac{2}{\rho D^3 f_0^2}\int_{\Gamma(t)} p \nb \cdot \bm{e}_\text{mb}\dif s
\end{equation}
where $\bm{e}_\text{mb}$ is the direction of motion of the moving boundary $\Gamma(t)$ and $f_0$ is the frequency of the motion of the moving boundary.

\section{Unfitted Eulerian Finite Element Method}
\label{sec.method}

We consider an unfitted finite element method known as cutFEM~\cite{BCH14}. To tackle the moving domain nature of the problem at hand, we use an approach extending the discrete solution and enabling finite difference discretizations to the time derivatives as proposed in~\cite{LO19}. To achieve optimal high-order geometry approximations, we use the isoparametric approach for unfitted finite elements introduced in \cite{Leh16} and extended to moving domain problems in \cite{LL21}.

\subsection{Preliminaries}
\label{sec.method:subsec.prelim}

We focus on a cutFEM discretisation with inf-sup stable Taylor-Hood elements~\cite{BH14,MLLR14,GO17}, as used in the moving domain context in \cite{vWRL21,von2021falling}. While a number of other unfitted discretisations are possible and available in the moving domain context, we will focus on this element as the one of the most popular. Nevertheless, we also discuss some alternative finite element choices in \Cref{sec.numerics:subsec.other-methods} below. We begin by covering some of the notation necessary to define the method.

\subsubsection{Meshes}
Let $\Tht$ be a simplicial, shape-regular and quasi-uniform triangulation of the background domain $\widetilde{\Omega}$.
At each time-step, the physical domain $\Omega(t^n)$ is approximated by $\Ohn$. The moving part of the boundary of this domain is denoted as $\Ghn$. 
At time step $n$ the \emph{active mesh} includes all simplexes overlapping with $\Ohn$ with possibly a few additional layers of simplexes. More precisely, we define  
the active mesh as 
\begin{equation*}
  \Thnd \coloneqq \{K\in\Tht{}\mid \dist(\xb,\Ohn) \leq \delta_h, ~\text{for some}~\xb\in K\}
\end{equation*}
with 
\begin{equation*}
  \delta_h = 2\Delta t \uninf, 
\end{equation*}
where $\uninf=\max_t\|\mathbf{u}_\Gamma\|_{L^\infty(\Gamma(t))}$ is the maximal possible speed of the boundary motion. 
It is easy to see that this ensures that 
the active domains 
\begin{equation*}
  \Ohndt \coloneqq \{\xb\in K \mid K\in \Thnd\}.
\end{equation*} 
are such that the following crucial property holds:  
\begin{equation}\label{eq:emb}
\Omega_h^{n}\subset \widetilde{\Omega}_{h,\delta_h}^{n-2}\cap \widetilde{\Omega}_{h,\delta_h}^{n-1}.
\end{equation}

\begin{figure}
  \centering
  \includegraphics{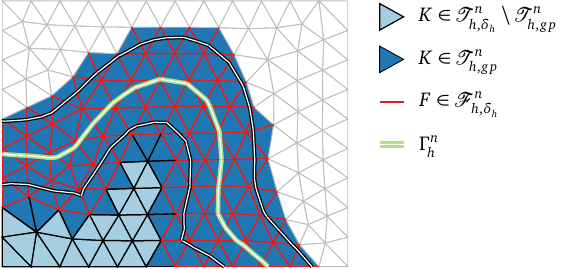}
  \caption{Sketch of different elements and facets used for the unfitted Eulerian time-stepping scheme.}
  \label{fig.mesh}
\end{figure}

For the stabilization and extension forms considered below, we define the elements in the $\delta_h$-strip around the moving interface
\begin{equation*}
  \ThnS \coloneqq \{K\in \Thnd \mid \dist(\xb, \Ghn) \leq \delta_h,~\text{for some}~\xb\in K \}
\end{equation*}
and the set of facets of these elements that are interior to the active mesh
\begin{equation}\label{eqn.extention.gp.facets}
  \Fhnd \coloneqq \{F = K_1\cap K_2 \mid K_1\in\ThnS,K_2\in\Thnd{}, K_1\neq K_2\}.
\end{equation}
Additionally, we define the set of elements intersected by the moving boundary
\begin{equation*}
  \ThnG \coloneqq \{K\in \Thnd \mid K\cap\Ghn\neq \emptyset \},
\end{equation*}
and the set of facets of these boundary elements that are interior to the cut mesh $\mathcal{T}_{h,0}^n$:
\begin{equation}\label{eqn.pre.gp.facets}
  \FhnG \coloneqq \{F = K_1\cap K_2 \mid K_1\in\ThnG,K_2\in \mathcal{T}_{h,0}^n, K_1\neq K_2\}.
\end{equation}
Finally, the set of active elements which have at least one facet in $\Fhnd$ are denoted by
\begin{equation*}
  \mathcal{T}_{h,gp}^{n}\coloneqq \{K\in \Thnd \mid \partial K \supset F \in \Fhnd \}.
\end{equation*}

\subsubsection{Weak forms}
For our discretization method (or variants based on other pairs of finite element spaces), we need the bilinear form for the viscosity term with weakly imposed Dirichlet boundary conditions
\begin{align*}
  a_h^n(\ub, \vb) &\coloneqq \nu \int_{\Ohn}\nabla \ub :\nabla\vb\dif\xb + \nu N_h^n(\ub, \vb)\\
  N_h^n(\ub, \vb) &\coloneqq N_{h,c}^n(\ub, \vb) + N_{h,c}^n(\vb, \ub) + N_{h,s}^n(\ub, \vb)
\end{align*}
and the Nitsche consistency and penalty terms given by
\begin{equation*}
  N_{h,c}^n(\ub, \vb) = \int_{\Ghn} -\dn\ub \cdot \vb \dif s
  \qquad\text{and}\qquad
  N_{h,s}^n(\ub, \vb) = \int_{\Ghn}\frac{\sigma k^2}{h}\ub\cdot\vb\dif s,
\end{equation*}
respectively. For the form to be coercive, the penalty parameter $\sigma>0$ needs to be chosen sufficiently large. The convective term is given in the convective form as
\begin{equation*}
  c_h^n(\ub, \vb, \wb) \coloneqq \int_{\Ohn} (\ub\cdot\nabla\vb)\cdot\wb\dif\xb.
\end{equation*}
Finally, for the velocity pressure coupling, we define
\begin{equation*}
  b_{h,0}^n(p, \vb) =  -\int_{\Ohn} p\diver(\vb)\dif\xb \qquad\text{and}\qquad
  b_h^n(p,\vb) = b_{h,0}^n(p, \vb) + \int_{\Ghn} p \vb\cdot\nb\dif s.
\end{equation*}

\subsection{Discrete Formulation}
\label{sec.method:subsec.taylor-hood}

This method is based on the work \cite{vWRL21} for the time-dependent Stokes problem in a moving domain.
On the active and cut meshes, we consider the inf-sup stable Taylor-Hood finite element pair
\begin{subequations}
\begin{align}
  V_h^n &\coloneqq \{v\in C(\Ohndt) \mid \restr{v}{K}\in\mathbb{P}^k(K)~\forall K\in \Thnd  \}^d,\\
  Q_h^n &\coloneqq \{q\in C(\Ohnt) \mid \restr{q}{K}\in\mathbb{P}^{k-1}(K)~\forall K\in \Thn  \},
\end{align}
\end{subequations}
with $k\geq 2$. For the inf-sup stability on unfitted meshes with ghost-penalty stabilization, we refer to \cite{GO17}.

The fully discrete variational formulation with a BDF2 time-discretization is given by
\begin{multline}\label{eqn.discrete.method:TH}
  \int_{\Ohn}\frac{3\uhn - 4\ub_h^{n-1} + \ub_h^{n-2}}{2 \Delta t} \vh \dif \xb + a_h^n(\uhn,\vh) + c_h^n(\uhn, \uhn, \vh) + b_h^n(p_h^n, \vh) + b_h^n(q_h,\uhn)\\
   + s_{gp}^{n}((\uhn, p_h^n), (\vh, q_h)) + s_{gd}^n(\uhn, \vh) = 0
\end{multline}
for all $(\vh, q_h)\in V_h^n\times Q_h^n$.
The finite difference in the first term is well-defined due to the domains  embedding \eqref{eq:emb}.
To improve the mass-conservation properties of the scheme, we consider the grad-div stabilization operator~\cite{olshanskii2002low}
\begin{equation}\label{eqn.grad-div-stabilization}
  s_{gd}^n(\ub,\vb) \coloneqq \ggd \int_{\Ohndt} \diver(\ub)\diver(\vb)\dif\xb,
\end{equation}
with the stabilization parameter $\ggd\geq 0$. We note that in \cite{vWRL21}, no grad-div stabilization was included.
Finally, the discrete extension enabling the Eulerian time-stepping scheme is realized by the ghost-penalty stabilization operator acting on the velocity and pressure
\begin{equation*}
  s_{gp}^{n}((\uhn, p_h^n), (\vh, q_h)) = \gpv L\big(\nu + \nu^{-1}\big)i_h^n(\ub,\vb) - \gpp(\nu + \ggd)^{-1} j_h^n(p,q)
\end{equation*}
where, $L=\left\lceil\frac{\delta_h}{h}\right\rceil$ is the strip width given in terms of the number of elements. This is needed to increase the stiffness in the ghost-penalty operator, when more facets must be crossed to reach an uncut element from an extension element, cf. \cite[Assumption~5.3 and Remark~5.4]{LO19}. Detailed choices for these operators will be given next.

\subsection{Ghost-penalty stabilization}
\label{sec.method:subsec.ghost-penalty}

As typical for unfitted finite elements on stationary domains, we have the so-called `small-cut' problem, where the stability of the solution and bounded condition numbers of the arising linear systems is lost when there are elements in the active mesh where only a small part  of the element is covered by the active domain. In cutFEM, ghost-penalty stabilization alleviates this problem. In the moving domain setting, the ghost-penalty operator has the additional role of realizing an extension of our discrete solution.

We consider the direct version of the ghost-penalty operator~\cite{Pre18,LO19}. Let $\pazocal{E}^\mathbb{P}\colon \mathbb{P}^k(K) \rightarrow \mathbb{P}^k(\mathbb{R}^d)$ be the canonical global extension operator of polynomials to $\mathbb{R}^d$, i.e., for $v\in\mathbb{P}^k(K)$ it holds that $\restr{\pazocal{E}^\mathbb{P}(v)}{K} = v$. For a facet $F=K_1\cap K_2$, $K_1\neq K_2$, we define the \emph{facet patch} as $\omega_F = K_1\cup K_2$. The volumetric patch jump $\jump{\cdot}_\omega$of a polynomial $\vh\in\mathbb{P}^{k}(\Tht)$ is then defined via
\begin{equation*}
  \restr{\jump{\vh}_\omega}{K_i} = \restr{\vh}{K_i} - \pazocal{E}^\mathbb{P}(\restr{\vh}{K_j}),\quad\text{for }i,j\in\{1,2\}\text{ and }i\neq j.
\end{equation*}
With this, we define the ghost-penalty operators
\begin{equation}\label{eqn.ghost-penalty.def}
  i_h^n(\ub, \vb) \coloneqq \sum_{F\in\Fhnd}\frac{1}{h^2}\int_{\omega_F}\jump{\ub}_\omega\jump{\vb}_\omega\dif \xb
  \qquad\text{and}\qquad
  j_h^n(p, q) \coloneqq \sum_{F\in\FhnG}\int_{\omega_F}\jump{p}_\omega\jump{q}_\omega\dif \xb,
\end{equation}
respectively.
We note that the different $h$-scalings are because the velocity extension needs control in the $H^1$-norm, while for the pressure extension the control in the $L^2$-norm is optimal; see, e.g.,~\cite{GO17}.

As discussed below, alternative choices for the set of ghost-penalty facets may be appropriate. To differentiate between the different choices, we call the above choice the \emph{narrow-band} stabilization choice, since $i_h^n(\cdot,\cdot)$ and $j_h^n(\cdot,\cdot)$ are acting only in a narrow band of triangles around $\Gamma_h^n$.

\subsubsection{Patch-wise choice of ghost-penalty facets}
\label{sec.method:subsec.ghost-penalty:subsubsec.patches}

The choice of ghost-penalty facets in \eqref{eqn.extention.gp.facets} and \eqref{eqn.pre.gp.facets} is sufficient for the stability and extension of the solution. However, it is not necessary. Indeed, a smaller, necessary set of facets can be chosen. The key to the ghost-penalty extension mechanism is the following assumption.

\begin{assumption}\label{assumption.ghost-penalty}
For every strip element $K\in\ThnS$ with $K\not\subset \Ohn$, there exists an uncut element $K'\in\ThnS$ with $K'\subset\Ohn$, which can be reached by a path that crosses a bounded number of facets $F\in\Fhnd$. The number of facets on each path can be bounded by $L \lesssim 1 + \frac{\delta_h}{h}$ and each uncut element $K'$ is at the end of at most $M$ such paths, with $M$ bounded independent of $h$ and $\Delta t$. 
\end{assumption}

This assumption essentially states, that each interior element in the strip domain supports at most $M$ elements that require an extension.

An alternative construction for $\Fhnd$, is to cluster elements in $\ThnS$ together into agglomerated elements, such that the area inside and outside of the physical domain is approximately equal, and then to only consider the facets interior to every patch. This is then a subset of the above choice $\Fhnd$ and \Cref{assumption.ghost-penalty} remains valid. This choice of stabilizing facets is known as \emph{weak AggFEM}~\cite{BNV22}, since the idea of AggFEM~\cite{BVM18} is to consider functions that are polynomials within each patch. 

In \cite{BNV22}, it is argued that the patch-wise choice for the ghost-penalty operator is preferable, since it allows to avoid
certain locking effect for large stabilization parameters. Below, we will therefore investigate if this the choice of ghost-penalty facets has any impact on oscillations of the pressure in this moving domain context.

\subsubsection{Global choice of ghost-penalty facets}
We will see in the mathematical considerations in \Cref{sec.math-persective} below, one key aspect preventing us from extending the stability analysis of the pressure drag coefficient from the fixed domain to the moving domain setting is the fact that the finite element spaces and divergence condition vary between time steps. One possibility to allow us to use the same discrete space $V_h^n\equiv V_h$ and make divergence constraint more uniform in time is to replace $\Fhnd$ and $\FhnG$ in \eqref{eqn.ghost-penalty.def} with
\begin{equation}\label{eqn.all-facets}
  \Fht \coloneqq \{F=K_1\cap K_2 \mid K_1,K_2\in\Tht, \meas_{d-1}(F)>0\},
\end{equation}
i.e., apply ghost-penalty stabilization on the entire active mesh. This introduces a consistency error, however, we have
\begin{equation*}
  i_h(w,w) \lesssim h^{2m}\norm{w}_{H^{m+1}(\widetilde{\Omega})}^2
\end{equation*}
for all $w\in H^{m+1}(\widetilde{\Omega})$, cf.\ \cite[Lemma 5.8]{LO19}. Consequently, the consistency error decreases with the appropriate order such that final error estimates are not affected. While we do not have any numerical analysis for this choice, we note that in the context of unfitted space-time methods, this choice was successfully applied in \cite{AB22}. In this setting $\delta_h$ does not play a role, and consequently, we set $L=1$. Below, we will call the resulting method the \emph{global} ghost-penalty stabilization choice.

\subsection{Geometry Approximation}
\label{sec.method:subsec.geom.approx}

Separating the geometry from the computational mesh introduces the challenge of accurately approximating the geometry, requiring robust numerical integration techniques on cut elements. Several techniques have been proposed in the literature; see, e.g.,~\cite{Leh16,MKO13,Say15,OS16,FOSS17,HSK17}.

In cutFEM, the standard approach is to use a piecewise linear interpolation of the level set function. On simplices, this allows for the creation of quadrature rules by dividing cut elements into sub-simplices and applying standard Gau\ss quadrature~\cite{BCH14}. However, this piecewise linear approximation introduces a geometry error:
\begin{equation*}
  \dist(\Omega, \Omega_h) \lesssim h^2.
\end{equation*}
For Taylor–Hood elements with fine enough meshes, this error dominates the $L^2$ velocity error for $k=2$, and for higher-order elements, it affects $L^2$ and $H^1$ velocity errors as well as the $L^2$ pressure error.

The isoparametric unfitted finite element method~\cite{Leh16}, analyzed in the moving domain context~\cite{LL21}, reduces the geometry error to $\mathcal{O}(h^{k+1})$. This approach computes a small mesh deformation $\Theta_h$ to map the piecewise linear level set onto a higher-order approximation of the exact level set. In the moving domain context, the challenge lies in projecting the velocity $\ub_h^{n-1}$, computed with respect to $\Theta_h^{n-1}$, to the deformation $\Theta_h^n$ in \eqref{eqn.discrete.method:TH} to preserve optimal convergence~\cite{LL21}.

An alternative, simpler approach is to subdivide cut elements for quadrature and apply the piecewise linear level set approximation on these subdivisions. While this reduces geometry error, it becomes computationally infeasible as finer mesh refinements require more subdivisions to maintain optimal convergence. Nonetheless, it can be useful for isolating pure geometry approximation error, distinct from the projections required in the isoparametric approach.

\subsection{Properties of the Method}
\label{sec.method:subsec.analysis}

We summarize the key relevant results known for the presented approaches. 
The Taylor--Hood version and an equal order variant of the unfitted FEMs for (linearized) incompressible fluid problems in evolving domains were recently analyzed in \cite{vWRL21,NO23} and \cite{BFM22}, respectively. 

The best available error estimate for the  Taylor--Hood method is given in \cite[Theorem~4.8]{NO23} for the BDF1 counterpart of \eqref{eqn.discrete.method:TH} with a linearized $c$-term (these simplifications are not essential for our discussion). We review this result here:  Assume that the mesh $\Tht$ is quasi-uniform and the mesh size and time step satisfy 
\begin{equation}\label{cond1}
c_0 h^2\le \Delta t \le C_0 h,    
\end{equation}
with some positive $c_0, C_0$ independent of the discretization parameters and domain motion. Then a Taylor--Hood $P^k$--$P^{k-1}$ unfitted FEM satisfies the error bound 
\begin{equation} \label{est1}
    \max_{1\le n\le N}\|\bu(t^n)-\bu_h^n\|^2_{L^2(\Ohn)} +\Delta t \sum_{n=1}^N \|\nabla(\bu(t^n)-\bu_h^n)\|^2_{L^2(\Ohn)} +
    \Delta t \sum_{n=1}^N h^2\|\nabla(p(t^n)-p_h^n)\|^2_{L^2(\Ohn)}
    \le C(|\Delta t|+ h^{q} +h^{k})^2,
\end{equation}
where $C$ may depend on $u$ and $p$, solutions to \eqref{NSE1}--\eqref{NSE2},  but not on the discretization parameters and domain motion through the background mesh; $q$ is the order of the geometry approximation.

For equal order $P^k$--$P^k$ elements, an estimate similar to \eqref{est1} under assumption~\eqref{cond1} was proved in \cite[Theorem~5.3]{BFM22} for the Stokes problem in an evolving domain. The paper also provided the (sub-optimal) estimate for the $L^2$-norm of the pressure error, 
\begin{equation} \label{est2}
    |\Delta t|^2 \sum_{n=1}^N \|p(t^n)-p_h^n\|^2_{L^2(\Ohn)}
    \le C(|\Delta t| +h^{m})^2.
\end{equation}
There is no $h^{q}$ term in \eqref{est2} since the geometric error was not quantified in this study. 

From \eqref{cond1}, we observe that the restriction preventing the time step from approaching zero also arises in the available analysis. However, even with this restriction, the pressure error bounds in \eqref{est1}--\eqref{est2} are not optimal as $h, \Delta t \to 0$. The authors of the three referenced papers \cite{vWRL21, NO23, BFM22} have identified a key obstacle to achieving optimal-order error analysis: the loss of the weakly divergence-free property of $\bu_h^{n-1}$ at time step $n$.  

This challenge motivated us to introduce grad-div stabilization \eqref{eqn.grad-div-stabilization} across the entire active domain to improve the divergence of the velocity. Additionally, we considered using strongly divergence-free elements; see \Cref{sec.numerics:subsec.other-methods}. Neither of these modifications resulted in a significant reduction in spurious pressure oscillations, prompting us to seek alternative analytical insights in the next section.

\section{Spurious pressure force oscillations: mathematical perspective}
\label{sec.math-persective}

The trace theorem provides us with the bound of the boundary pressure force by a volumetric norm of the pressure: 
\begin{equation}\label{p_bound}
  |F_p(t)| \coloneqq \Big|\int_{\Gamma(t)} p \bm{n} \,\dif s\Big|\le \|p\|_{L^1{(\Gamma(t)})}\le C\,\|p\|_{W^{1,1}(\Omega(t))},
\end{equation}
where $C$ depends only on $\Omega(t)$ and for a smoothly deforming  domain $C$ is uniformly bounded for $t\in[0,T]$. 
Therefore, control of spurious oscillations in $F_p(t)$ (and so in $C_p(p)$) follows from the unconditional stability of the finite element pressure in the  (discrete analogue of)  $L^\infty(0,T; W^{1,1}(\Omega(t)))$-norm: 
\begin{equation}\label{p_W11}
  \max_{n} \Big( \|p^n_h\|_{L^1(\Omega(t^n))}+\|\nabla p^n_h\|_{L^1(\Omega(t^n)))}\Big)\le C, 
\end{equation}
with some $C$ independent of the discretization parameters. 
Such uniform in time bounds for $p_h^n$ mean the stability of the \emph{instantaneous discrete pressure} and instantaneous pressure gradient.  

Since the Lebesgue $L^1$-norm is less convenient for analysis then $L^2$-norm and $L^2(\Omega)\subset L^1(\Omega)$ for bounded domains, one can look instead for the pressure stability in the $L^\infty(H^1)$-norm. Furthermore, by the finite element inverse inequality and the extension property of the form $j^n_h(\cdot,\cdot)$, it holds that
\begin{equation*}
  \max_n \|p_h^n\|_{H^1(\Omega(t^n))}\le C h^{-1} \max_n \left( \|p_h^n\|_{L^2(\Omega(t^n))} + s_h^{\frac12}(p_h^n,p_h^n)\right).
\end{equation*}
Therefore, as far as we are concern with  stability of the instantaneous pressure being uniform in $\Delta t$ (or the specific CFL number), the condition \eqref{p_W11} can be weaken to the pressure stability in $L^\infty(L^2)$-norm: 
\begin{equation}\label{p_L2}
 \max_n \Big(\|p_h^n\|_{L^2(\Omega(t^n))}^2 + s_h(p_h^n,p_h^n)\Big)\le C,
\end{equation}
with $C$ independent of mesh parameters.

We conclude that pursuing an oscillation-free discretization method can be viewed as developing an instantaneous pressure-stable method in the sense of \eqref{p_W11} or \eqref{p_L2}. With this more mathematical perspective in mind, let us consider an argument that demonstrates how a standard (fitted) finite element (FE) method in a steady domain satisfies \eqref{p_L2}. Understanding where this argument fails for the unfitted FEM in an evolving domain provides further insight into the cause and potential remedies for spurious boundary force oscillations.

We recall the inf-sup (also known as LBB) stability condition satisfied by the fitted Taylor--Hood element:
\begin{equation*}
    \beta\|p_h^n\|_{L^2(\Omega)} \le \sup_{\bv_h\in V_h}\frac{ (p_h^n,\nabla\cdot\bv_h)}{\|\nabla \bv_h\|_{L^2(\Omega)}},
\end{equation*}
with some $\beta>0$ independent of $h$.

To avoid unnecessary technical details, we shall consider the unsteady Stokes system with $\nu=1$, $T=1$, $\Delta t=1/N$, $\ub|_{t=0}=0$ in $\Omega$, $\ub=0$ on $\Gamma$, and a smooth body force $\blf$.
Given a velocity--pressure pair of FE spaces $V_h\times Q_h$ satisfying the inf-sup stability condition, the 
first order implicit method takes the form
\begin{equation} \label{eq1}
    (D_t^n\bu_h, \bv_h)+ (\nabla\bu^n_h,\nabla\bv_h) - (p_h^n,\nabla\cdot\bv_h)+ (q_h^n,\nabla\cdot\bu_h) = (\blf^n,\bv_h)
\end{equation}
for all $\bv_h\in V_h$, $q_h\in Q_h$, $n=1,\dots,N$.  Here and further we use the notation $(\cdot,\cdot)$ for $L^2$ scalar product and
$D_t^n\phi=(\phi^n-\phi^{n-1})/\Delta t$ for the BDF1 derivative of a quantity $\phi=\{\phi^0, \phi^1,\dots\}$.

Letting $q_h=0$, $\bv_h=D_t^n\bu_h$ in \eqref{eq1}, and noting 
\begin{equation} \label{eq2}
    (p_h^n,\nabla\cdot D_t^n\bu_h) =0,
\end{equation}
leads to the identity and the inequality
\begin{equation}\label{aux552} 
    \|D_t^n\bu_h\|^2_{L^2(\Omega)}+ \frac12\left(D^n_t\|\nabla\bu_h\|^2_{L^2(\Omega)}+\Delta t \|\nabla D^n_t\bu\|^2_{L^2(\Omega)}\right) = (\blf^n,D_t^n\bu_h) \le \frac12 \left(\|\blf^n_h\|^2_{L^2(\Omega)}+\|D^n_t\bu^n\|^2_{L^2(\Omega)}\right).
\end{equation}
Multiplying the above inequality by $\Delta t$ and summing up for $n=1,\dots K$, with any $K\le N$, gives after straightforward calculations the following stability bound:
\begin{equation}\label{stab1} 
   \sum_{n=1}^N\Delta t \|D_t^n\bu_h\|^2_{L^2(\Omega)}+ \max_{1\le n\le N}\|\nabla\bu^n_h\|^2_{L^2(\Omega)} \le  \sum_{n=1}^N\Delta t\|\blf^n_h\|^2_{L^2(\Omega)}=:C_1,
\end{equation}
where we also used $\|\nabla\bu^0_h\|_{L^2(\Omega)}=0$.   

The estimate \eqref{stab1} provides $L^\infty$ control over time for the gradient of the finite element velocity but only $L^2$ control for the time derivative, which is insufficient for our purposes. Therefore, we proceed as follows: Let $\bd^n_h = D_t^n \bu_h$ and $r_h^n = D_t^n p_h$. By subtracting \eqref{eq1} from itself with an index shift, we find that $\bd^n_h$ and $r_h^n$ satisfy
\begin{equation} \label{eq1d}
    (D_t^n\bd_h, \bv_h)+ (\nabla\bd^n_h,\nabla\bv_h) - (r_h^n,\nabla\cdot\bv_h)+ (q_h^n,\nabla\cdot\bd^n_h) = (D_t^n\blf,\bv_h)
\end{equation}
for all $\bv_h\in V_h$, $q_h\in Q_h$, $n=2,\dots,N$. Testing with $q_h=r_h$, $\bv_h=\bd_h$ and using similar arguments as above, we arrive at the estimate
\begin{equation}\label{stab2} 
  \max_{1\le n\le N}\|\bd^n_h\|^2_{L^2(\Omega)} \le  c_f\sum_{n=1}^N\Delta t\|D_t^n\blf\|^2_{L^2(\Omega)}+\|\bd^1_h\|^2_{L^2(\Omega)},
\end{equation}
where $c_f$ is the Poincare--Friedrichs constant for $\Omega$. From \eqref{aux552} for $n=1$, we have the bound $\|\bd^1_h\|_{L^2(\Omega)}\le \|\blf^1_h\|_{L^2(\Omega)}$ (one uses $\|\nabla\bu^0_h\|_{L^2(\Omega)}=0$ to show it). 
Thus, the smoothness of $\blf$ and \eqref{stab2} provides the control of the $L^\infty$-norm of the velocity time derivative:
\begin{equation}\label{stab3} 
  \max_{1\le n\le N}\|D_t^n\bu_h\|_{L^2(\Omega)} \le  C_2,
\end{equation}
with a constant $C_2$ independent of $h$ and $\Delta t$.  

The pressure inf-sup stability condition and \eqref{eq1} implies
\begin{equation}\label{eq:aux615}
\begin{split}
    \beta\|p_h^n\|_{L^2(\Omega)} &\le \sup_{\bv_h\in V_h}\frac{ (p_h^n,\nabla\cdot\bv_h)}{\|\nabla \bv_h\|_{L^2(\Omega)}}
    =\sup_{\bv_h\in V_h}\frac{(D_t^n\bu_h, \bv_h)+ (\nabla\bu^n_h,\nabla\bv_h) - (\blf^n,\bv_h)}{\|\nabla \bv_h\|_{L^2(\Omega)}}\\
    &\le c_f\|D_t^n\bu_h\|_{L^2(\Omega)}+\|\nabla\bu^n_h\|_{L^2(\Omega)}+c_f\|\blf^n\|_{L^2(\Omega)},
\end{split}
\end{equation}
with the inf-sup constant $\beta>0$ independent of $h$ and $\Delta t$. 
Applying \eqref{stab1} and \eqref{stab3} to the right-hand side of the above inequality we get the $L^\infty$ stability bound for the pressure:
\begin{equation}\label{stab4} 
  \max_{1\le n\le N}\|p_h^n\|_{L^2(\Omega)} \le  C.
\end{equation}

The upper bound in \eqref{stab4} is the proper analogue of \eqref{p_L2} and  it provides an unconditional stability estimate for the instantaneous discrete pressure. As we already mentioned, by the finite element inverse inequality, \eqref{stab4} is also sufficient to control the instantaneous pressure gradient uniformly in time, though with an $h$-dependent upper bound:
$\max_{1\le n\le N}\|\nabla p_h^n\|_{L^2(\Omega)}\le  C h^{-1}$. Achieving control of the instantaneous pressure gradient that is robust in both $\Delta t$ and $h$ is also possible but would require additional technical details, which are beyond the scope of this paper.
\smallskip

In examining the proof of \eqref{stab4}, we identify two key arguments that do not extend to the case of an evolving domain $\Omega(t)$ and the unfitted Eulerian FEM, such as in \eqref{eqn.discrete.method:TH}: 
\begin{description}
    \item[(i)] The  divergence conditions for $\bu_h^n$ and $\bu_h^{n-1}$ differ, meaning that $D_t^n\bu$ does not satisfy an analogue of \eqref{eq2}.

    Indeed, even without pressure-stabilizing terms, the condition $b_h^{n-1}(q_h,\bu_h^{n-1})=0$ for all $q_h \in Q_h^{n-1}$ does not necessarily imply $b_h^{n}(q_h,\bu_h^{n-1})=0$ for all $q_h \in Q_h^{n}$. As a result, the discrete time-derivative $(\bu_h^{n}-\bu_h^{n-1})/\Delta t$ is not weakly divergence-free, which in turn may introduce a spurious potential in the finite element momentum equation. This spurious potential is then scaled by $|\Delta t|^{-1}$, leading to degradation in the pressure variable for smaller time steps.
    
    \item[(ii)]  More generally, the spaces $V_h^n$ and $V_h^{n-1}$ are not necessarily the same, preventing us from deriving the equation \eqref{eq1d}. 
\end{description}

We address the above issues by introducing the global ghost penalty stabilization in \Cref{sec.method:subsec.ghost-penalty:subsubsec.patches}. This ensures that the finite element spaces remain the same across time steps, i.e., $V_h^n = V_h^{n-1}$, and, perhaps more importantly, sets the stabilization forms $j^n_h(\cdot,\cdot)$ in the  divergence condition and
$i^n_h(\cdot,\cdot)$ in the momentum equation to be the same for all $n$.

\section{Numerical Results}
\label{sec.numerics}

All numerical examples in this section are implemented using \texttt{netgen}/\texttt{NGSolve}~\cite{Sch97,Sch14} and the add-on \texttt{ngsxfem}~\cite{LHPvW21} for unfitted finite element discretizations.

In addition to standard statistics like the error between the true and the discrete solutions, we will be interested in quantifying spurious boundary pressure force oscillations.

\subsection{Measuring pressure oscillations}
\label{sec.num:subsec.pre-oscil}

To measure the presence and scale of any non-physical high-frequency oscillations present in the pressure drag coefficient \eqref{eqn.pressure-drag}, we filter out a smooth component of the pressure drag coefficient history and compare this with the discrete data
\begin{equation*}
  e_h^n = \overline{C}_p^n - C_{p, h}^n,
\end{equation*}
where $C_{p, h}^n \coloneqq C_p(p_h^n)$ as defined in \eqref{eqn.pressure-drag}, and $\overline{C}_p^n$ is the post-processed signal at time $t^n$ computed through a second order Butterworth digital high-pass filter. We then consider both the $L^\infty(\mathbb{T})$-norm and the $L^1(\mathbb{T})$-norm of the signal error to measure the maximal amplitude and the amount of noise, respectively. We note that $\mathbb{T}\subset [0,T]$ is not the complete computational time interval to avoid edge effects inherent in signal processing due to our noisy signal not being periodic.

Oscillations in the pressure commonly depend on the ratio between the time-step and the mesh size. To indicate this, we consider the specific CFL number from  \eqref{eq:CFL} with $V_\Gamma=\Vert \ub_\Gamma\Vert_\infty$
as suggested in \cite{XCM24}. 

\subsection{Example 1: Convergence Study}
\label{sec.numerics:subsec.convergence}

To test the general convergence properties of the method, we consider a test case with a synthetic  analytical solution as described in \cite{vWRL21}. 
The domain is $\Omega(t)=\{\xb\in\mathbb{R}^2\;|\; (\xb_1-t)^2 + \xb_2^2 < 1/2\}$, and the analytical solution is given by
\begin{equation*}
  \ub_{ex}(t) =
  \begin{pmatrix}
    2\pi \xb_2\cos(\pi((\xb_1-t)^2 +  \xb_2^2))\\
    -2\pi \xb_1\cos(\pi((\xb_1-t)^2 + \xb_2^2))
  \end{pmatrix}
  \quad\text{and}\quad 
  p_{ex}(t)=\sin(\pi((\xb_1-t)^2 + \xb_2^2)) - \frac{2}{\pi}.
\end{equation*}
The error between the numerical velocity and pressure and the analytical solution are then computed using the discrete norms
\begin{equation*}
  \Vert e_h \Vert_{L^2(X)} \coloneqq \Big(\Delta t\sum_{i=1}^n\Vert e_h^i \Vert_X^2 \Big)^{1/2}
  \quad\text{and}\quad
  \Vert e_h \Vert_{L^\infty(X)} \coloneqq \max_{i=1,\dots,n}\Vert e_h^i \Vert_X.
\end{equation*}

We consider the background domain $\widetilde{\Omega} = (-1,2) \times (-1,1)$ with a time interval of $[0, 1]$ and $\nu=0.01$. We choose a uniform triangulation of $\widetilde{\Omega}$ with $h = 0.2$ for the initial mesh. The initial time step is $\Delta t = 0.1$. The mesh and time step are both refined uniformly and simultaneously five times for the convergence study. Consequently, the width of the strip $\ThnS$ (in terms of the number of elements) is $L = 1$ in this case.

\subsubsection*{Results}

\begin{figure}
  \centering
  \includegraphics{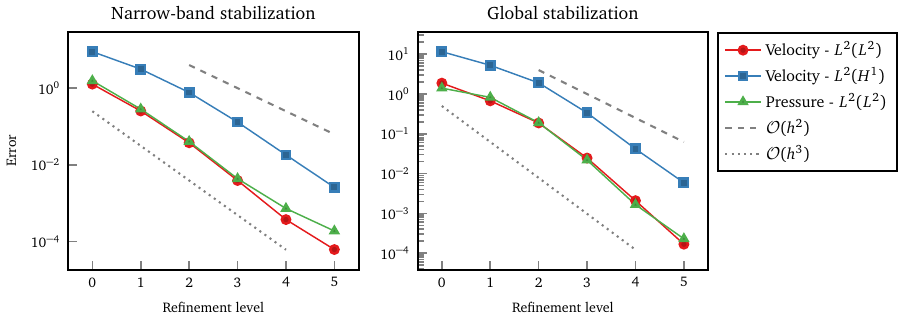}
  \caption{\hyperref[sec.numerics:subsec.convergence]{Example 1} -- Convergence behavior of the error of the velocity and pressure under uniform mesh refinement and combined mesh and time-step refinement for the isoparametric Taylor-Hood elements with the sufficient and global ghost-penalty facet choices.}
  \label{fig.convergence}
\end{figure}

The convergence results for the isoparametric unfitted Taylor-Hood method with BDF2 time-stepping and the narrow band and global ghost-penalty stabilization choices are presented in \Cref{fig.convergence} for the velocity error in the $L^2(L^2)$- and $L^2(H^1)$-norms, respectively, and in the $L^2(L^2)$-norm for the pressure error. After some pre-asymptotic behavior in the case of the global stabilization choice, we see optimal order convergence with both variants of stabilization. The results also suggest  that the spatial error is dominant in this set up and that the global ghost-penalty choice introduces slightly larger consistency error.

\begin{figure}
  \centering
  \includegraphics{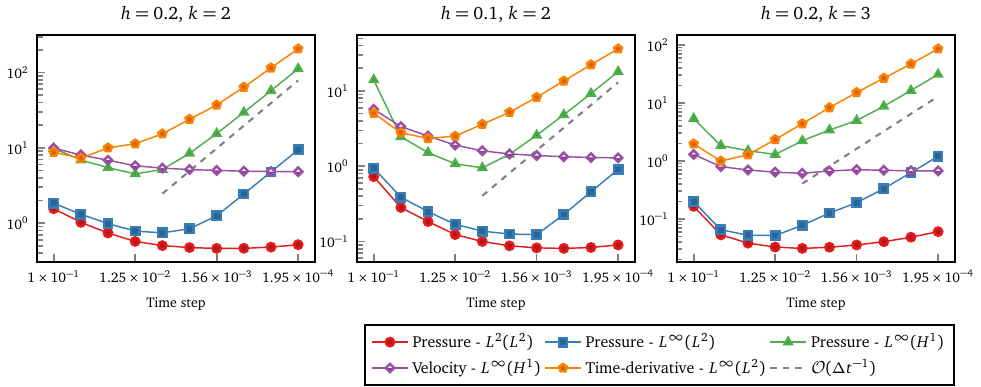}
  \caption{\hyperref[sec.numerics:subsec.convergence]{Example 1} -- Behavior of pressure, velocity, and discrete velocity time derivative norms under uniform time-step refinement for the isoparametric Taylor-Hood elements with the narrow-band ghost-penalty stabilization.  }
  \label{fig.pressure-norm.time-step.refinement.sufficient}
\end{figure}

\begin{figure}
  \centering
  \includegraphics{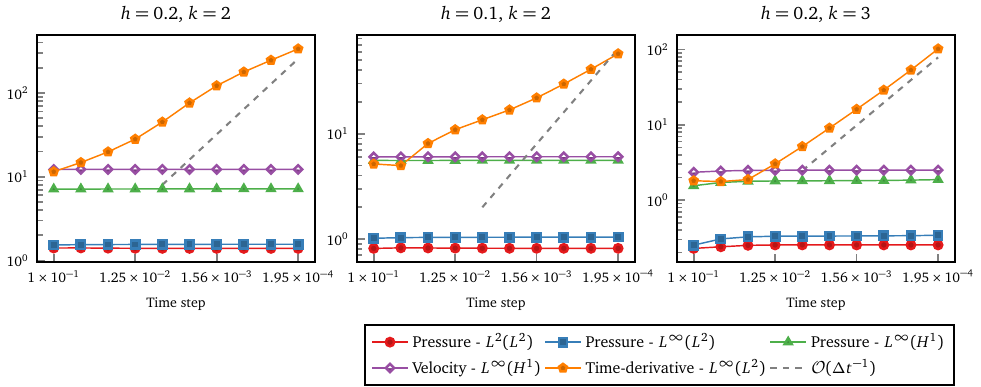}
  \caption{\hyperref[sec.numerics:subsec.convergence]{Example 1} -- Behavior of the pressure, velocity, and discrete velocity time derivative under uniform time-step refinement for the isoparametric Taylor-Hood elements with the global ghost-penalty stabilization.}
  \label{fig.pressure-norm.time-step.refinement.global}
\end{figure}

In line with analysis in \Cref{sec.math-persective}, we investigate the stability of the velocity and pressure with respect to the $L^\infty$-norm in time over a series of ten time-steps refinement. The results for the narrow-band stabilization are presented in \Cref{fig.pressure-norm.time-step.refinement.sufficient} and for the global stabilization in \Cref{fig.pressure-norm.time-step.refinement.global}.

For the narrow-band stabilization, we see that the pressure diverges in the $L^\infty(L^2)$- and $L^\infty(H^1)$-norm at a rate of $O(\Delta t^{-1})$ for small time steps, while the $L^2(L^2)$-norm error remains bounded. This indicates that the method is not unconditionally stable in either the $L^\infty(H^1)$-norm or the weaker $L^\infty(L^2)$-norm for pressure, although it is stable in the $L^2(L^2)$ pressure norm.

In the middle plot of \Cref{fig.pressure-norm.time-step.refinement.sufficient}, we observe that refining the spatial mesh once delays the increase in the $L^\infty(L^2)$- and $L^\infty(H^1)$-norm pressure errors by one level of time-step refinement. In the right plot of \Cref{fig.pressure-norm.time-step.refinement.sufficient}, we see that increasing the polynomial degree may improve accuracy but does not significantly influence the behavior with respect to the time step.

The key estimate \eqref{eq:aux615} for the instantaneous pressure involves the $L^\infty(H^1)$-norm of the velocity and the $L^\infty(L^2)$-norm of the velocity time derivative. The results in \Cref{fig.pressure-norm.time-step.refinement.sufficient} indicate that, of these two terms, it is the velocity time derivative that lacks a uniform stability bound in the $L^\infty(L^2)$-norm, whereas the velocity $L^\infty(H^1)$-norm remains uniformly bounded.

Examining the results with global stabilization, we find that both the velocity and pressure are uniformly bounded in both the $L^\infty(L^2)$-norm and $L^\infty(H^1)$-norm across all three spatial discretization scenarios previously investigated for the narrow-band stabilization case. This suggests that the global stabilization approach is a promising candidate for eliminating spurious temporal pressure oscillations.
Interestingly, the velocity time derivative still lacks a uniform stability bound, which implies that \eqref{eq:aux615} may not be sharp.

\subsection{Example 2: Moving cylinder in a cross flow}
\label{sec.numerics:subsec.ex2}

We now investigate the behavior of the pressure drag coefficient in more detail. This example, adapted from \cite{ruberg2014fixed,XCM24}, involves an oscillating cylinder in a cross-flow. The background domain is $\widetilde\Omega=(-10,10)\times(-5,5)$. The level set of the moving cylinder is given by
\begin{equation}\label{eqn.ex2.lset}
  \phi(t) = 1 - \sqrt{(\xb_1+5)^2 + (\xb_2 - x_c(t))^2} \quad\text{with}\quad x_c(t)\coloneqq x_0\cos(2\pi f_0 t),
\end{equation}
where $x_0=0.5$ is the initial vertical position of the cylinder, $f_0=0.2237$ is the oscillation frequency and the cylinder velocity is given by
\begin{equation}\label{eqn.ex2.bc}
  \ub_\Gamma(t) = (0, -u_0\sin(2\pi f_0 t))^T, \quad\text{with the maximal speed}\quad u_0=2\pi f_0 x_0.
\end{equation}
On the moving (unfitted) boundary, we impose the no-slip condition $\ub(t) = \ub_\Gamma(t)$. At the left (fitted) boundary of the domain, a inflow Dirichlet boundary condition $\ub = (1, 0)^T$ is applied, while the top and bottom (fitted) boundaries have a free-slip condition. At the right (fitted) boundary, we apply a free outflow `do-nothing' condition. The initial condition is given by a stationary Stokes flow in the domain at $t = 0$.

The time interval considered is $[0, 10]$. To avoid edge effects in signal processing, we analyze the results in the interval $\mathbb{T} = [1,9]$ to quantify spurious oscillations; see \Cref{sec.num:subsec.pre-oscil}.

Each computation is performed with three different time steps, resulting in CFL numbers of $1$, $0.1$, and $0.01$. Given that the method also involves numerous parameters likely to impact performance, we examine the effects of viscosity, ghost-penalty stabilization parameter, grad-div stabilization, and geometry approximation.

\subsubsection*{Results}

The results for the oscillatory component of the pressure drag force are presented in \Cref{tab.ex2.th.viscosity,tab.ex2.th.gamma-gd,tab.ex2.th.gp-parameters,tab.ex2.th.gp+gd,tab.ex2.th.geometry,tab.ex2.stabil+ho} for the narrow band stabilisation choice and \Cref{tab.ex2.th.global-stabil} for the global ghost-penalty stabilization. To gain insight into both the amplitude and frequency of potential spurious oscillations, we use the $L^\infty$- and $L^1$-norms. Based on these results, we make the following observations:

\paragraph{Viscosity dependence}
\Cref{tab.ex2.th.viscosity} presents results for viscosities $\nu \in \{1, 0.1, 0.01\}$ with a fixed spatial discretization and no parameter tuning (we shall call the choice of  $\gpv=\gpp=0.1$, $\ggd=0$ \emph{default}). We observe that the magnitude and frequency of the pressure drag oscillations increase as viscosity decreases. Additionally, for higher viscosity (essentially a Stokes flow), the magnitude of oscillations only increases by a factor of 1.7-2.4, with a 100-fold decrease in the CFL number. Even with the default set of parameters, the $\nu=1$ results are a significant improvement over the state-of-the art results achieved with the weighted shifted boundary method for the same setup  of a Stokes flow, cf.~\cite[Fig.~27]{XCM24}. Finally, the largest oscillation, $\Vert{e_h}\Vert_{L^\infty} = 1.17$, for the case of $C_\text{CFL} = 0.01$ and $\nu = 0.01$, is about 5\% relative to the signal amplitude.

\begin{table}
  \centering
  \caption{\hyperref[sec.numerics:subsec.ex2]{Example 2} -- Norms of pressure drag oscillations for varying viscosities. Isoparametric Taylor-Hood elements with $h=0.1$, $k=2$, $\gpv=\gpp=0.1$, $\ggd=0$ and narrow-band ghost-penalty stabilization.}
  \label{tab.ex2.th.viscosity}
  \begin{tabular}{ccc@{\hspace*{8pt}}ccc@{\hspace*{8pt}}ccc@{\hspace*{8pt}}c}
    \toprule
      && \multicolumn{2}{c}{$\nu=1$} && \multicolumn{2}{c}{$\nu=0.1$} && \multicolumn{2}{c}{$\nu=0.01$}\\
    \cmidrule{3-10}
    $\Delta t$ && $\Vert e_h\Vert_{L^\infty(\mathbb{T})}$ & $\Vert e_h\Vert_{L^1(\mathbb{T})}$ && $\Vert e_h\Vert_{L^\infty(\mathbb{T})}$ & $\Vert e_h\Vert_{L^1(\mathbb{T})}$ && $\Vert e_h\Vert_{L^\infty(\mathbb{T})}$ & $\Vert e_h\Vert_{L^1(\mathbb{T})}$\\
    \midrule
    0.1400 && 0.01015 & 0.00349 && 0.02809 & 0.01026 && 0.17018 & 0.04400\\
    0.0140 && 0.01486 & 0.00503 && 0.03028 & 0.00771 && 0.27549 & 0.06265\\
    0.0014 && 0.01785 & 0.00481 && 0.06778 & 0.00609 && 1.17637 & 0.06584\\
    \bottomrule
  \end{tabular}
\end{table}
  
\paragraph{Grad-div stabilization}
\Cref{tab.ex2.th.gamma-gd} examines the effects of grad-div stabilization for $\gamma_{\text{gd}} \in \{0.01, 0.1, 1.0\}$ with $\nu = 0.01$. Comparing these results to those in \Cref{tab.ex2.th.viscosity}, we observe that penalizing the divergence of the velocity field reduces the pressure drag oscillations in both the $L^\infty$ and $L^1$-norms. For $C_{\text{CFL}} = 0.01$, the largest reduction occurs with $\gamma_{\text{gd}} = 0.1$. Quantitatively, for $C_{\text{CFL}} = 0.01$, the best decrease is approximately a factor of 2.9, whereas for $C_{\text{CFL}} = 0.1$, the improvement factor is only 1.3. The observation that increasing the grad-div  penalty parameter beyond this does not further reduce oscillations suggests that divergence conformity alone does not fully address the spurious temporal oscillations in pressure.

\begin{table}
  \centering
  \caption{\hyperref[sec.numerics:subsec.ex2]{Example 2} --  Norms of pressure drag oscillations for varying grad-div stabilization. Isoparametric Taylor-Hood elements with $h=0.1$, $k=2$, $\nu=0.01$, $\gpv=\gpp=0.1$ and narrow-band ghost-penalty stabilization.}
  \label{tab.ex2.th.gamma-gd}
  \begin{tabular}{ccc@{\hspace*{8pt}}ccc@{\hspace*{8pt}}ccc@{\hspace*{8pt}}c}
    \toprule
      && \multicolumn{2}{c}{$\ggd=0.01$} && \multicolumn{2}{c}{$\ggd=0.1$} && \multicolumn{2}{c}{$\ggd=1$}\\
    \cmidrule{3-10}
    $\Delta t$ && $\Vert e_h\Vert_{L^\infty(\mathbb{T})}$ & $\Vert e_h\Vert_{L^1(\mathbb{T})}$ && $\Vert e_h\Vert_{L^\infty(\mathbb{T})}$ & $\Vert e_h\Vert_{L^1(\mathbb{T})}$ && $\Vert e_h\Vert_{L^\infty(\mathbb{T})}$ & $\Vert e_h\Vert_{L^1(\mathbb{T})}$\\
    \midrule
    0.1400 && 0.12924 & 0.03541 && 0.10110 & 0.02852 && 0.08216 & 0.02553\\
    0.0140 && 0.27148 & 0.05519 && 0.21489 & 0.04332 && 0.21576 & 0.03451\\
    0.0014 && 0.60189 & 0.05194 && 0.41119 & 0.04702 && 0.57972 & 0.03522\\
    \bottomrule
  \end{tabular}
\end{table}

\paragraph{Ghost-penalty stabilization}
\Cref{tab.ex2.th.gp-parameters} studies the effects of varying the ghost-penalty parameters for velocity and pressure compared to the results in \Cref{tab.ex2.th.viscosity} with $\nu = 0.01$, i.e., without grad-div stabilization and narrow-band stabilization. Decreasing the pressure ghost-penalty parameter $\gpp$ while keeping the velocity ghost-penalty parameter $\gpv$ fixed does not consistently reduce the pressure drag oscillations across the three CFL numbers. However, reducing $\gpv$ while keeping $\gpp$ fixed has a positive impact on these oscillations. Further reducing both parameters has additional benefits, but only for small CFL numbers, and the effect is smaller than the initial decrease in $\gpv$. This reduction is significant only in the $L^\infty$ error, not in the $L^1$ error. The best reduction in the $L^\infty$-norm for $C_{\text{CFL}} = 0.01$ is a factor of 6, yielding an error of less than 1\% relative to the signal amplitude.

\begin{table}
  \centering
  \caption{\hyperref[sec.numerics:subsec.ex2]{Example 2} --  Norms of pressure drag oscillations for varying ghost-penalty stabilization parameters. Isoparametric Taylor-Hood elements with $h=0.1$, $k=2$, $\nu=0.01$, $\ggd=0.0$ and narrow-band ghost-penalty stabilization.}
  \label{tab.ex2.th.gp-parameters}
  \begin{tabular}{ccc@{\hspace*{8pt}}ccc@{\hspace*{8pt}}ccc@{\hspace*{8pt}}c}
    \toprule
      && \multicolumn{2}{c}{$\gpv=0.1,\gpp=0.01$} && \multicolumn{2}{c}{$\gpv=0.01,\gpp=0.1$} && \multicolumn{2}{c}{$\gpv=0.01,\gpp=0.01$}\\
    \cmidrule{3-10}
    $\Delta t$ && ~$\Vert e_h\Vert_{L^\infty(\mathbb{T})}$~ & ~$\Vert e_h\Vert_{L^1(\mathbb{T})}$~ && ~$\Vert e_h\Vert_{L^\infty(\mathbb{T})}$~ & ~$\Vert e_h\Vert_{L^1(\mathbb{T})}$~ && ~$\Vert e_h\Vert_{L^\infty(\mathbb{T})}$~ & ~$\Vert e_h\Vert_{L^1(\mathbb{T})}$~\\
    \midrule
    0.1400 && 0.17273 & 0.04391 && 0.05249 & 0.01524 && 0.06037 & 0.01956\\
    0.0140 && 0.37452 & 0.06733 && 0.10630 & 0.01910 && 0.09852 & 0.01946\\
    0.0014 && 0.62735 & 0.05895 && 0.26992 & 0.02452 && 0.19471 & 0.02186\\
    \bottomrule
  \end{tabular}
\end{table}

\paragraph{Combined parameter tuning}
\Cref{tab.ex2.th.gp+gd} explores whether further improvements can be achieved by optimizing both the grad-div and ghost-penalty parameters. Indeed, adding grad-div stabilization leads to an additional reduction in spurious oscillations. For $C_{\text{CFL}} = 0.01$, a larger pressure ghost-penalty parameter proves beneficial, achieving a reduction by a factor of up to 9.1 compared with the default parameter case and an error of around 0.55\% relative to the signal amplitude.

\begin{table}
  \centering
  \caption{\hyperref[sec.numerics:subsec.ex2]{Example 2} --  Norms of pressure drag oscillations for varying ghost-penalty stabilization parameters and grad-div stabilization enabled. Isoparametric Taylor-Hood elements with $h=0.1$, $k=2$, $\nu=0.01$, $\ggd=0.1$, $\gpv=0.01$, and narrow-band ghost-penalty stabilization.}
  \label{tab.ex2.th.gp+gd}
  \begin{tabular}{ccc@{\hspace*{8pt}}ccc@{\hspace*{8pt}}c}
    \toprule
      && \multicolumn{2}{c}{ $\gpp=0.1$} && \multicolumn{2}{c}{ $\gpp=0.01$}\\
    \cmidrule{3-7}
    $\Delta t$ && ~$\Vert e_h\Vert_{L^\infty(\mathbb{T})}$~ & ~$\Vert e_h\Vert_{L^1(\mathbb{T})}$~ && ~$\Vert e_h\Vert_{L^\infty(\mathbb{T})}$~ & ~$\Vert e_h\Vert_{L^1(\mathbb{T})}$~ \\
    \midrule
    0.1400 && 0.04648 & 0.01359 && 0.03932 & 0.01456\\
    0.0140 && 0.08134 & 0.01655 && 0.07915 & 0.01822\\
    0.0014 && 0.12923 & 0.01995 && 0.15431 & 0.02339\\
    \bottomrule
  \end{tabular}
\end{table}

\begin{table}
  \centering
  \caption{\hyperref[sec.numerics:subsec.ex2]{Example 2} -- Norms of pressure drag oscillations for different geometry approximations. (Isoparametric) Taylor-Hood elements with $h=0.1$, $k=2$, $\nu=0.01$, $\gpv=\gpp=0.1$, $\ggd=0$ and narrow-band ghost-penalty stabilisation.}
  \label{tab.ex2.th.geometry}
  \begin{tabular}{ccc@{\hspace*{8pt}}cc@{\hspace*{8pt}}cc@{\hspace*{8pt}}cc@{\hspace*{8pt}}c}
    \toprule
    \multicolumn{2}{r}{Geom. Approx.:}& \multicolumn{2}{c}{Isoparametric} & \multicolumn{2}{c}{P1} & \multicolumn{2}{c}{P1 - 1 Subdivision} & \multicolumn{2}{c}{P1 - 2 Subdivisions}\\
    \cmidrule{3-10}
    $\Delta t$ & $k$ & $\Vert e_h\Vert_{L^\infty(\mathbb{T})}$ & $\Vert e_h\Vert_{L^1(\mathbb{T})}$ & $\Vert e_h\Vert_{L^\infty(\mathbb{T})}$ & $\Vert e_h\Vert_{L^1(\mathbb{T})}$ & $\Vert e_h\Vert_{L^\infty(\mathbb{T})}$ & $\Vert e_h\Vert_{L^1(\mathbb{T})}$ & $\Vert e_h\Vert_{L^\infty(\mathbb{T})}$ & $\Vert e_h\Vert_{L^1(\mathbb{T})}$\\
    \midrule
    0.1400 & 2 & 0.17018 & 0.04400 & 0.16026 & 0.04326 & 0.16194 & 0.04355 & 0.16255 & 0.04363\\
    0.0140 & 2 & 0.27549 & 0.06265 & 0.31770 & 0.07258 & 0.32314 & 0.07292 & 0.32318 & 0.07295\\
    0.0014 & 2 & 1.17637 & 0.06584 & 1.16602 & 0.06816 & 1.17069 & 0.06829 & 1.16952 & 0.06832\\[2pt]
    0.1400 & 3 & 0.11748 & 0.03948 & 0.11269 & 0.04053 & 0.12144 & 0.04198 & 0.12238 & 0.04217\\
    0.0140 & 3 & 0.06412 & 0.01314 & 0.06191 & 0.01347 & 0.06308 & 0.01368 & 0.06321 & 0.01371\\
    0.0014 & 3 & 0.10590 & 0.01089 & 0.11795 & 0.01216 & 0.12619 & 0.01218 & 0.12722 & 0.01221\\
    \bottomrule
  \end{tabular}
\end{table}

\paragraph{Geometry approximation}
\Cref{tab.ex2.th.geometry} studies both the use of higher-order elements and the effect of geometry approximation error. The isoparametric results for $k = 2$ repeat those in \Cref{tab.ex2.th.viscosity} for $\nu = 0.01$. We do not observe consistent increases or decreases when comparing the results with isoparametric elements to those with a $P^1$ level set approximation. Furthermore, the absolute differences are relatively small for $k = 2$ and $k = 3$, and reducing the geometric error by using a composite quadrature with additional subdivisions does not significantly affect the results. This suggests that the preservation of area is not a major source of spurious pressure-force oscillations in this Eulerian cutFEM approach.

\begin{figure}
  \centering
  \includegraphics{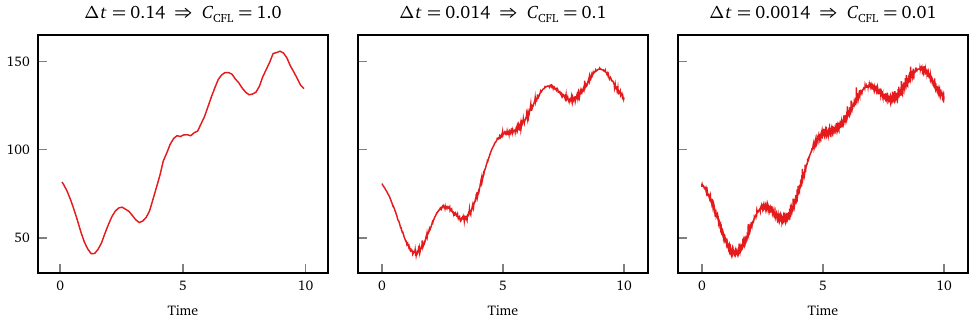}
  \captionof{figure}{\hyperref[sec.numerics:subsec.ex2]{Example 2} -- $p - p_{mean}$ integrated over the moving boundary for the isoparametric Taylor-Hood method using narrow-band ghost-penalty stabilization and parameters $\nu=0.01$, $h=0.1$, $k=2$, $\gpv=\gpp=0.1$ and $\ggd=0$.}
  \label{fig.ex2-pressure-mean}
\end{figure}

To confirm that the observed oscillations in the pressure drag coefficient are localized, rather than simply due to spurious oscillations in the mean pressure (which could suggest a lack of conservation of total area), we analyze $\int_{\Gamma(t)} (p - p_{\text{mean}}) \, \mathrm{d}s$, with $p_{\text{mean}} = \int_{\Omega(t)} p \, \mathrm{d}x$, in \Cref{fig.ex2-pressure-mean}. As expected, we also observe increasing spurious oscillations in this (normalized) quantity.

\paragraph{Higher-order elements}
Using higher-order ($k = 3$) elements, i.e., improving the spatial approximation of the method in \Cref{tab.ex2.th.geometry}, significantly improves the results for small CFL numbers. Indeed, the oscillations decrease by a factor of 10.8 for $C_{\text{CFL}} = 0.01$ compared with the the results in \Cref{tab.ex2.th.viscosity}, which is greater than the effect of tuning parameters for $k=2$.

\paragraph{Ghost-penalty facets}
\Cref{tab.ex2.stabil+ho} investigates if choosing a smaller set of ghost-penalty facets (in the spirit of weak AggFEM in the narrow-band setting) has a positive effect (for $L = 1$). We observe that the magnitude of the pressure oscillations increases.

\paragraph{Finer meshes}
We also investigate the effect of decreasing the mesh size (thereby increasing the CFL number) in \Cref{tab.ex2.stabil+ho}. We again observe a significant reduction in the spurious oscillations in both the $L^\infty$- and $L^1$-norms. Additionally, we do not see an increase in the oscillations for decreasing CFL numbers, similar to the results for $k = 3$ in \Cref{tab.ex2.th.geometry}. This aligns with our convergence results in \Cref{sec.numerics:subsec.convergence}, where the loss of control over the pressure in the $L^\infty(H^1)$-norm was delayed by increasing the spatial resolution.

\begin{table}
  \centering
  \caption{\hyperref[sec.numerics:subsec.ex2]{Example 2} --  Norms of pressure drag oscillations using finer meshes and patch-wise ghost-penalty stabilization. Isoparametric Taylor-Hood elements with $k=2$, $\nu=0.01$, $\gpv=\gpp=0.1$, $\ggd=0$.}
  \label{tab.ex2.stabil+ho}
  \begin{tabular}{ccc@{\hspace*{8pt}}ccc@{\hspace*{8pt}}ccc@{\hspace*{8pt}}c}
    \toprule
    && \multicolumn{2}{c}{Patch, $h=0.1$}&& \multicolumn{2}{c}{Band, $h=0.1$} && \multicolumn{2}{c}{Band, $h=0.05$}\\
    \cmidrule{3-10}
    $\Delta t$ && ~$\Vert e_h\Vert_{L^\infty(\mathbb{T})}$~ & ~$\Vert e_h\Vert_{L^1(\mathbb{T})}$~ && ~$\Vert e_h\Vert_{L^\infty(\mathbb{T})}$~ & ~$\Vert e_h\Vert_{L^1(\mathbb{T})}$~&& ~$\Vert e_h\Vert_{L^\infty(\mathbb{T})}$~ & ~$\Vert e_h\Vert_{L^1(\mathbb{T})}$~ \\
    \midrule
    0.1400 && n/a     & n/a     && 0.17018 & 0.04400 && 0.10390 & 0.03145\\
    0.0140 && 0.31563 & 0.07067 && 0.27549 & 0.06265 && 0.07208 & 0.01480\\
    0.0014 && 1.82717 & 0.09786 && 1.17637 & 0.06584 && 0.08262 & 0.01335\\
    \bottomrule
  \end{tabular}
\end{table}

\begin{table}
  \centering
  \caption{\hyperref[sec.numerics:subsec.ex2]{Example 2} --  Norms of pressure drag oscillations for varying viscosities and grad-div stabilization. Isoparametric Taylor-Hood elements with $h=0.1$, $k=2$, $\gpv=\gpp=0.1$ and the \emph{global stabilization} choice.}
  \label{tab.ex2.th.global-stabil}
  \begin{tabular}{cc@{\hspace*{5pt}}c@{\hspace*{8pt}}cc@{\hspace*{5pt}}c@{\hspace*{8pt}}cc@{\hspace*{5pt}}c@{\hspace*{8pt}}cc@{\hspace*{5pt}}c@{\hspace*{8pt}}c}
    \toprule
      && \multicolumn{2}{c}{$\nu=1,\ggd=0$} && \multicolumn{2}{c}{$\nu=0.1,\ggd=0$} && \multicolumn{2}{c}{$\nu=0.01,\ggd=0$} && \multicolumn{2}{c}{$\nu=0.01,\ggd=0.1$}\\
    \cmidrule{3-13}
    $\Delta t$ && $\Vert e_h\Vert_{L^\infty(\mathbb{T})}$ & $\Vert e_h\Vert_{L^1(\mathbb{T})}$ && $\Vert e_h\Vert_{L^\infty(\mathbb{T})}$ & $\Vert e_h\Vert_{L^1(\mathbb{T})}$ && $\Vert e_h\Vert_{L^\infty(\mathbb{T})}$ & $\Vert e_h\Vert_{L^1(\mathbb{T})}$ && $\Vert e_h\Vert_{L^\infty(\mathbb{T})}$ & $\Vert e_h\Vert_{L^1(\mathbb{T})}$\\
    \midrule
    0.1400 && 0.01504 & 0.00428 && 0.01358 & 0.00514 && 0.01311 & 0.00401 && 0.02754 & 0.00803\\
    0.0140 && 0.03248 & 0.00717 && 0.02641 & 0.00639 && 0.02047 & 0.00532 && 0.03549 & 0.00951\\
    0.0014 && 0.03292 & 0.00704 && 0.02799 & 0.00629 && 0.02555 & 0.00542 && 0.03999 & 0.00948\\
    \bottomrule
  \end{tabular}
\end{table}

\paragraph{Global stabilization}
\Cref{tab.ex2.th.global-stabil} presents the results for the global stabilization approach \eqref{eqn.all-facets} with $\nu \in \{1, 0.1, 0.01\}$ without grad-div stabilization and for $\nu = 0.01$ with grad-div stabilization. Consistent with our expectations from \Cref{sec.numerics:subsec.convergence}, the pressure-drag oscillations are very small, approximately 0.1\% relative to the signal amplitude, and show almost no dependence on $\Delta t$. Moreover, we observe that the oscillations do not increase with higher Reynolds numbers. Interestingly, the grad-div stabilization appears to have a slight deteriorating effect on the pressure oscillations in this case.
\begin{remark}[Computational effort]
 Compared to the narrow-band stabilization approach, the global stabilization approach incurs additional computational costs due to increased coupling across all interior facets and the inclusion of more degrees of freedom. For example, when factorizing the matrix of the discrete Stokes problem with global stabilization, the direct solver \texttt{pardiso} requires approximately twice the time needed for the same problem with narrow-band stabilization. At the same time, since the set of active degrees of freedom remains constant in the global stabilization approach, the factorized Jacobian matrix within Newton's method can be reused more frequently across different time steps compared to the narrow-band stabilization approach.
\end{remark}

\begin{figure}
  \centering
  \includegraphics{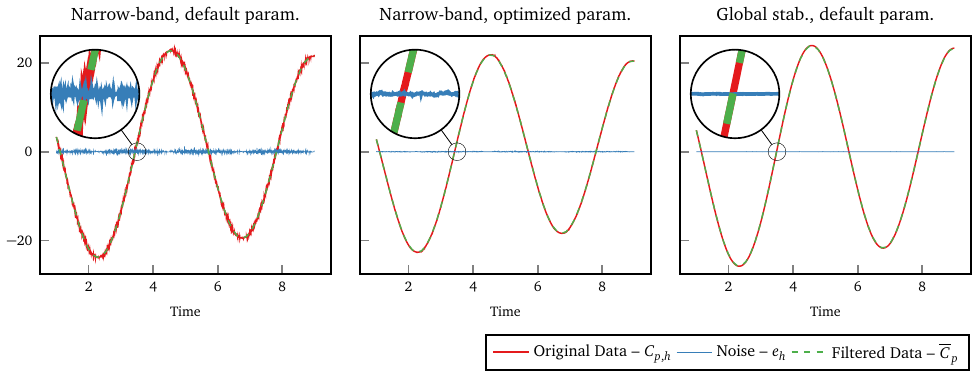}
  \caption{\hyperref[sec.numerics:subsec.ex2]{Example 2} -- Pressure drag coefficient, filtered signal and numerical noise. Computed using isoparametric Taylor-Hood elements with $k=2$, $h=0.1$, $\Delta t=0.0014$, viscosity $\nu=0.01$. Default parameters: $\gpv=\gpp=0.1$, $\ggd=0$. Optimized parameters: $\gpv=0.01$, $\gpp=\ggd=0.1$}
  \label{fig.ex2.k2ho}
\end{figure}

We visualize the data, post-processed signal and resulting noise in \Cref{fig.ex2.k2ho} for $\nu=0.01$, $C_{CFL}=0.01$ and the narrow band method with default parameters and optimized parameters in the case of the narrow-band stabilization and for the global ghost-penalty stabilization. Here, we see that the noise is clearly visible for the default parameter case, almost not noticeable for the optimized parameters and not visible to the naked eye for the global stabilization variant of the method.

\begin{figure}
  \centering
  \begin{minipage}[b]{14cm}
  \centering
  $t=0.0798$
  \hspace*{78pt}
  $t=0.0812$
  \hspace*{78pt}
  $t=0.0826$\\[2pt]
  \includegraphics[width=4.5cm, trim={0 2cm 0 1.5cm}, clip]{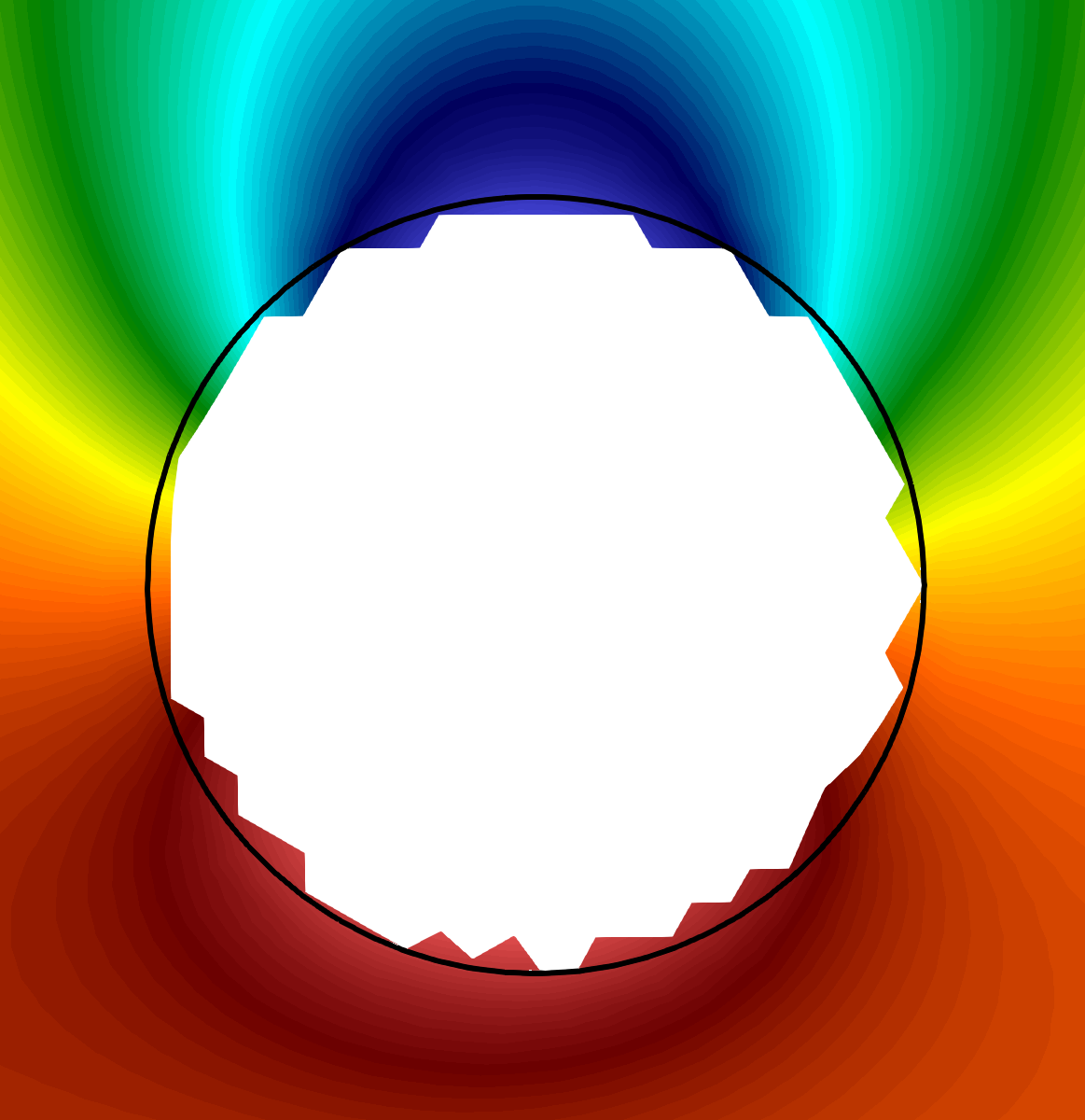}
  \includegraphics[width=4.5cm, trim={0 2cm 0 1.5cm}, clip]{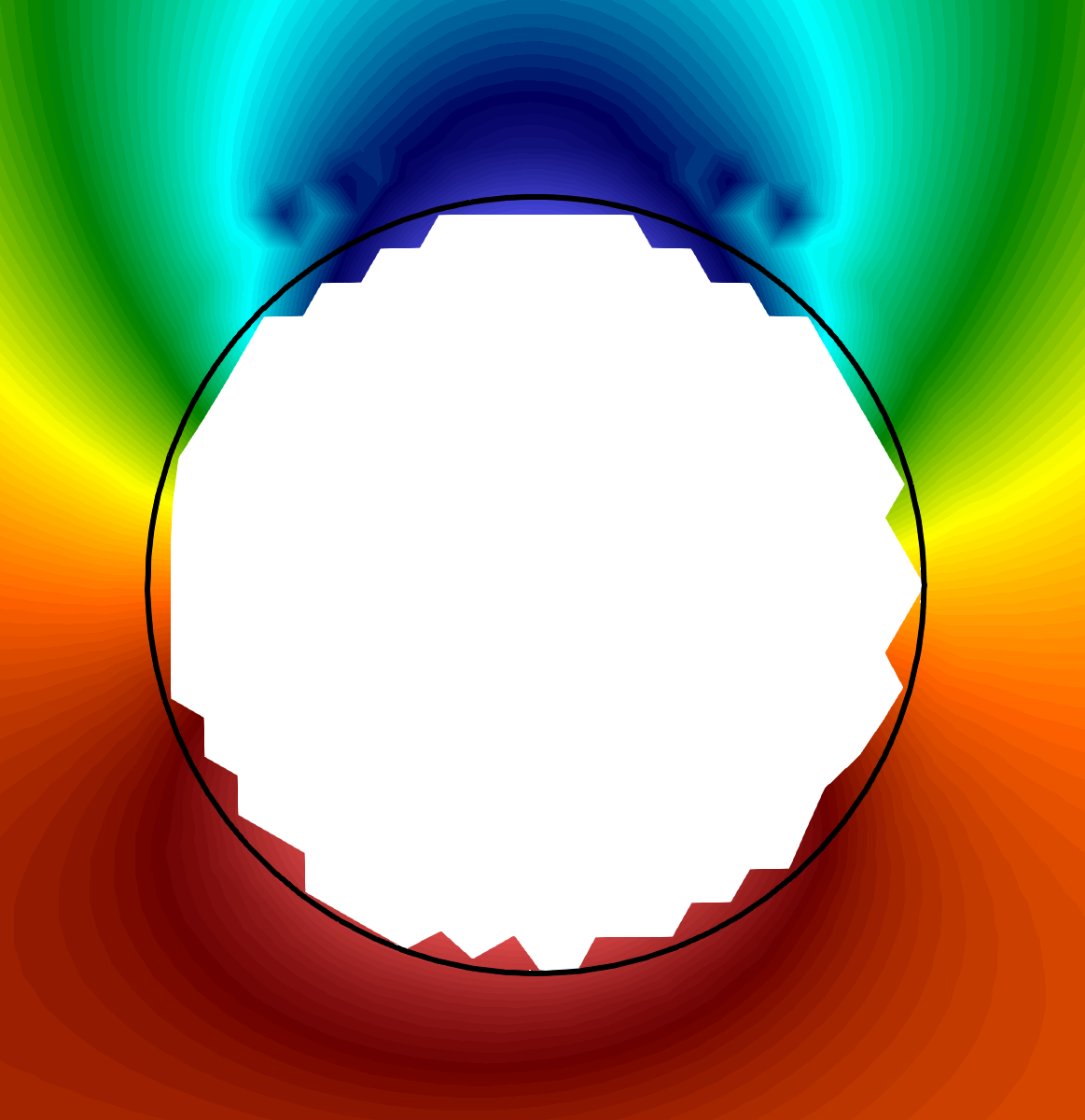}
  \includegraphics[width=4.5cm, trim={0 2cm 0 1.5cm}, clip]{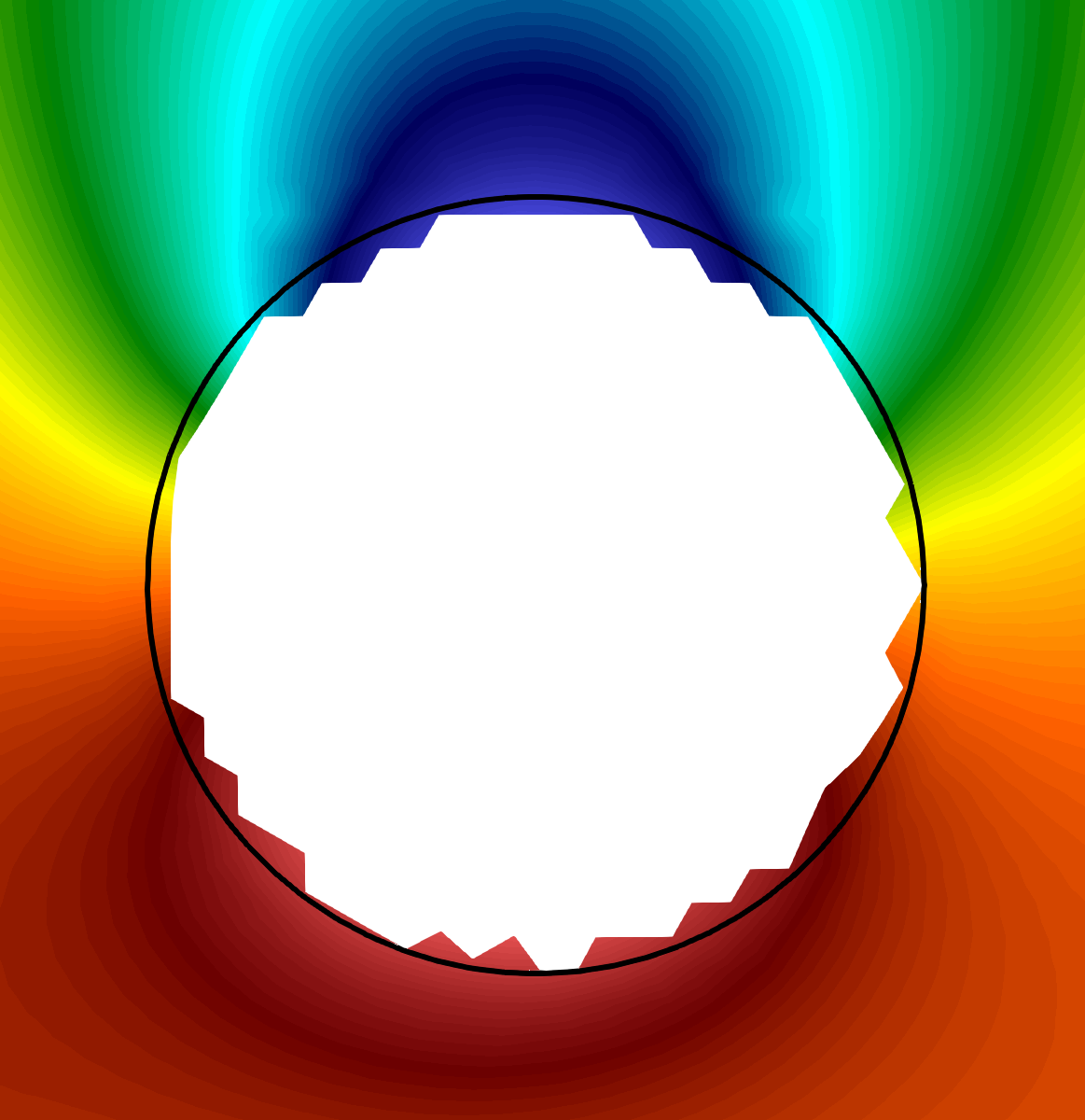}\\[2pt]
  \includegraphics[width=4.5cm, trim={0 2cm 0 1.5cm}, clip]{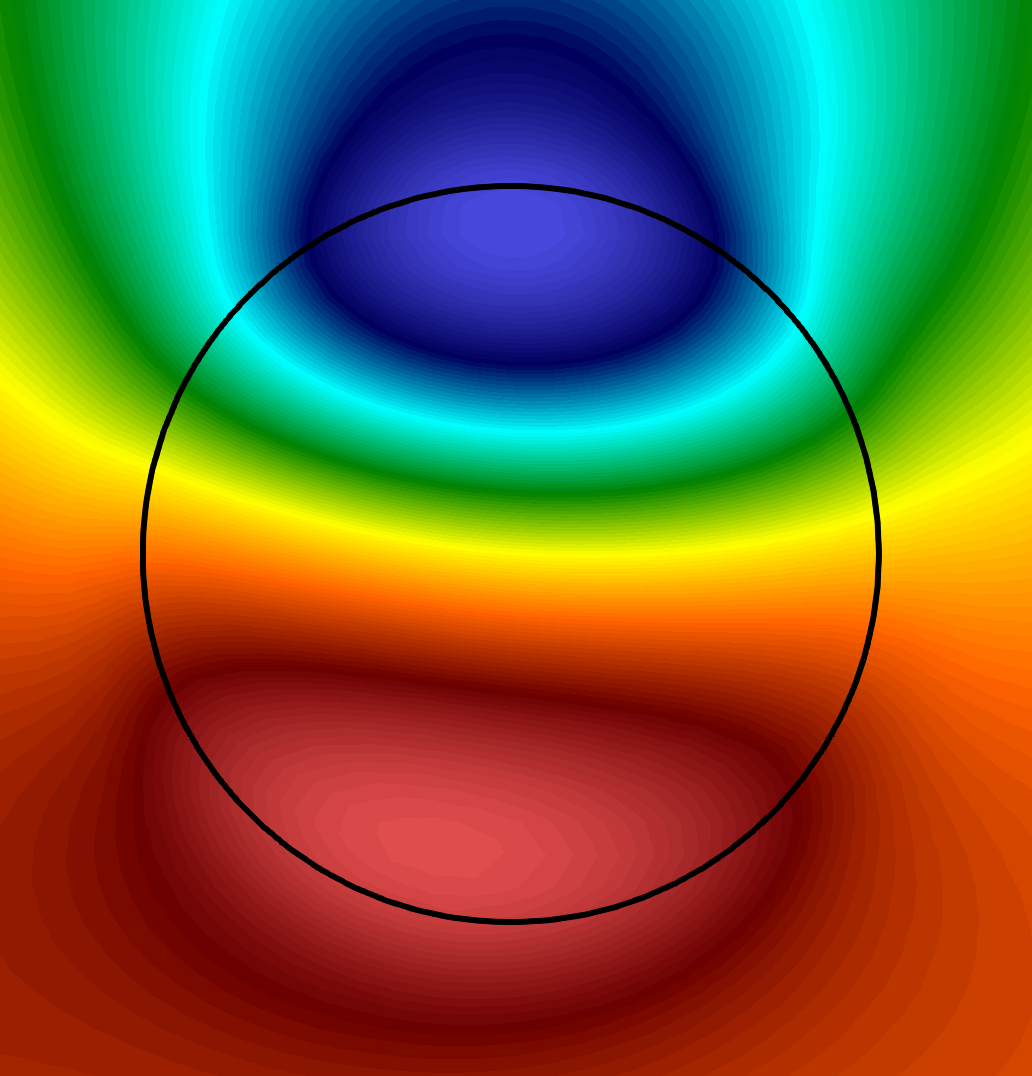}
  \includegraphics[width=4.5cm, trim={0 2cm 0 1.5cm}, clip]{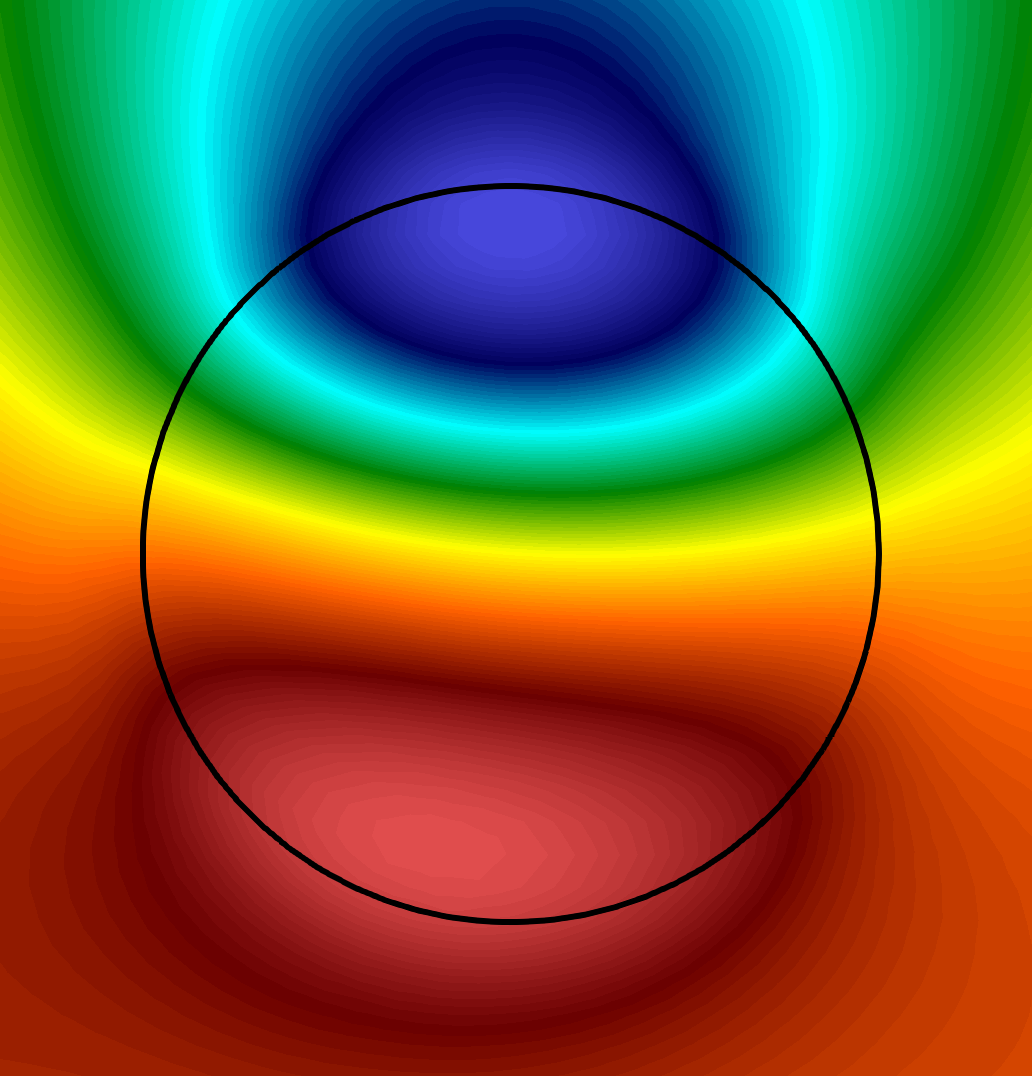}
  \includegraphics[width=4.5cm, trim={0 2cm 0 1.5cm}, clip]{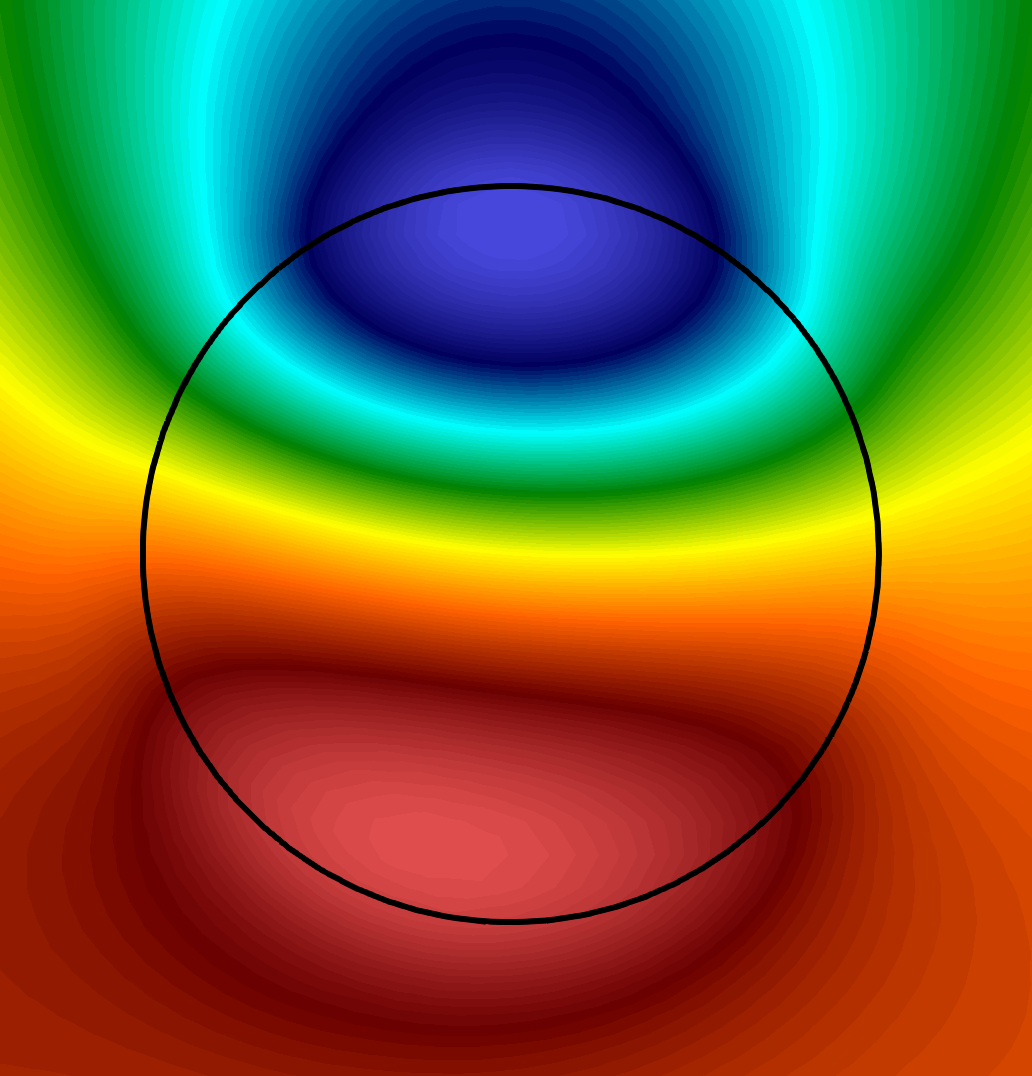}
  \end{minipage}
  \includegraphics[width=1.8cm]{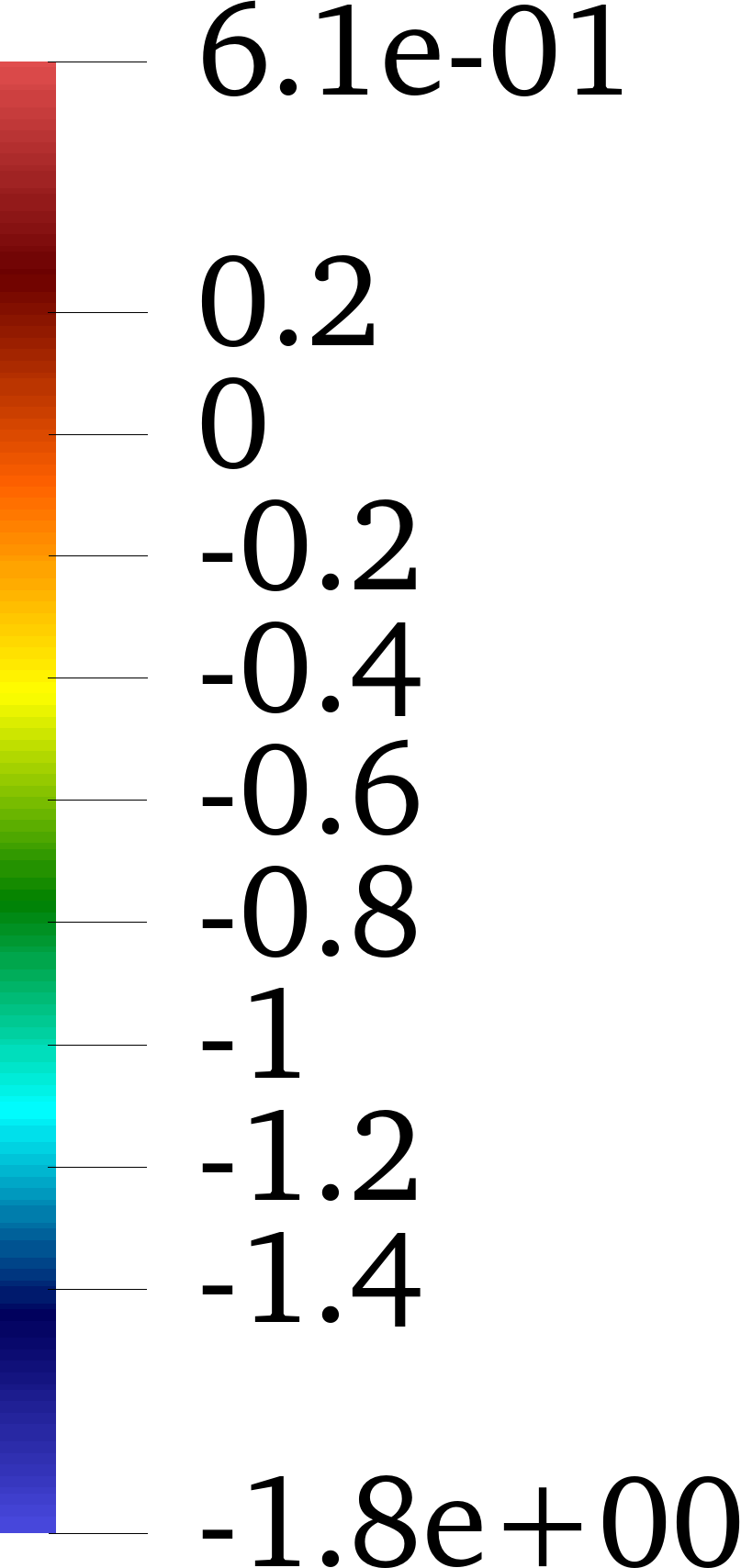}
  \caption{\hyperref[sec.numerics:subsec.ex2]{Example 2} -- Pressure solution near the moving cylinder at times $t=0.0798, 0.0812, 0.0826$. Computed with isoparametric Taylor-Hood elements with $k=2$, $h=0.1$, $\Delta t=0.0014$, $\nu=0.01$, $\gpv=\gpp=0.1$, $\ggd=0$. Top: Narrow-band stab.; Bottom: Global stab.}
  \label{fig.ex2.picture.pre}
\end{figure}

\begin{figure}
  \centering
  \begin{minipage}[b]{14cm}
  \centering
  $t=0.0798$
  \hspace*{78pt}
  $t=0.0812$
  \hspace*{78pt}
  $t=0.0826$\\[2pt]
  \includegraphics[width=4.5cm, trim={0 2cm 0 1.5cm}, clip]{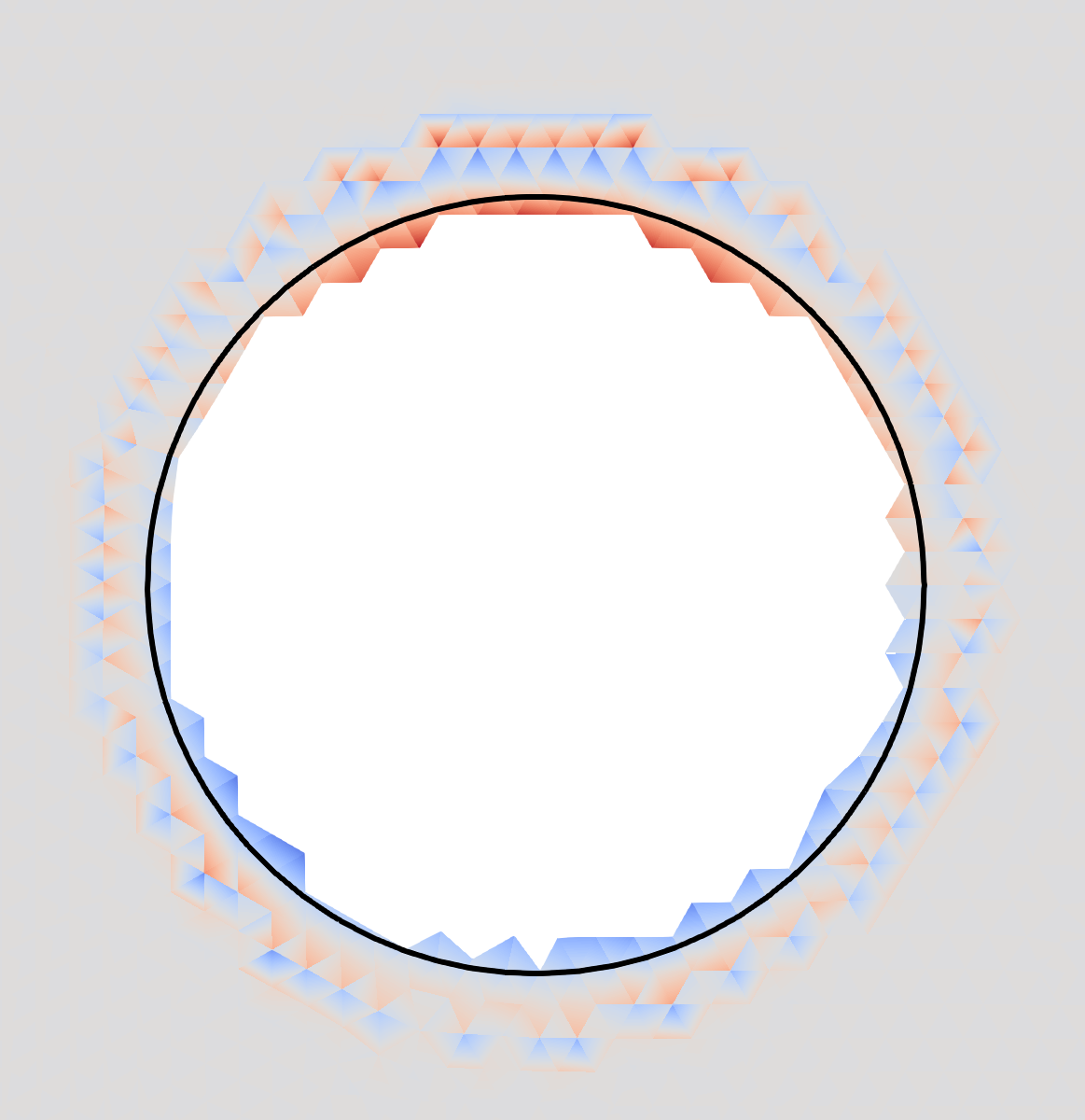}
  \includegraphics[width=4.5cm, trim={0 2cm 0 1.5cm}, clip]{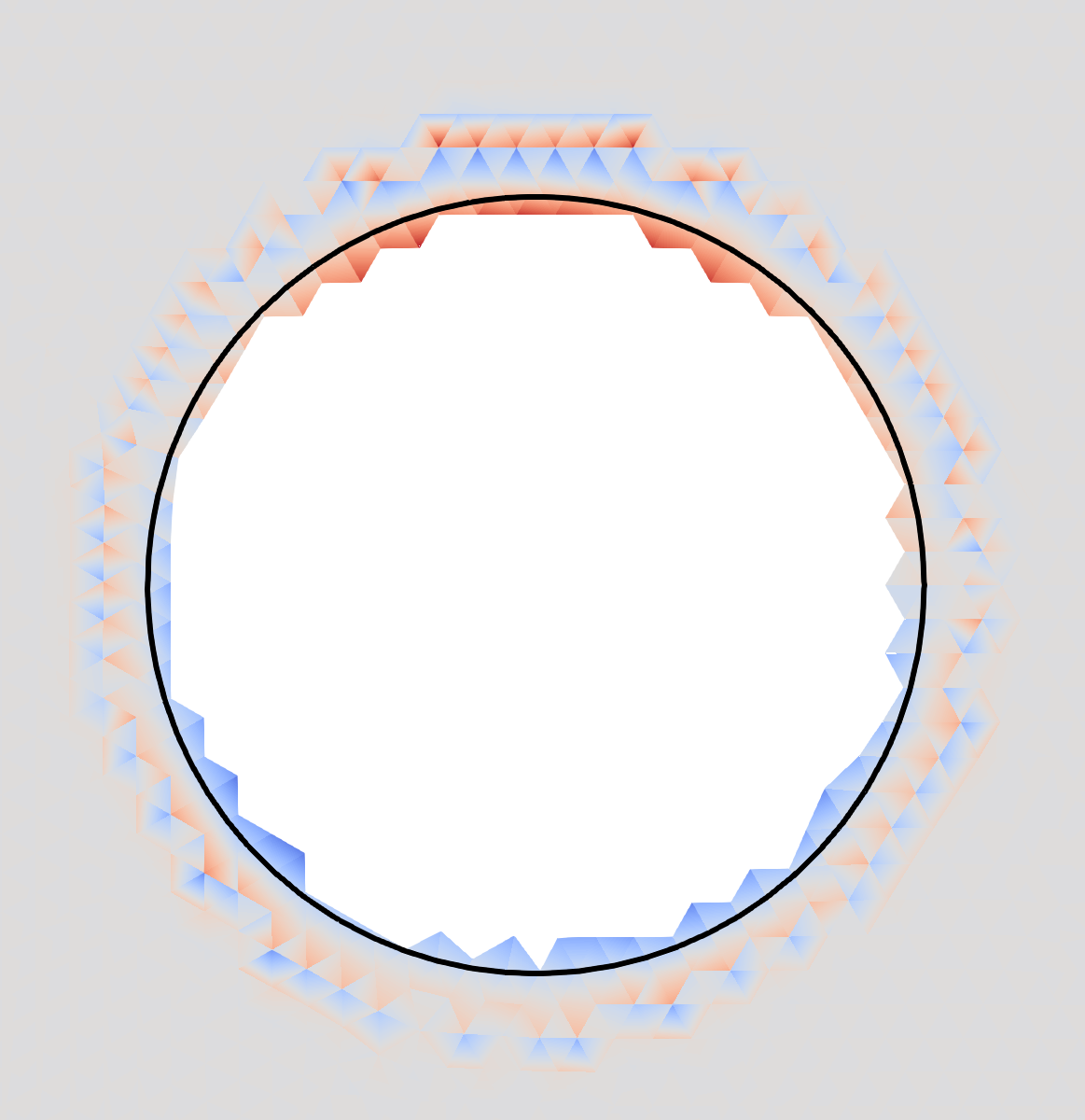}
  \includegraphics[width=4.5cm, trim={0 2cm 0 1.5cm}, clip]{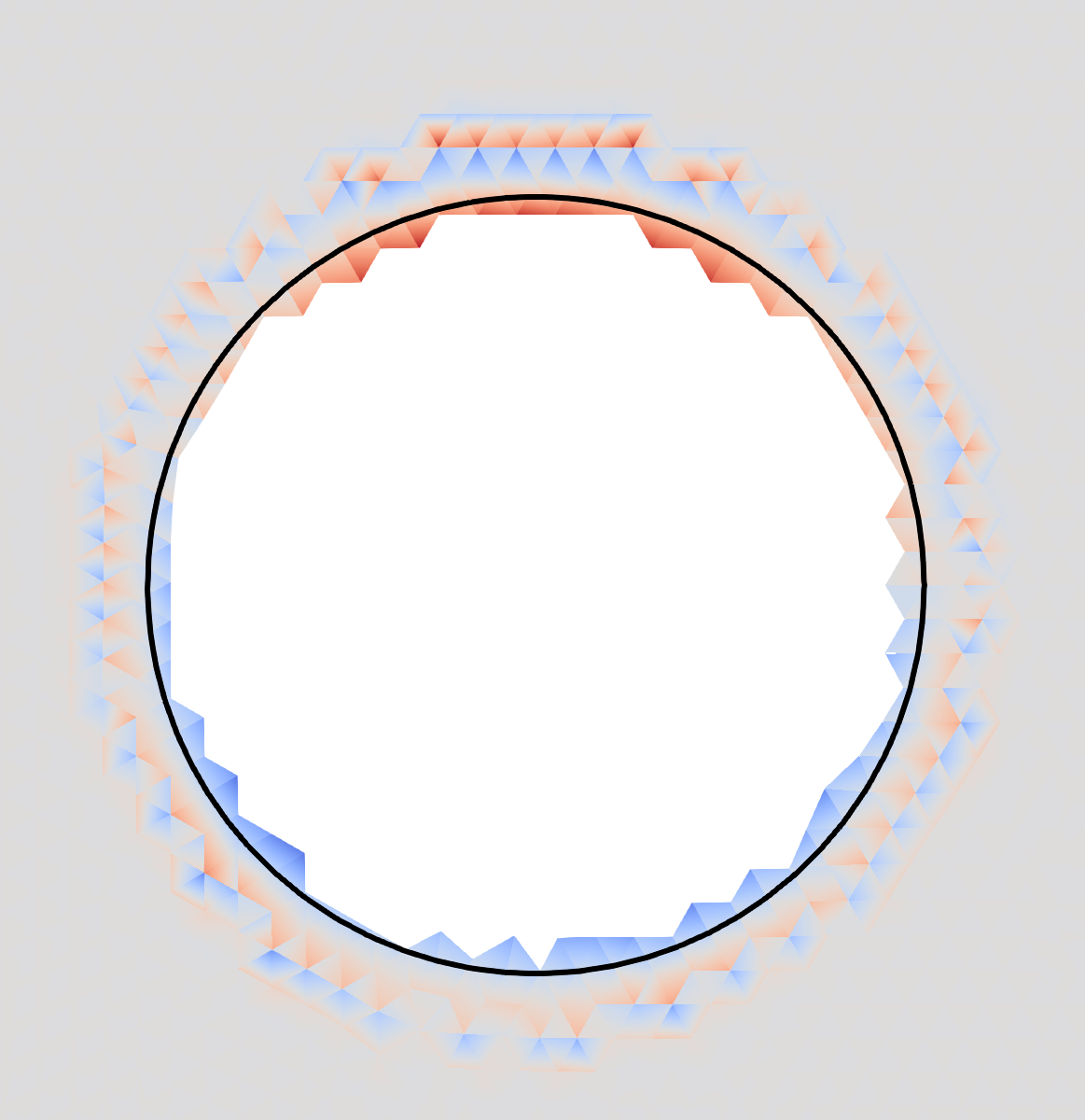}\\[2pt]
  \includegraphics[width=4.5cm, trim={0 2cm 0 1.5cm}, clip]{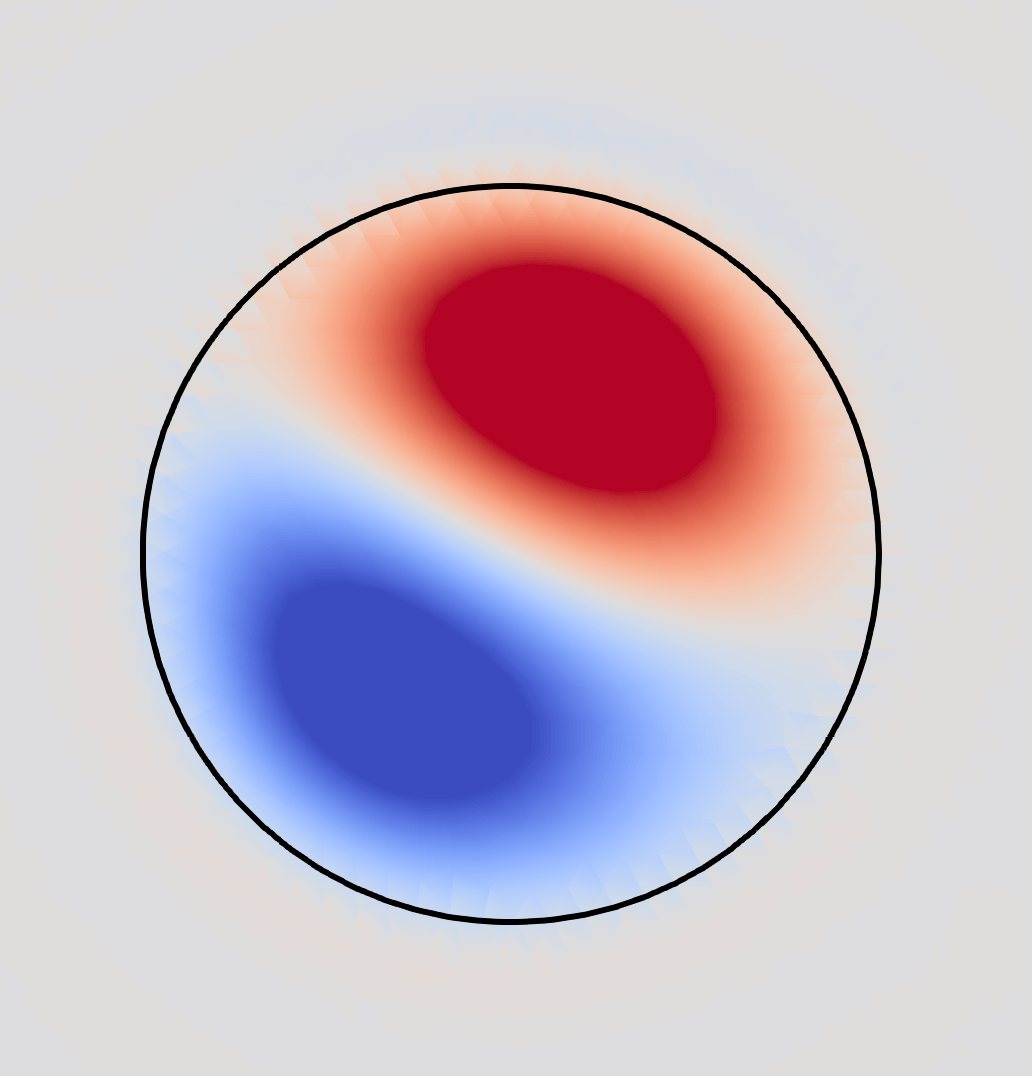}
  \includegraphics[width=4.5cm, trim={0 2cm 0 1.5cm}, clip]{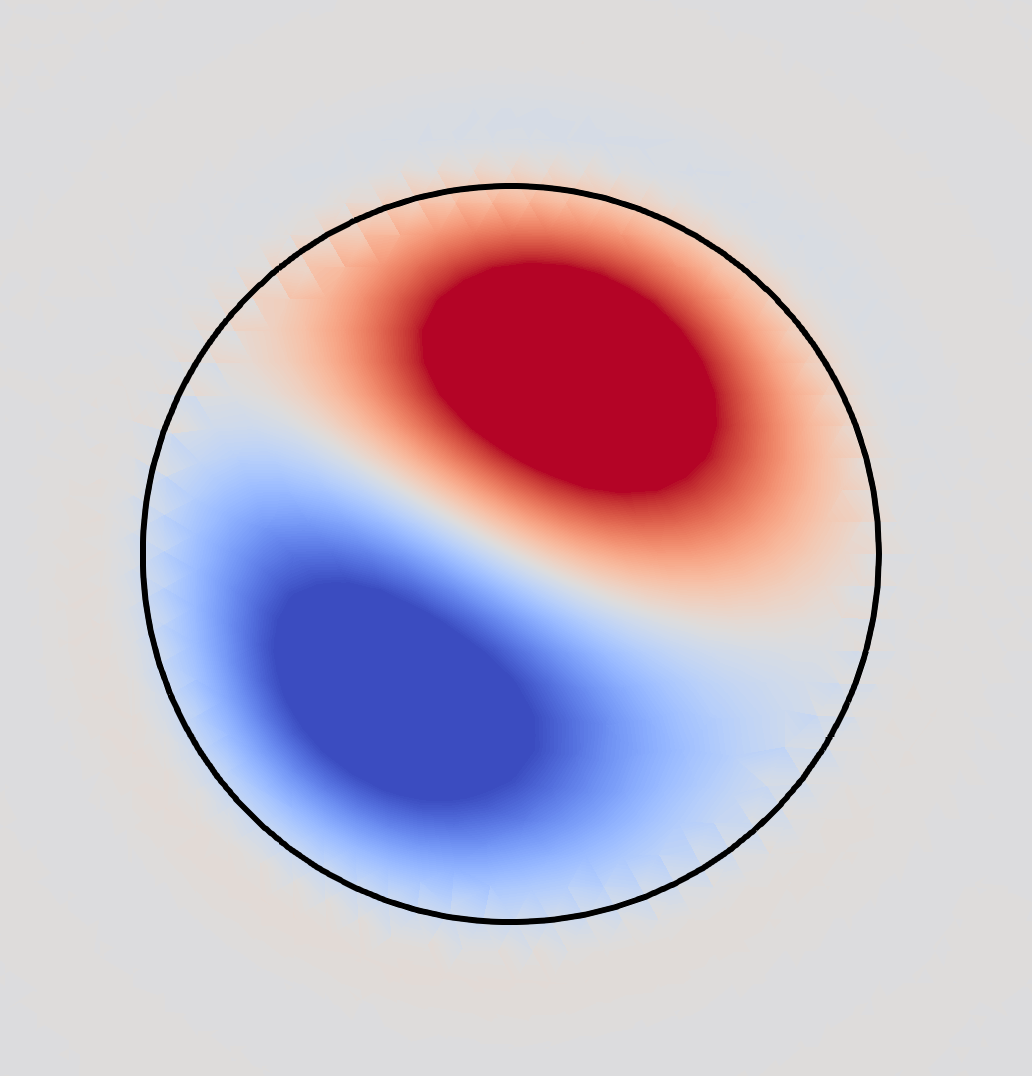}
  \includegraphics[width=4.5cm, trim={0 2cm 0 1.5cm}, clip]{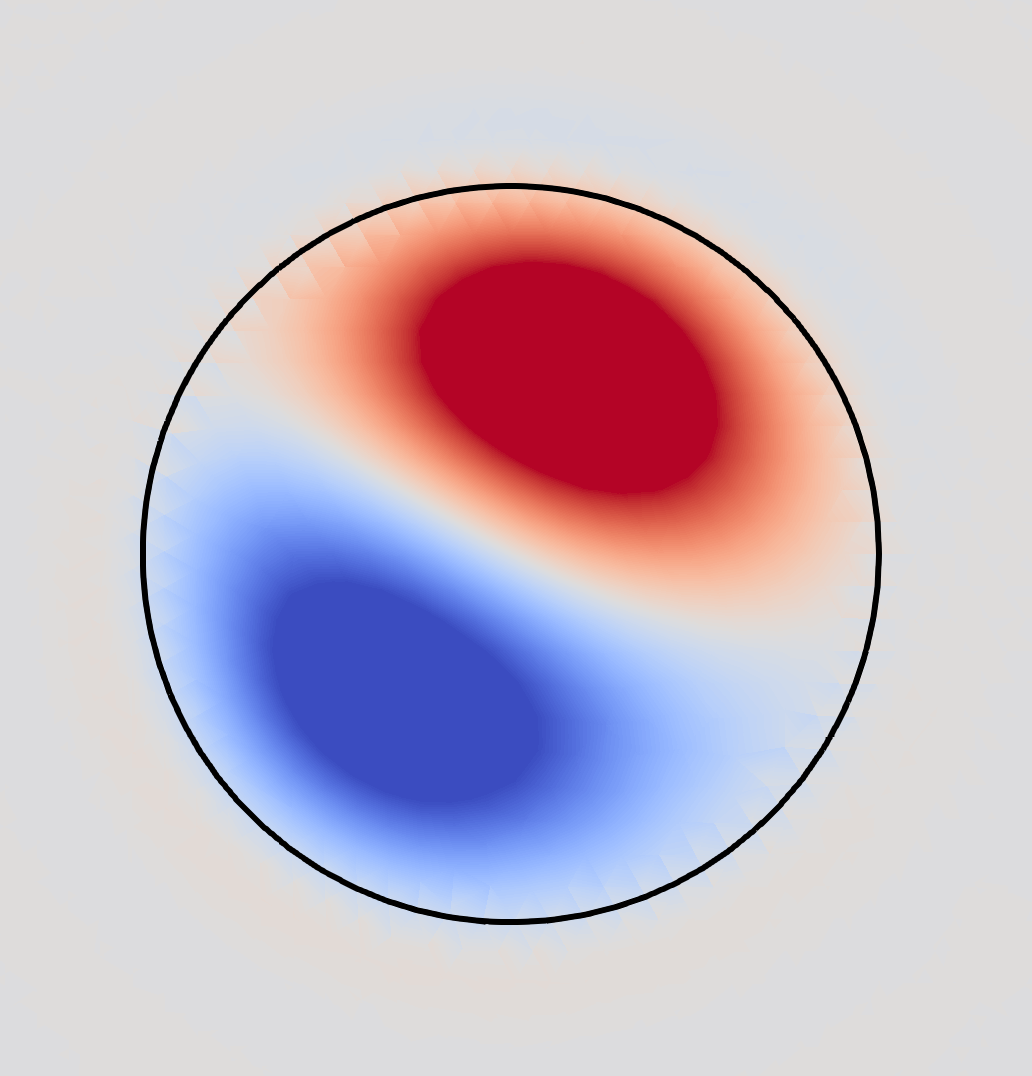}
  \end{minipage}
  \includegraphics[width=1.8cm]{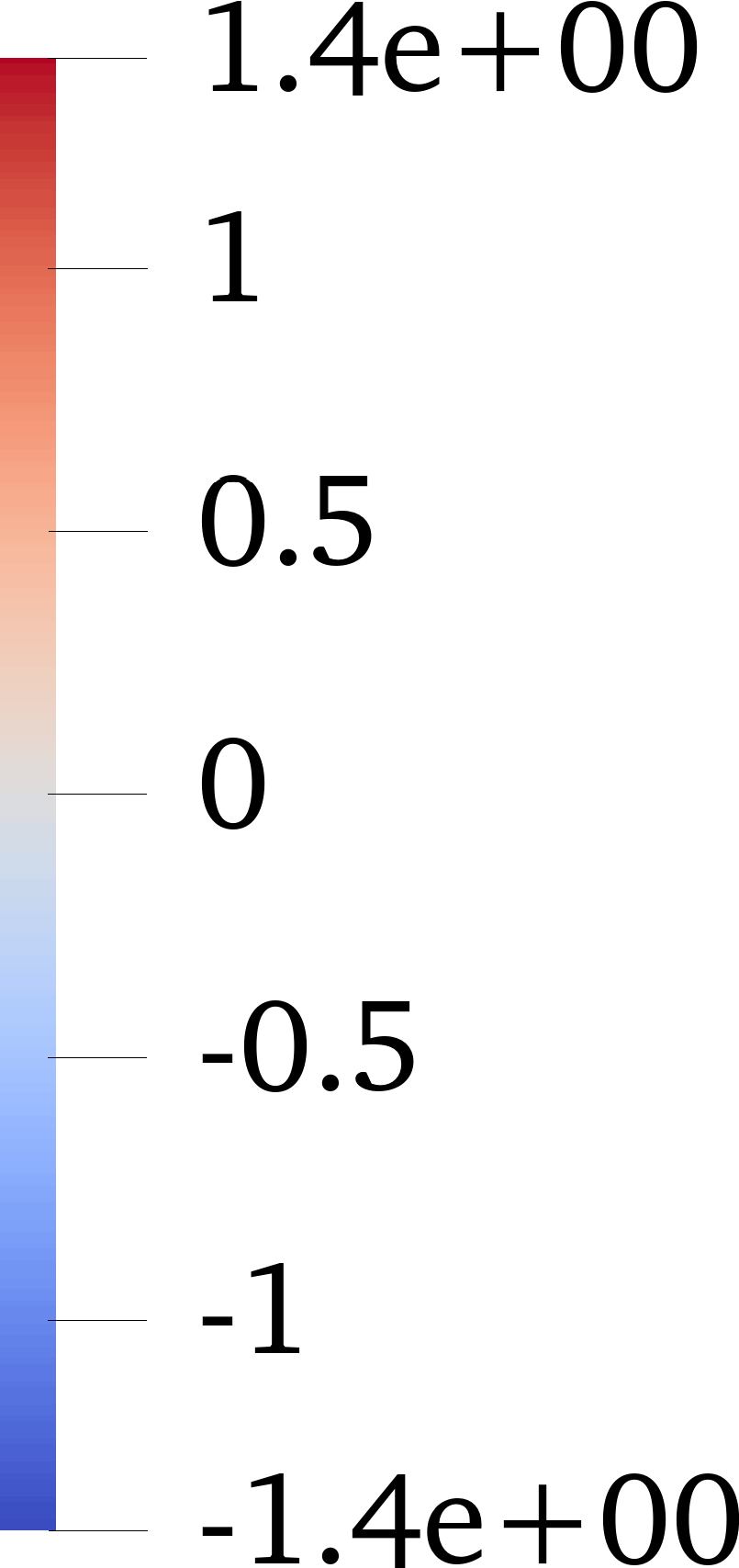}
  \caption{\hyperref[sec.numerics:subsec.ex2]{Example 2} -- Divergence of the velocity solution near the moving cylinder at times $t=0.0798, 0.0812, 0.0826$. Computed using isoparametric Taylor-Hood elements with $k=2$, $h=0.1$, $\Delta t=0.0014$, $\nu=0.01$, $\gpv=\gpp=0.1$, $\ggd=0$. Top: Narrow-band stab.; Bottom: Global stab.}
  \label{fig.ex2.picture.div}
\end{figure}

\begin{figure}
  \centering
  \begin{minipage}[b]{14cm}
  \centering
  $t=0.0798$
  \hspace*{78pt}
  $t=0.0812$
  \hspace*{78pt}
  $t=0.0826$\\[2pt]
  \includegraphics[width=4.5cm, trim={0 2cm 0 1.5cm}, clip]{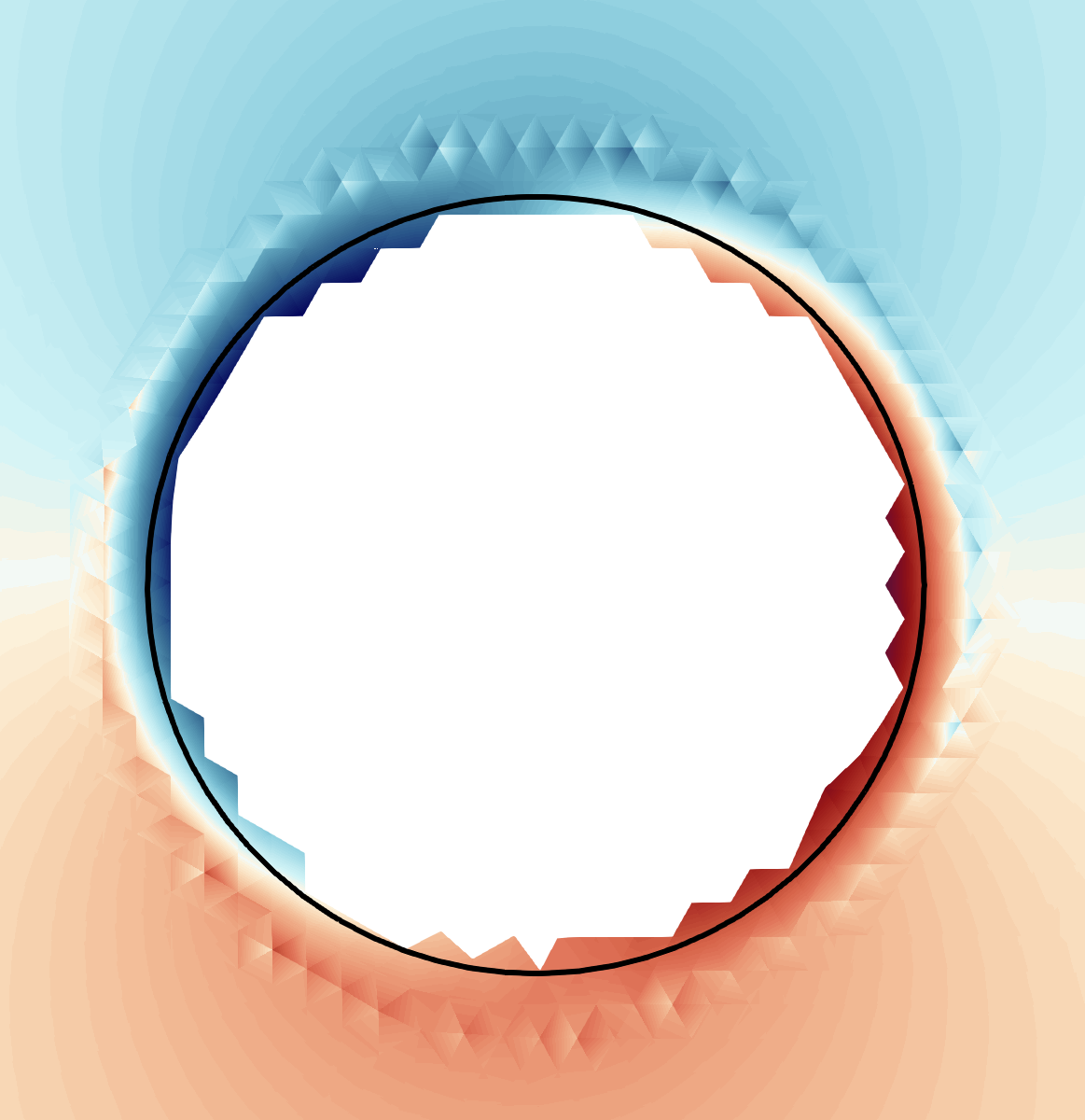}
  \includegraphics[width=4.5cm, trim={0 2cm 0 1.5cm}, clip]{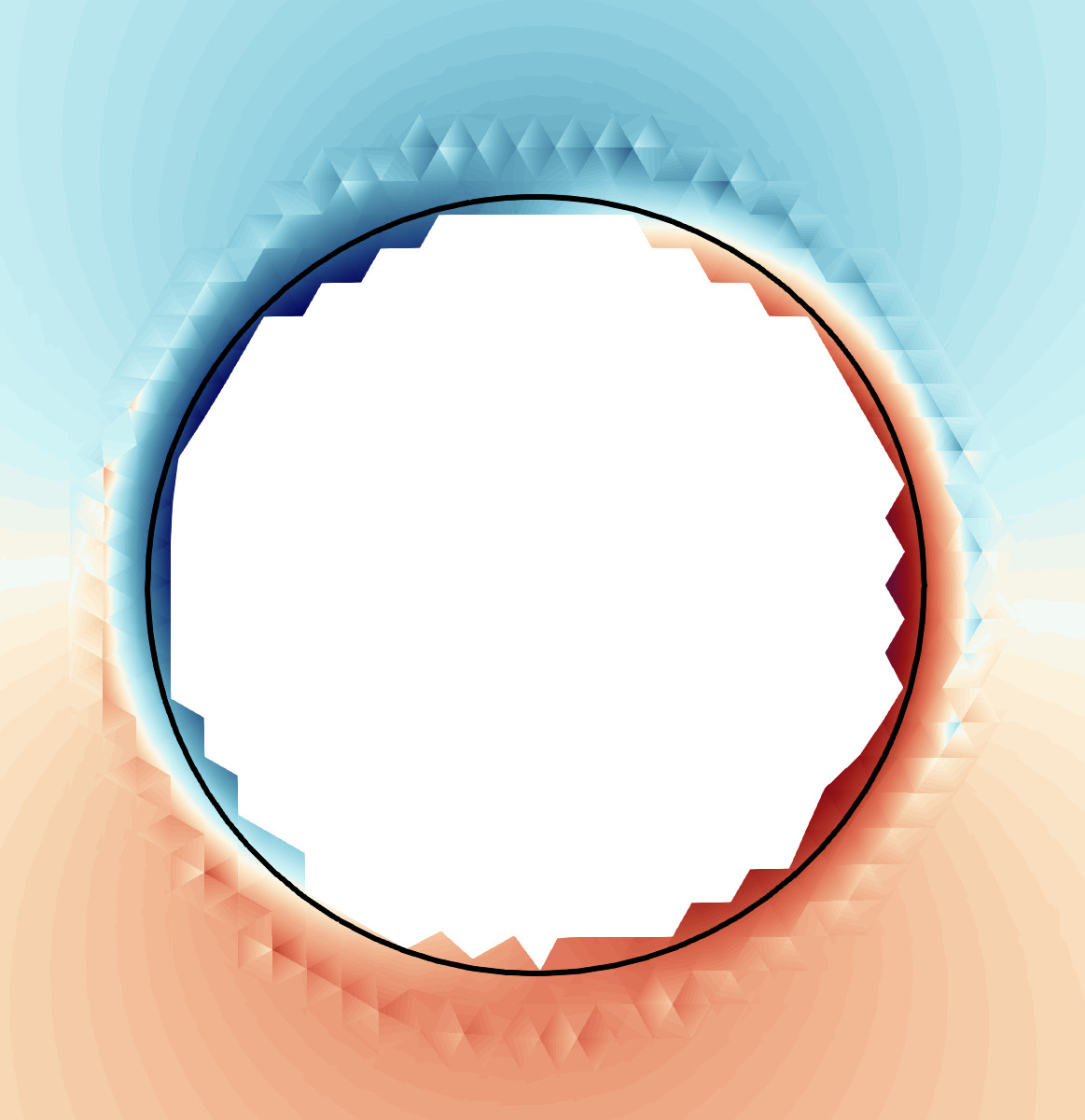}
  \includegraphics[width=4.5cm, trim={0 2cm 0 1.5cm}, clip]{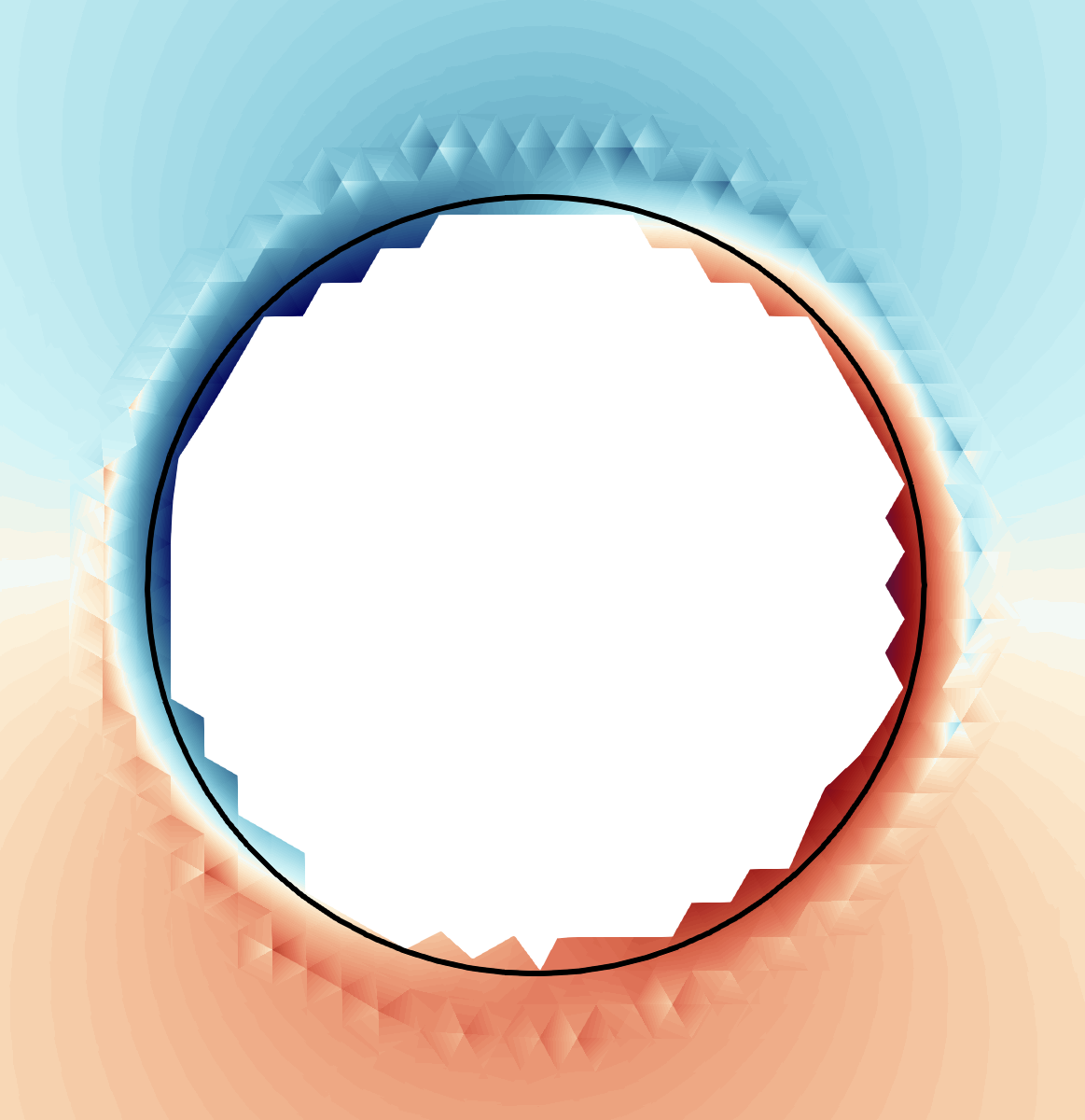}\\[2pt]
  \includegraphics[width=4.5cm, trim={0 2cm 0 1.5cm}, clip]{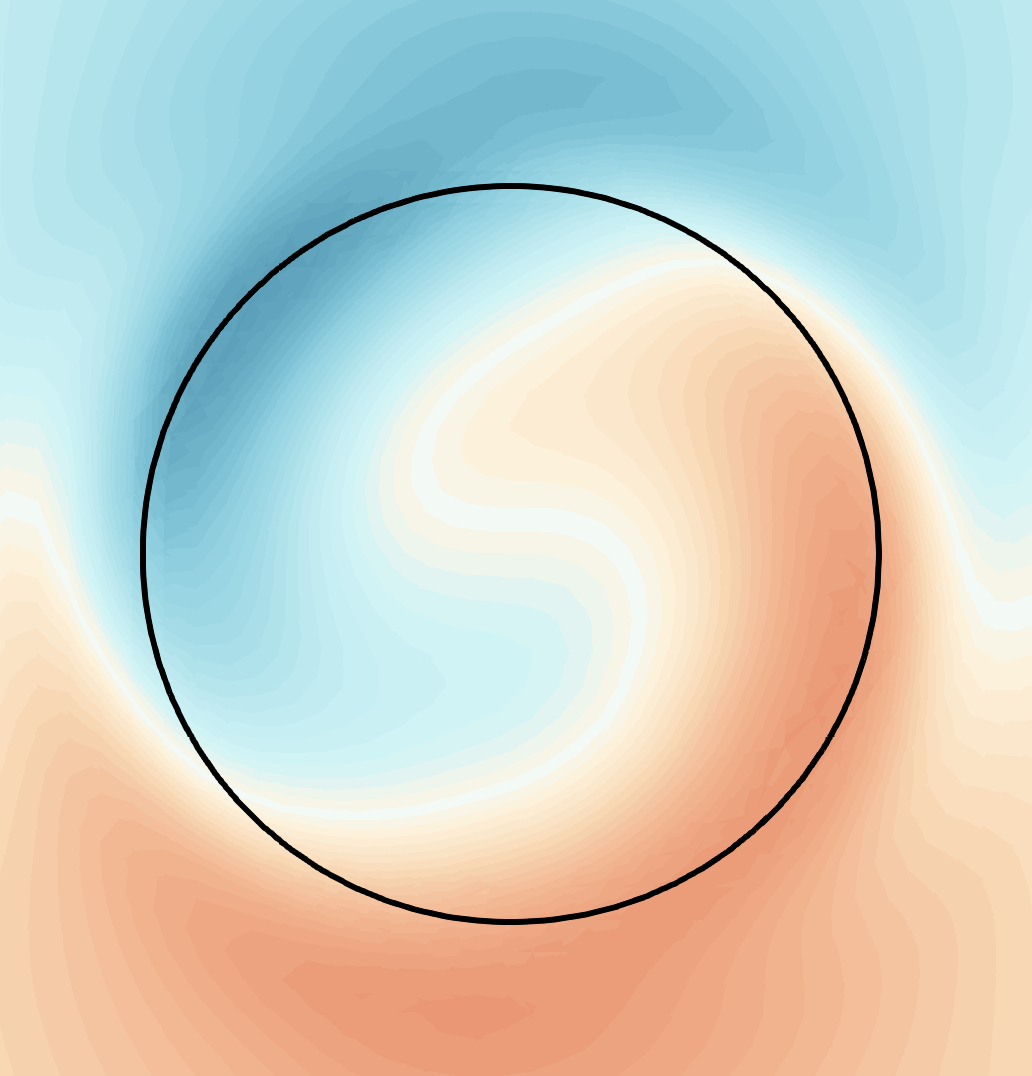}
  \includegraphics[width=4.5cm, trim={0 2cm 0 1.5cm}, clip]{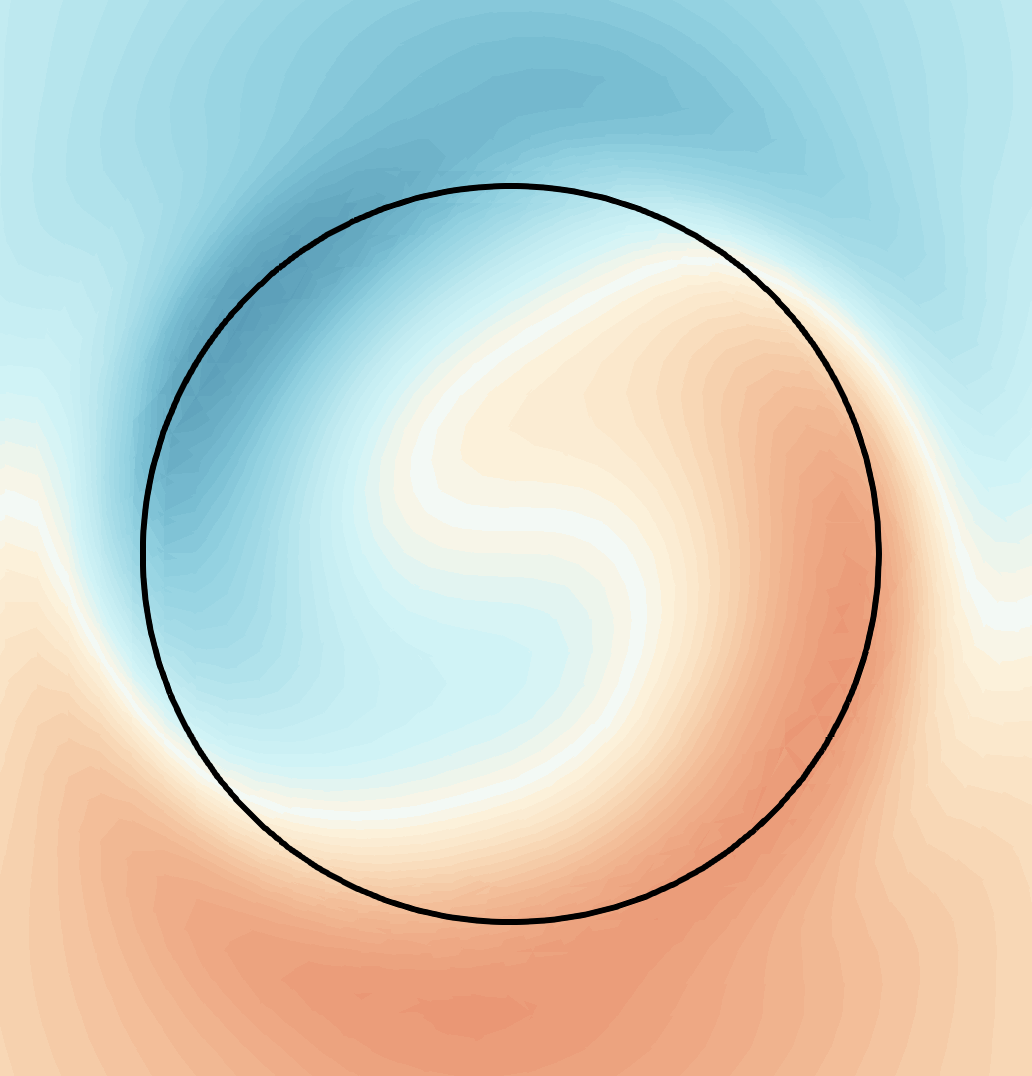}
  \includegraphics[width=4.5cm, trim={0 2cm 0 1.5cm}, clip]{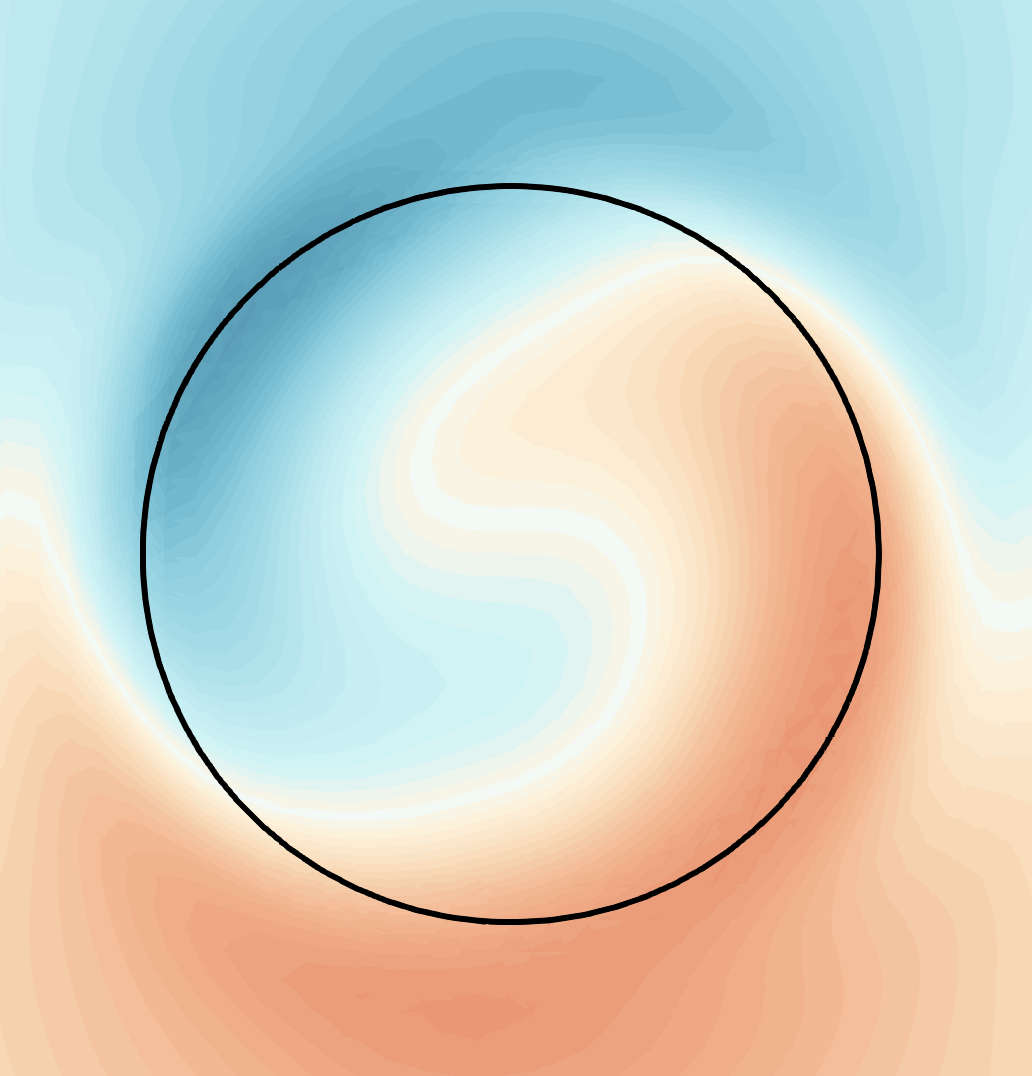}
  \end{minipage}
  \includegraphics[width=1.8cm]{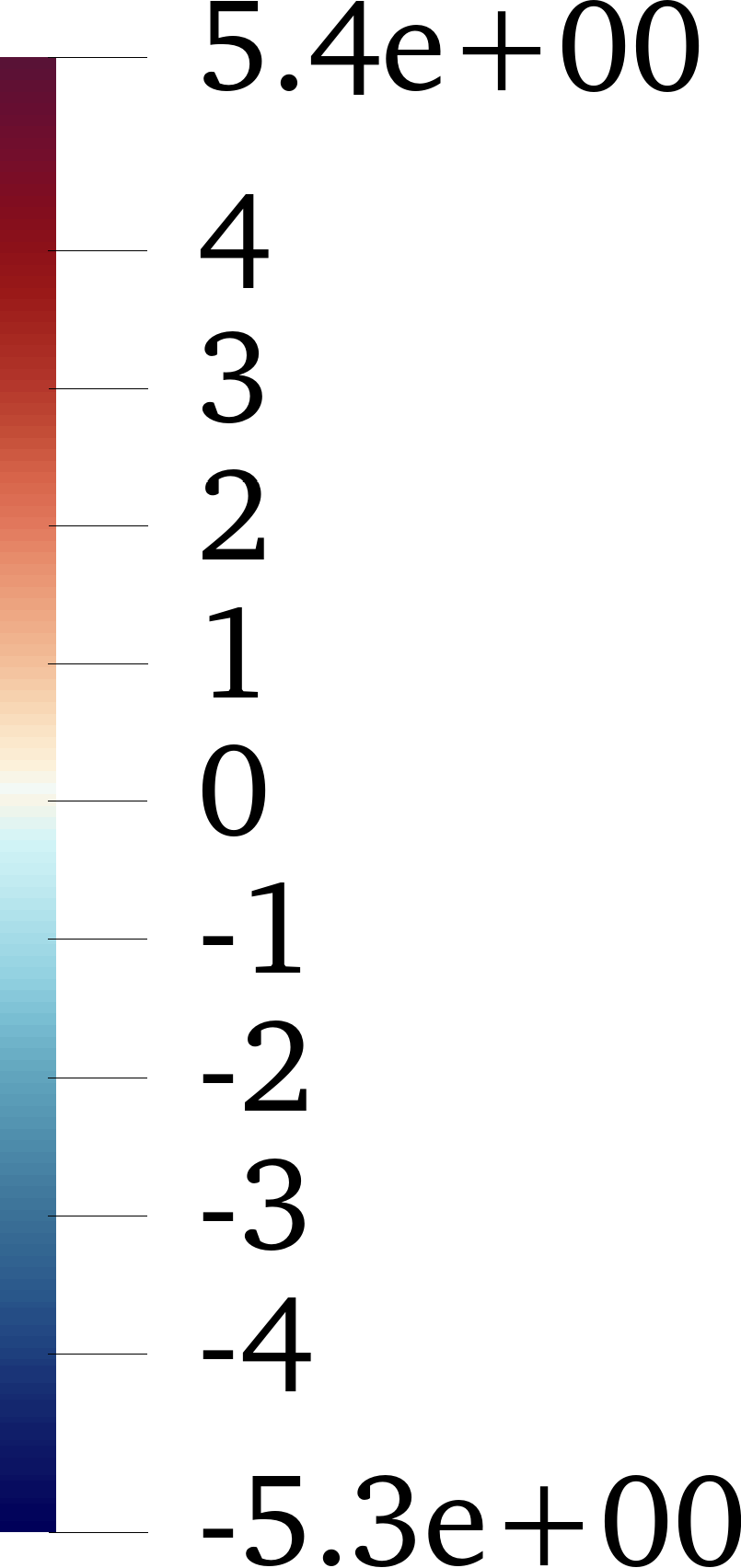}
  \caption{\hyperref[sec.numerics:subsec.ex2]{Example 2} -- Vorticity solution near the moving cylinder at times $t=0.0798, 0.0812, 0.0826$. Computed using isoparametric Taylor-Hood elements with $k=2$, $h=0.1$, $\Delta t=0.0014$, $\nu=0.01$, $\gpv=\gpp=0.1$, $\ggd=0$. Top: Narrow-band stab.; Bottom: Global stab.}
  \label{fig.ex2.picture.vort}
\end{figure}

Finally, we examine the solution visually. In \Cref{fig.ex2.picture.pre,fig.ex2.picture.div,fig.ex2.picture.vort}, we present the pressure, divergence, and vorticity solutions, respectively, near the moving cylinder at three time steps where spurious oscillations were observed in the pressure drag coefficient for the narrow-band stabilization. For the narrow-band case, the pressure solution appears relatively smooth; however, significant errors in the divergence and vorticity are evident in the extension strip. At the second time step, noticeable changes occur in the pressure solution near the top of the moving cylinder, where new pressure elements have become active. By the third time step, no new pressure elements are activated, and the pressure solution appears smoother once again.

In contrast, for the global stabilization choice shown in \Cref{fig.ex2.picture.pre,fig.ex2.picture.div,fig.ex2.picture.vort}, both pressure and vorticity are smooth within the physical domain and extend smoothly into the exterior of the fluid domain. No significant changes are observed between the three time steps. Additionally, in \Cref{fig.ex2.picture.div}, we see that with global stabilization, the velocity divergence is no longer significantly polluted within the fluid domain, particularly in the strip around the moving interface.

\begin{figure}
  \centering
  \begin{minipage}[b]{0.325\textwidth}
    \centering
      Pressure\\[4pt]
      \includegraphics[width=4.8cm, trim={0 2cm 0 1.5cm}, clip]{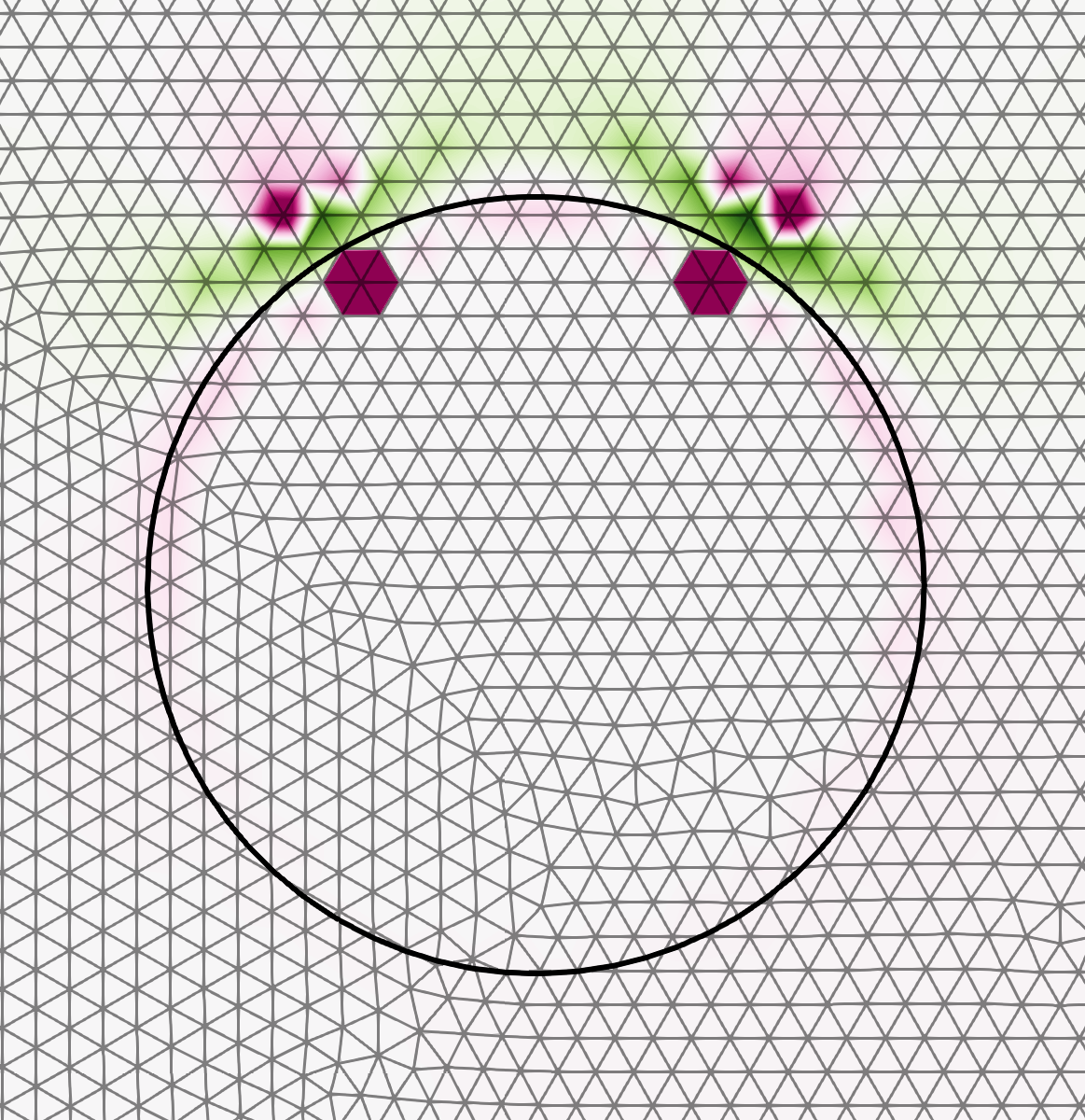}\\[8pt]
      \includegraphics[width=5cm]{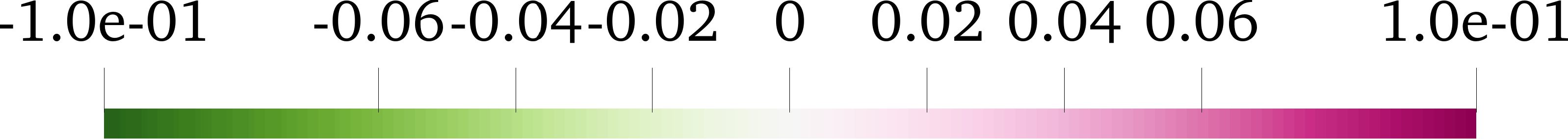}
  \end{minipage}
  \hfill
  \begin{minipage}[b]{0.325\textwidth}
    \centering
    Divergence\\[4pt]
    \includegraphics[width=4.8cm, trim={0 2cm 0 1.5cm}, clip]{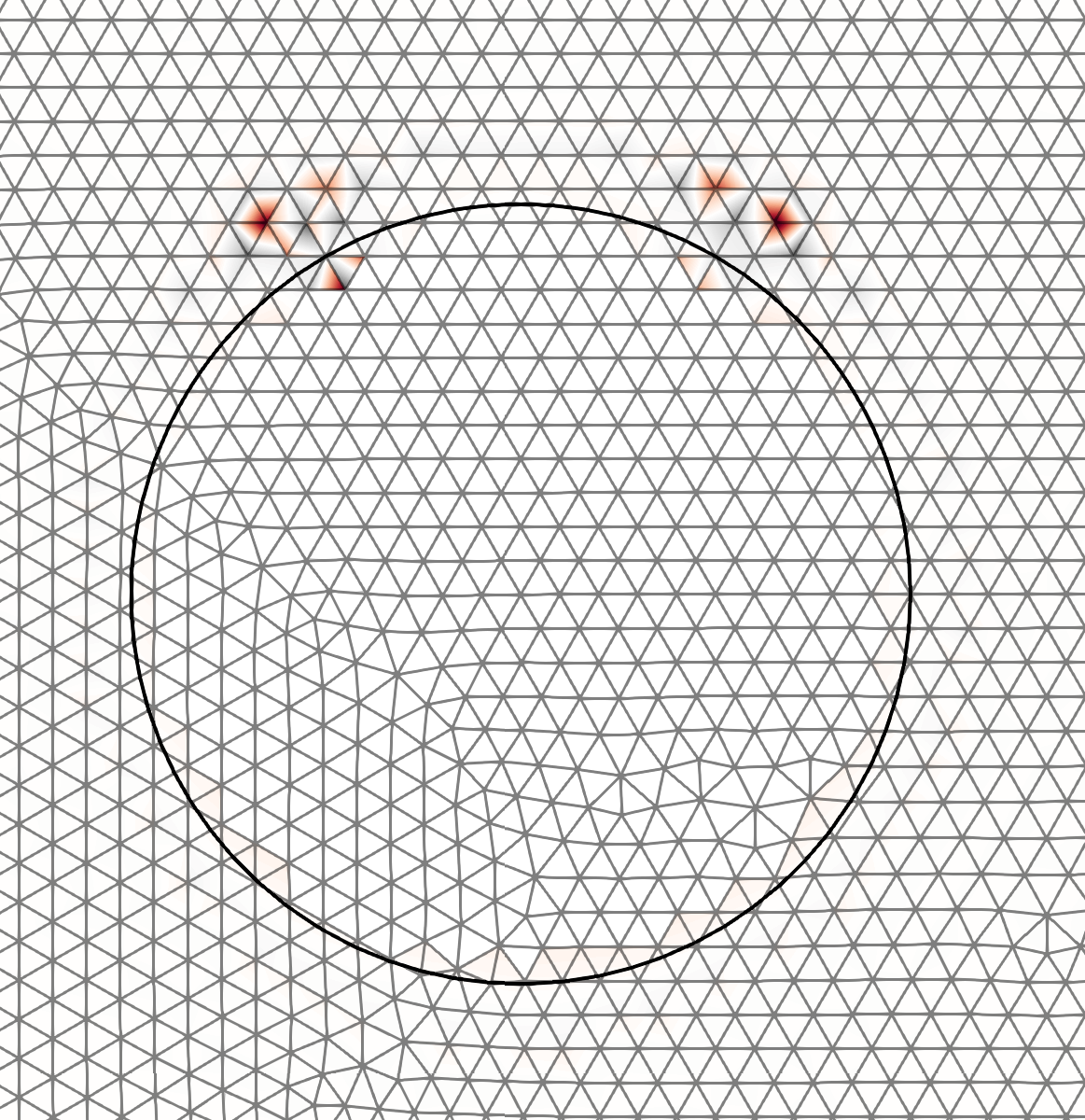}\\[8pt]
    \includegraphics[width=5cm]{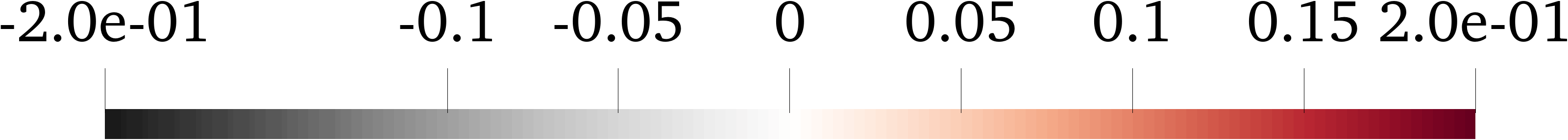}
  \end{minipage}
  \hfill
  \begin{minipage}[b]{0.325\textwidth}
    \centering
    Vorticity\\[4pt]
    \includegraphics[width=4.8cm, trim={0 2cm 0 1.5cm}, clip]{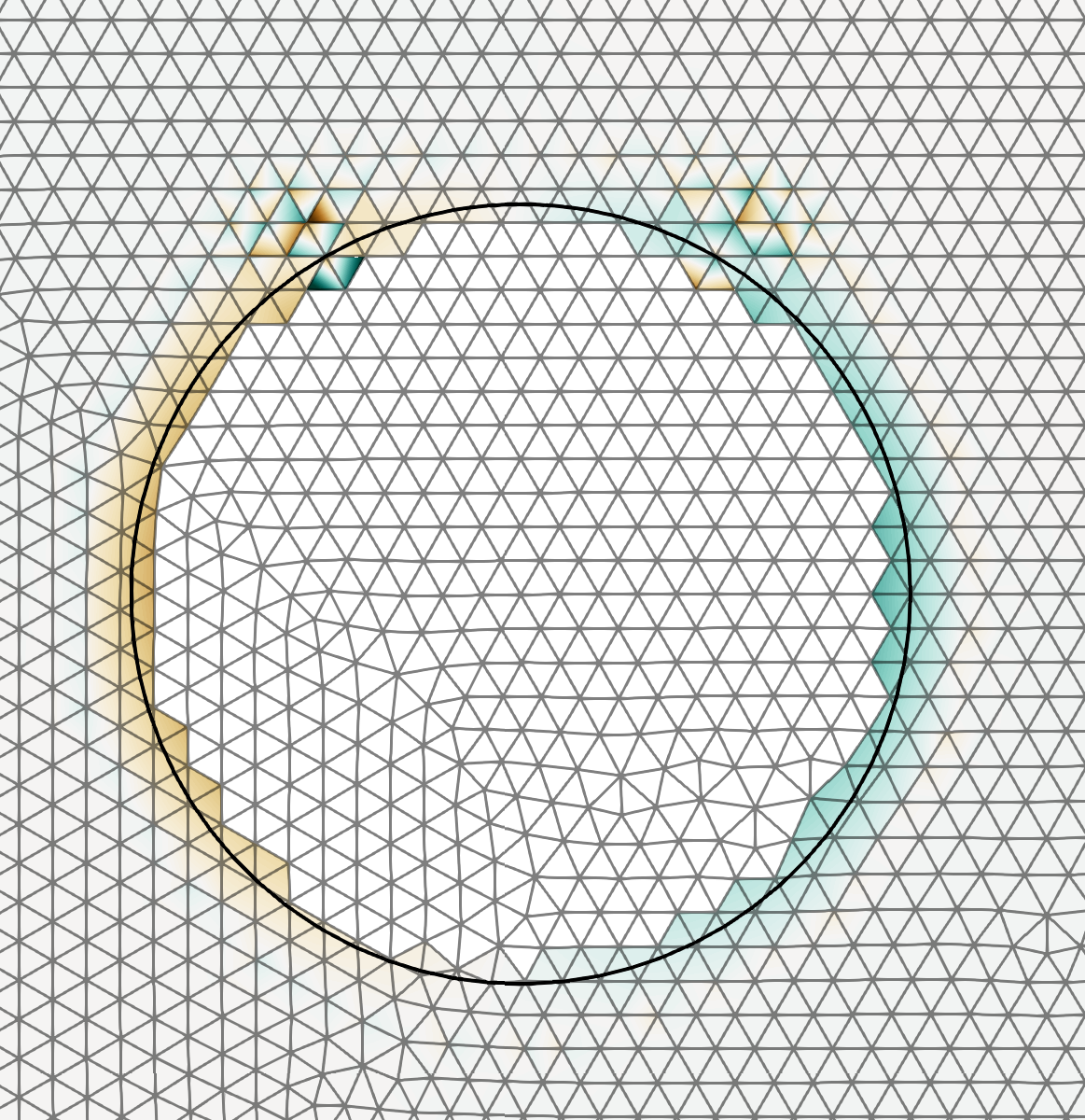}\\[8pt]
    \includegraphics[width=5cm]{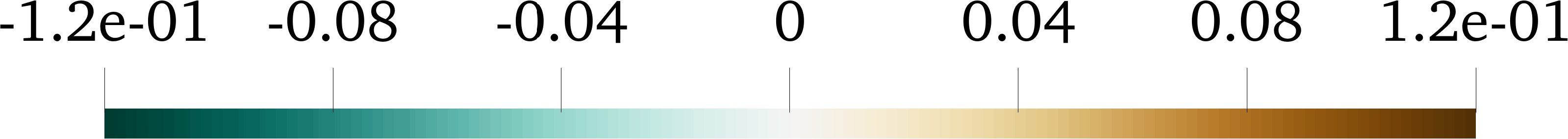}
  \end{minipage}
  \caption{\hyperref[sec.numerics:subsec.ex2]{Example 2} -- Difference between the solution at $t=0.0798$ and $0.0812$. Computed with isoparametric Taylor-Hood elements with $k=2$, $h=0.1$, $\Delta t=0.0014$, $\nu=0.01$, $\gpv=\gpp=0.1$, $\ggd=0$ .}
  \label{fig.ex2.picture_diff}
\end{figure}

To further analyze the change in the solution for the narrow-band stabilization between the first two time steps, we examine the difference between these solutions in \Cref{fig.ex2.picture_diff}. Here, the pressure oscillations are clearly \emph{localized} in space, with the velocity divergence also affected in this region. 

\subsection{Example 3: Moving square in a cross flow}
\label{sec.numerics:subsec.ex3}
As a more challenging test, we consider a modification of the previous example by replacing the cylinder with a square of side length 2, following \cite[Section 5.6]{XCM24}. Introducing corners into the geometry reduces the regularity of the exact solution. The moving geometry is described by two level set functions:
\begin{equation*} 
\phi_1(t) \coloneqq 1 - \vert \xb_1 + 5\vert \quad\text{and}\quad \phi_2(t) \coloneqq 1 - \vert \xb_2 - x_c(t)\vert 
\end{equation*}
with $x_c(t)$ as in \eqref{eqn.ex2.lset}. The domain is defined as the complement of the region where both level set functions are positive. Boundary condition on the moving boundary remains as specified in \eqref{eqn.ex2.bc}. The problem is computed for the time interval $[0, 10]$ and the results are again analyzed in the interval $\mathbb{T} = [1,9]$.

Since the exact zero level-sets are linear in this example, the resulting quadrature is exact, eliminating any geometry approximation error. Thus, any numerical fluctuations in the total fluid volume are absent. Although previous results suggested that geometry approximation error does not significantly affect the spurious pressure boundary force, here it is entirely excluded.

\subsubsection*{Results}

We again consider a series of time steps $\Delta t \in \{0.14, 0.014, 0.0014\}$, resulting in decreasing values of $C_\text{CFL}$. For each time step, we analyze the effects of the stabilization parameters and spatial resolution. The results are presented in \Cref{tab.ex3.th}. From these results, we observe the following:

\begin{itemize}
\item This example is more challenging, with larger spurious pressure oscillations overall, compared to \hyperref[sec.numerics:subsec.ex2]{Example 2}. The observed spurious oscillations here, despite exact geometry handling, reinforce the conclusion that geometry approximation error is not a primary cause of these oscillations in the method.

\item Grad-div stabilization continues to be effective, reducing the oscillation magnitude by a factor of 3.47 for $C_\text{CFL}=0.01$. As in the previous example, the choice $\gpv = 0.01$, $\gpp = 0.1$, and $\ggd = 0.1$ produces the smallest spurious oscillations among the parameters considered. However, the improvement over the default parameter case is only by a factor of 3.4, compared to 9.1 in \hyperref[sec.numerics:subsec.ex2]{Example 2}.

\item Higher-order elements are helpful, though less so than in the previous example, likely due to the reduced regularity of the exact solution caused by the sharp corners of the geometry. On the other hand, refining the mesh (and thus increasing the CFL number) significantly reduces oscillations, as expected from the previous results.

\item The spurious oscillations for the narrow-band stabilization with default parameters are of a similar magnitude as those observed for the state-of-the-art weighted shifted boundary method applied to this configuration in the simpler Stokes setting with $\nu=1$, cf.~\cite[Fig.~31]{XCM24}.

\item Global ghost-penalty stabilization leads to a significant decrease in the spurious pressure  oscillations. The absolute values are larger than in the previous example. Nevertheless, they are relatively robust with respect to the decreasing time step with increase by a factor of 2.2 for a decrease in $C_\text{CFL}$ by a factor of 100. This compares to an increase by a factor of 18 in the narrow-band case for the same parameters.
\end{itemize}

\begin{table}
  \centering
  \caption{\hyperref[sec.numerics:subsec.ex2]{Example 3} -- Spurious pressure forces statistic $\Vert e_h\Vert_{L^\infty(\mathbb{T})}$. Taylor-Hood elements with $\gpv=\gpp=0.1$. The default choice of other parameters is  $k=2,\ggd=0, h=0.1$. Three columns in the middle show results when one of these parameters is changed.}
  \label{tab.ex3.th}
  \begin{tabular}{lc@{\hspace*{8pt}}c@{\hspace*{5pt}}c@{\hspace*{8pt}}c@{\hspace*{5pt}}c}
    \toprule
      & \multicolumn{4}{c}{Narrow-band stab.} & \multicolumn{1}{c}{Global stab.}\\
    \cmidrule{2-6}
    $\Delta t~\backslash$~Param.: & Default & $\ggd=0.1$& $h=0.05$ & $k=3$ &  Default \\
    \midrule
    0.1400 & 0.24235& 0.10097&0.40380&nan     &0.04710\\
    0.0140 & 1.71026& 0.73309&0.61956&0.61050 &0.09772\\
    0.0014 & 4.46424& 1.30030&0.72382&2.47181 &0.10378\\
    \bottomrule
  \end{tabular}
\end{table}

\begin{figure}
  \centering
  \begin{minipage}[b]{14cm}
  \centering
   $t=0.3206$
   \hspace*{78pt}
   $t=0.3220$
   \hspace*{78pt}
   $t=0.3234$\\[2pt]
  \includegraphics[width=4.5cm]{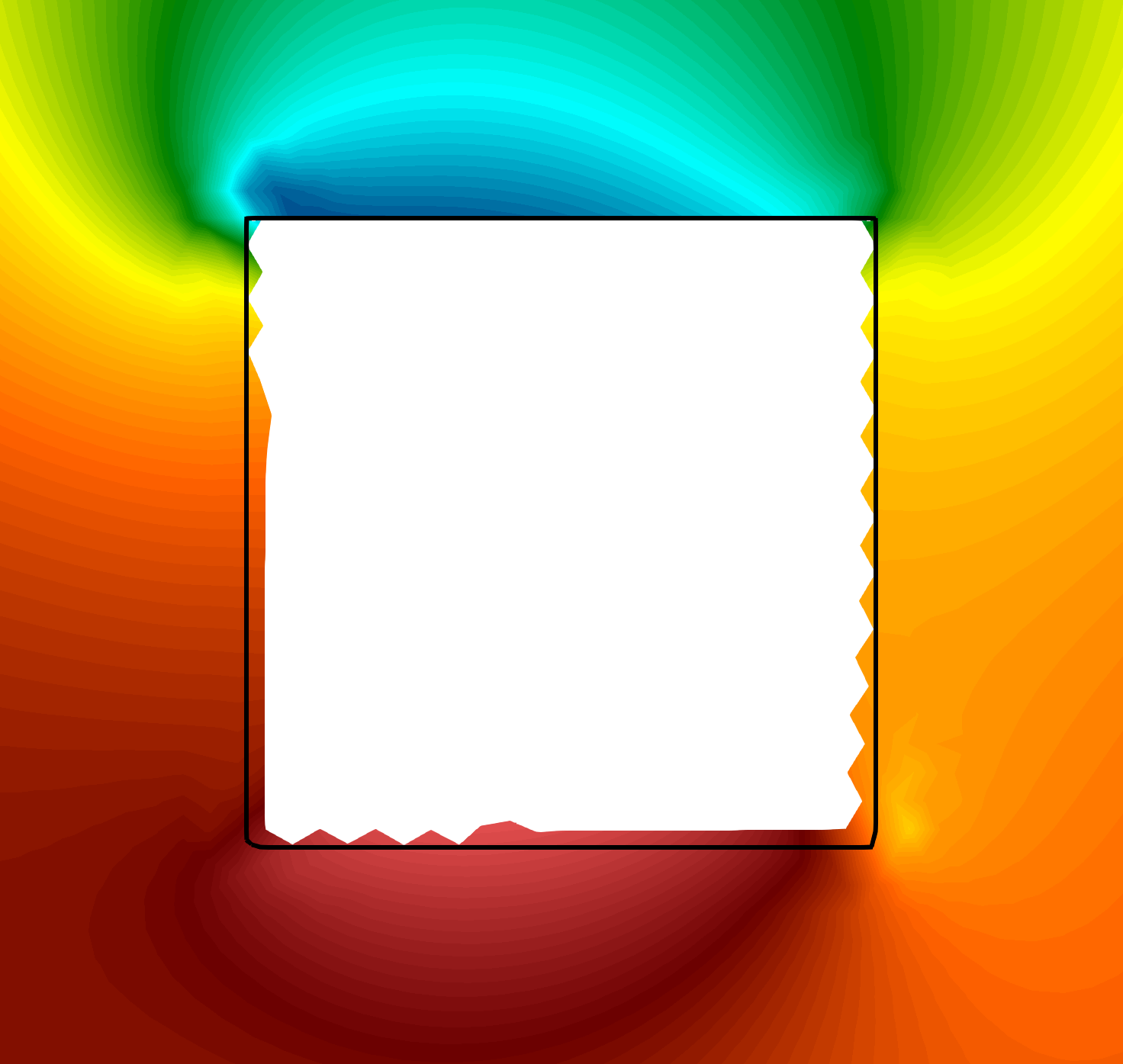}
  \includegraphics[width=4.5cm]{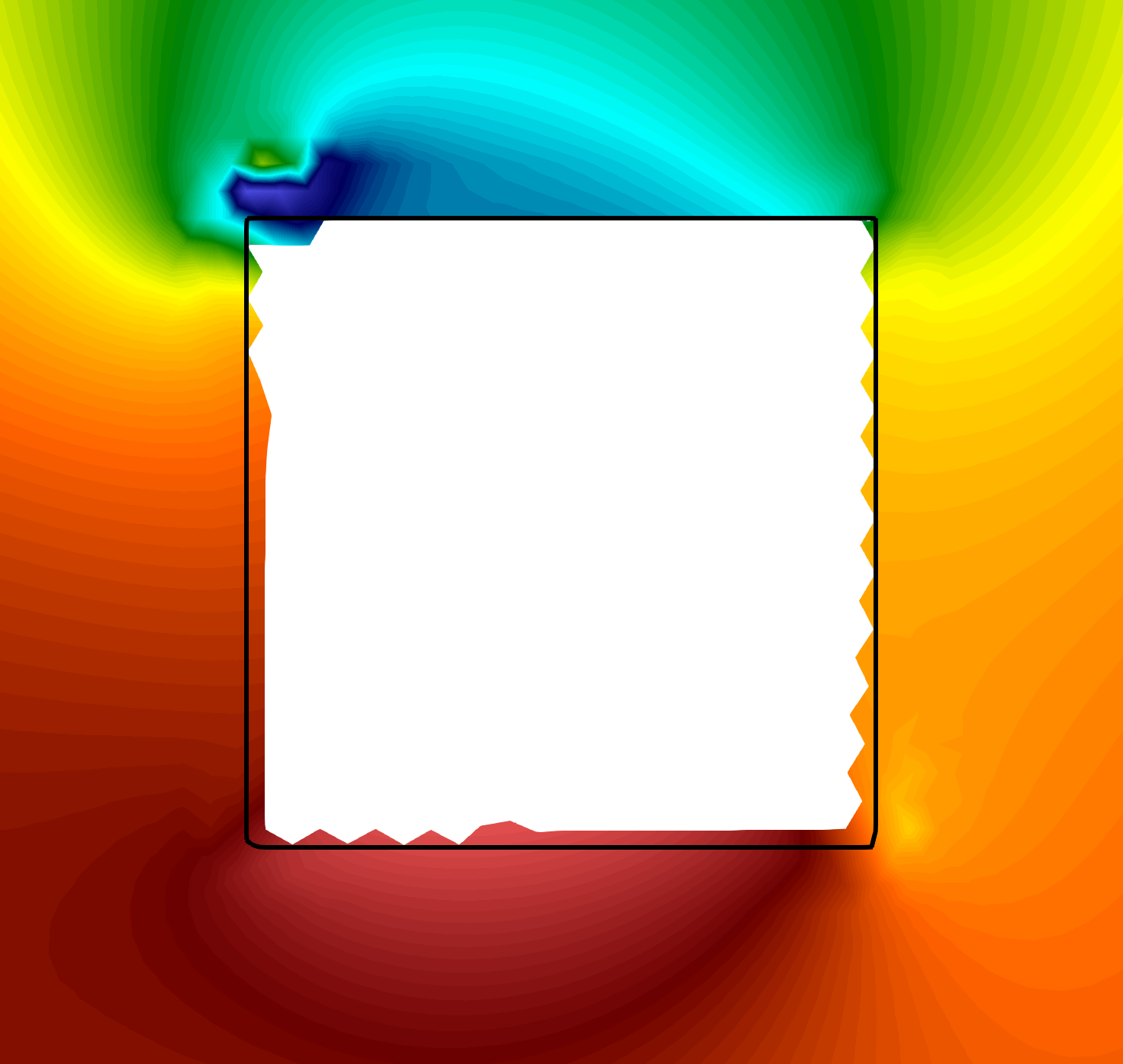}
  \includegraphics[width=4.5cm]{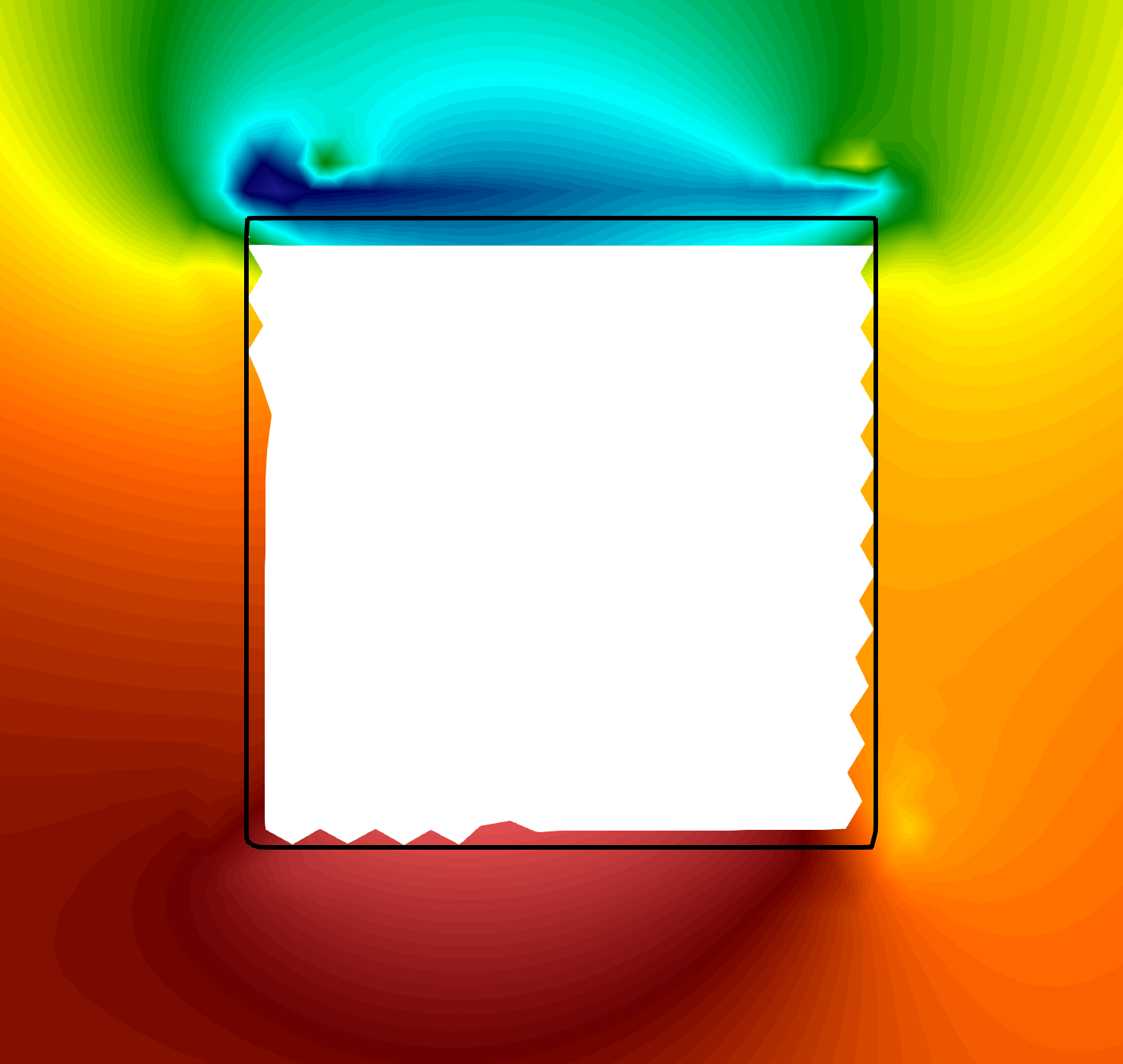}\\[2pt]
  \includegraphics[width=4.5cm]{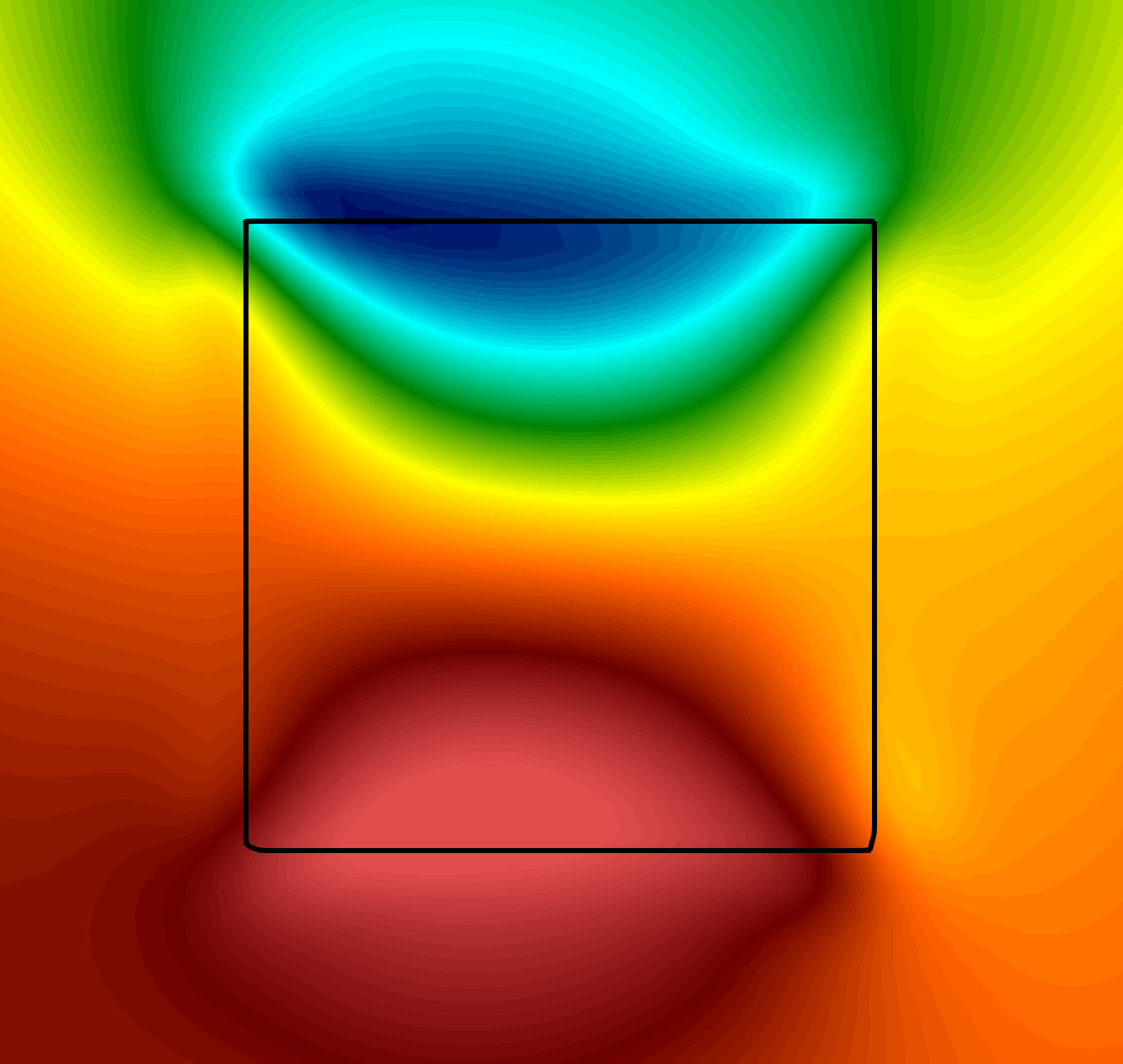}
  \includegraphics[width=4.5cm]{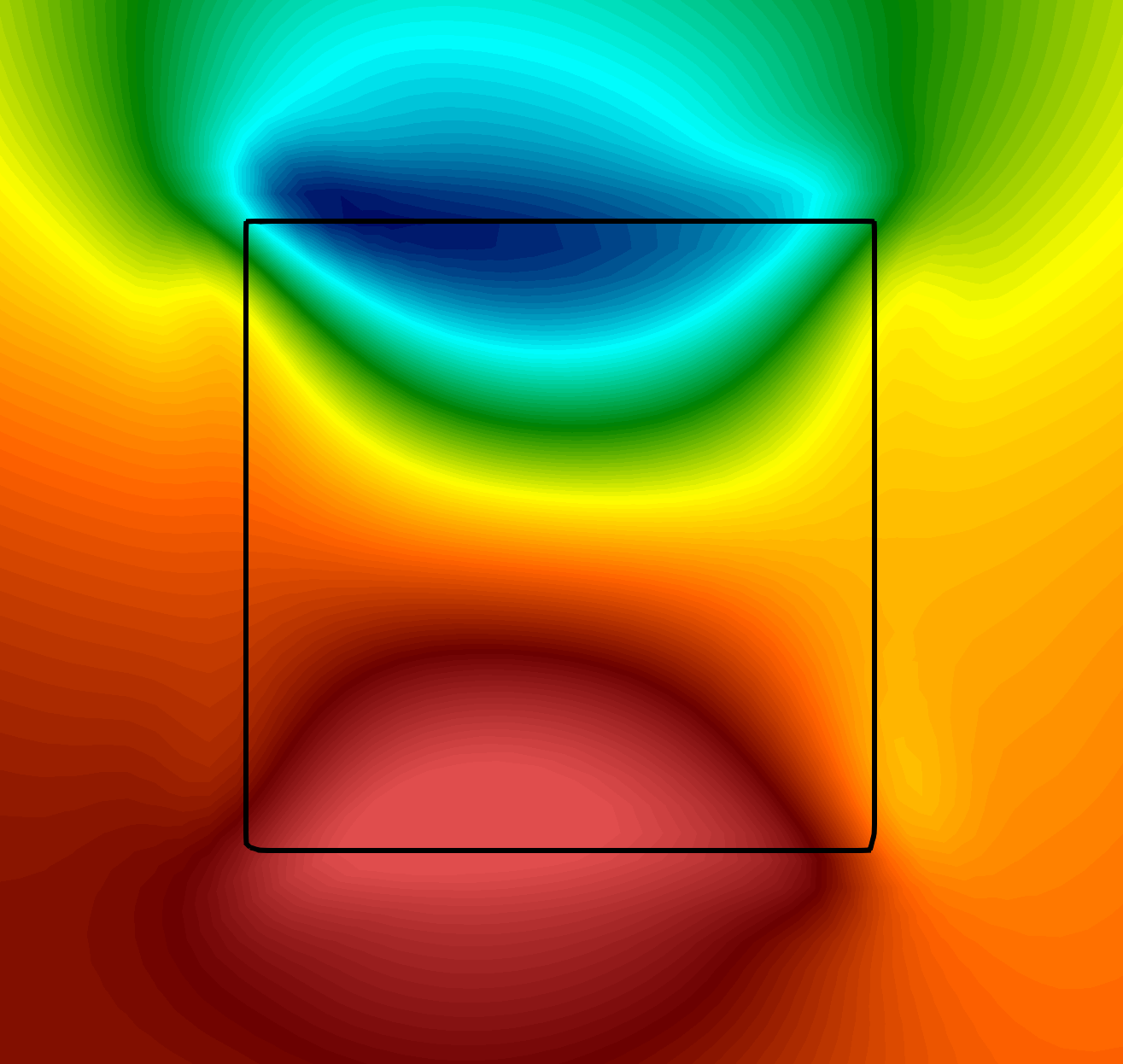}
  \includegraphics[width=4.5cm]{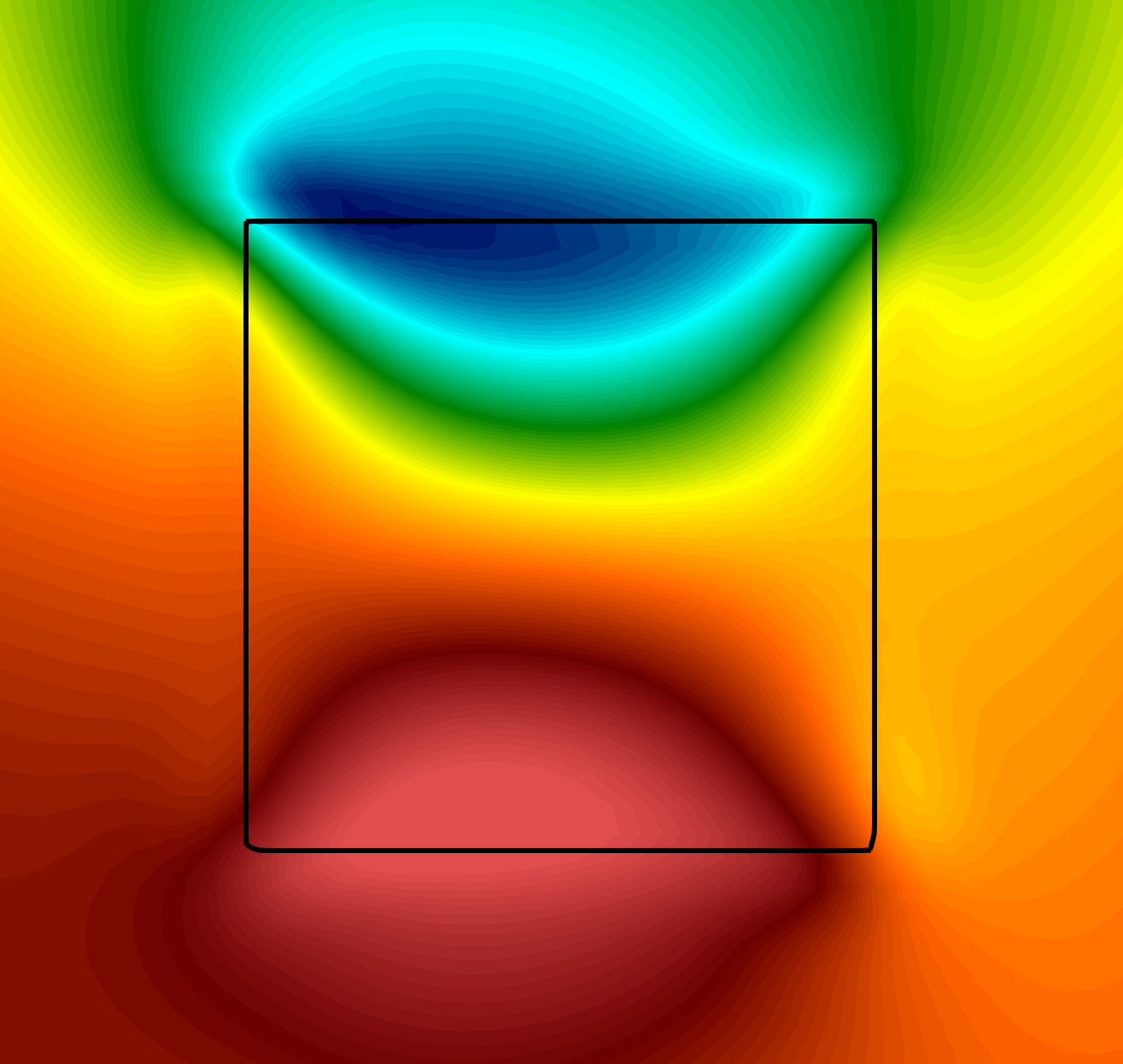}
  \end{minipage}
  \includegraphics[width=1.8cm]{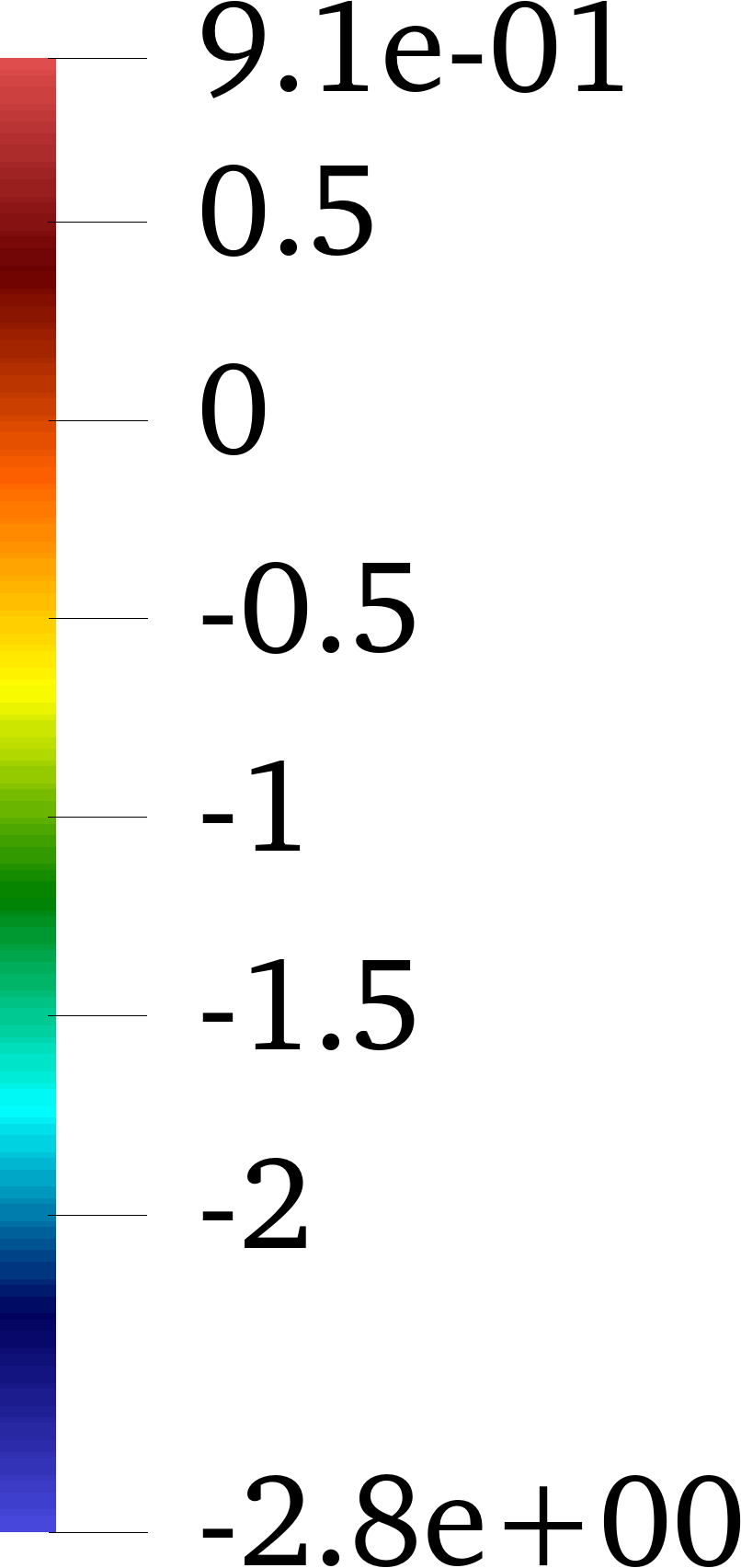}
  \caption{\hyperref[sec.numerics:subsec.ex3]{Example 3} -- Pressure solution near the moving cylinder at times $t=0.0798, 0.0812, 0.0826$. Computed with Taylor-Hood elements with $k=2$, $h=0.1$, $\Delta t=0.0014$, $\nu=0.01$, $\gpv=\gpp=0.1$, $\ggd=0$. Top: Narrow-band stab.; Bottom: Global stab.}
  \label{fig.ex3.picture.pre}
\end{figure}

\begin{figure}
  \centering
  \begin{minipage}[b]{14cm}
  \centering
   $t=0.3206$
   \hspace*{78pt}
   $t=0.3220$
   \hspace*{78pt}
   $t=0.3234$\\[2pt]
  \includegraphics[width=4.5cm]{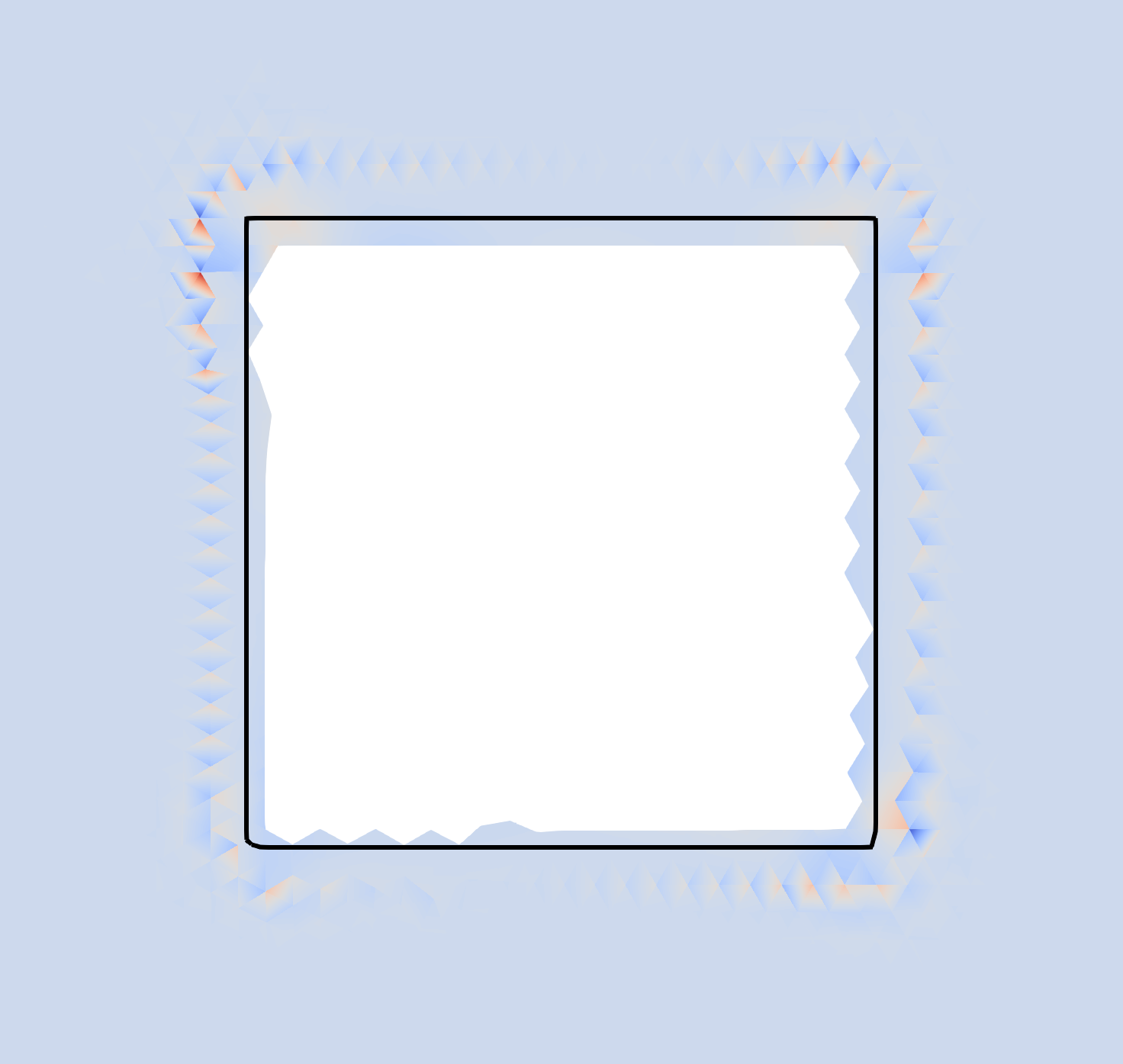}
  \includegraphics[width=4.5cm]{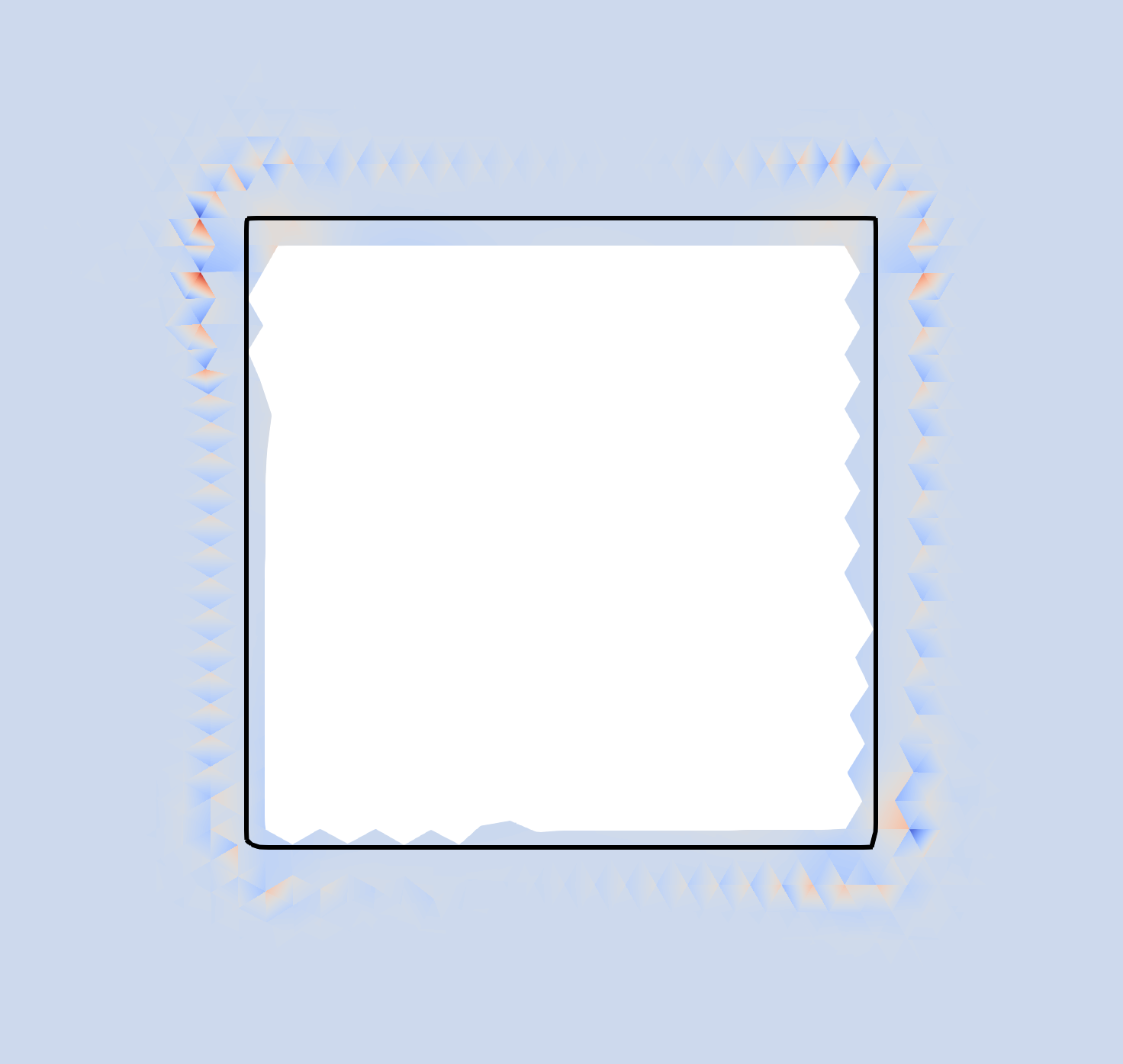}
  \includegraphics[width=4.5cm]{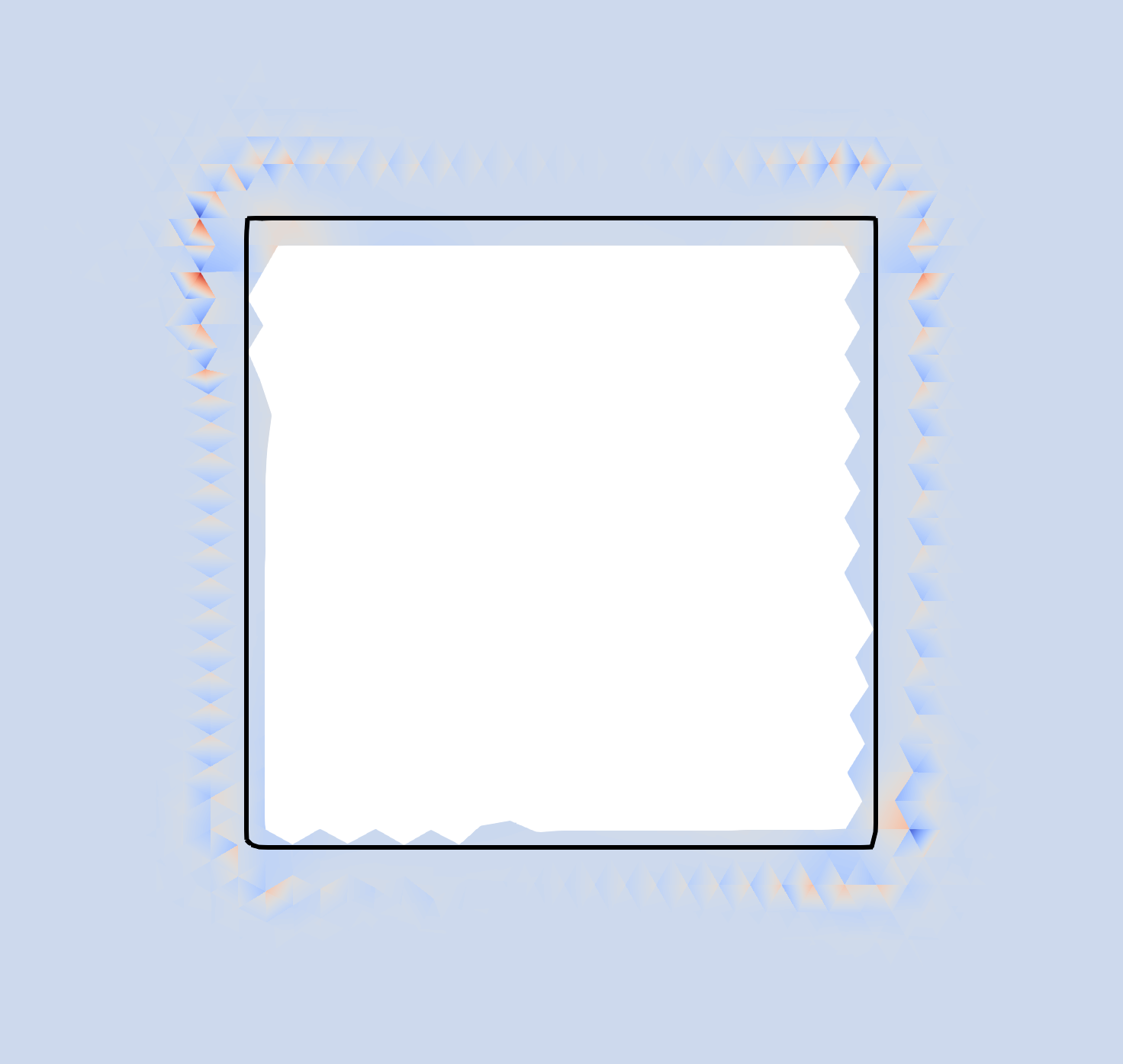}\\[2pt]
  \includegraphics[width=4.5cm]{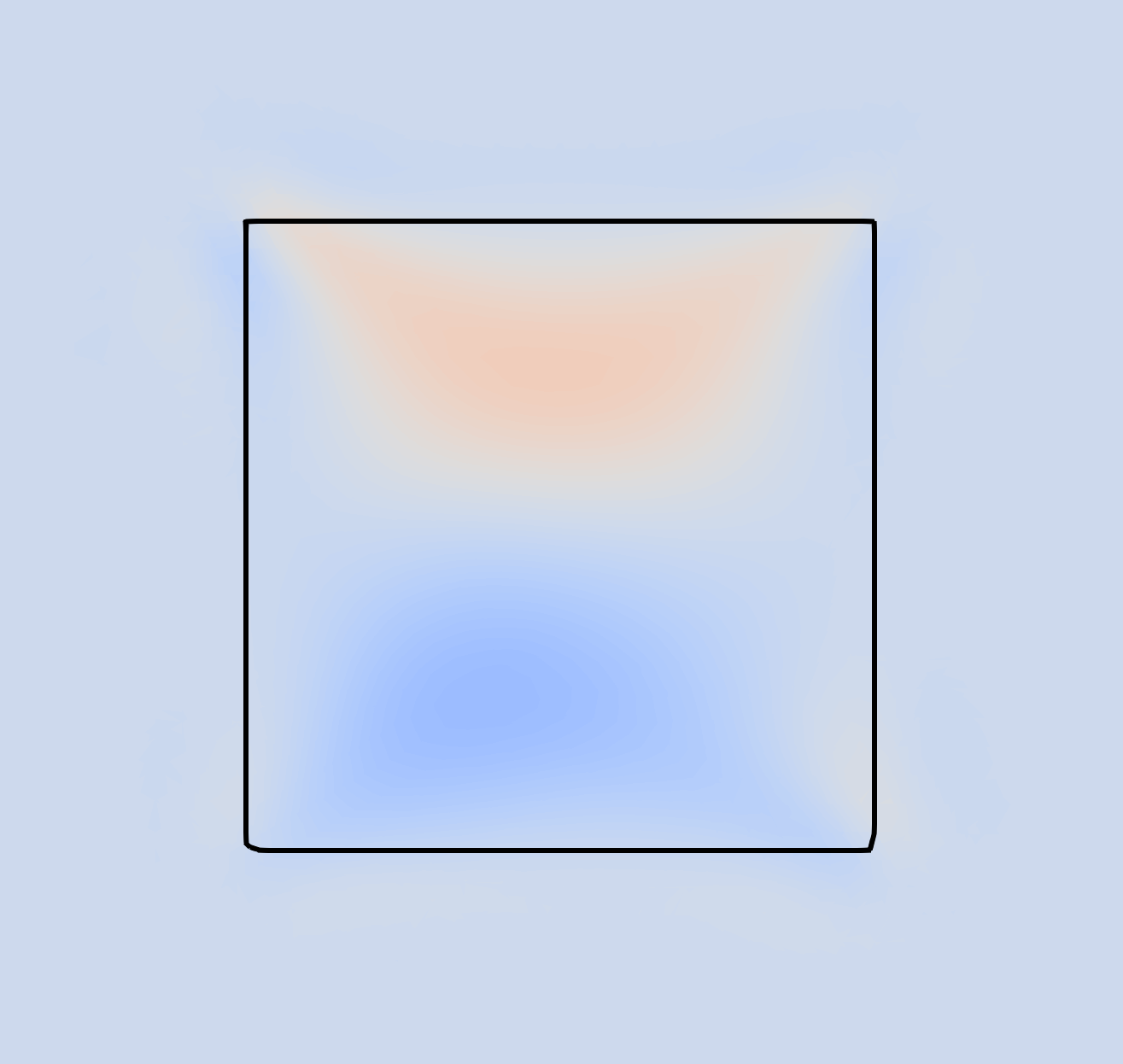}
  \includegraphics[width=4.5cm]{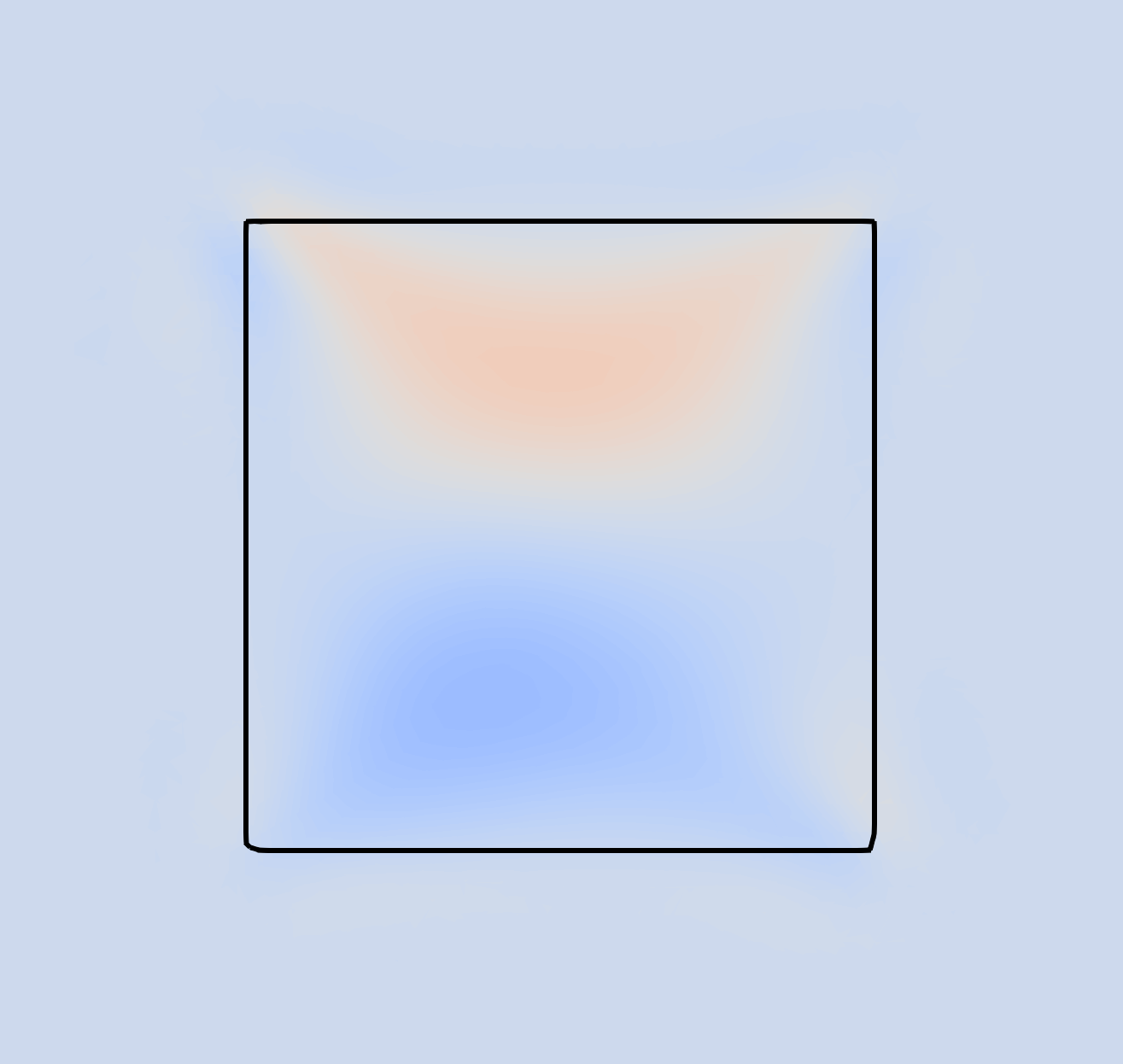}
  \includegraphics[width=4.5cm]{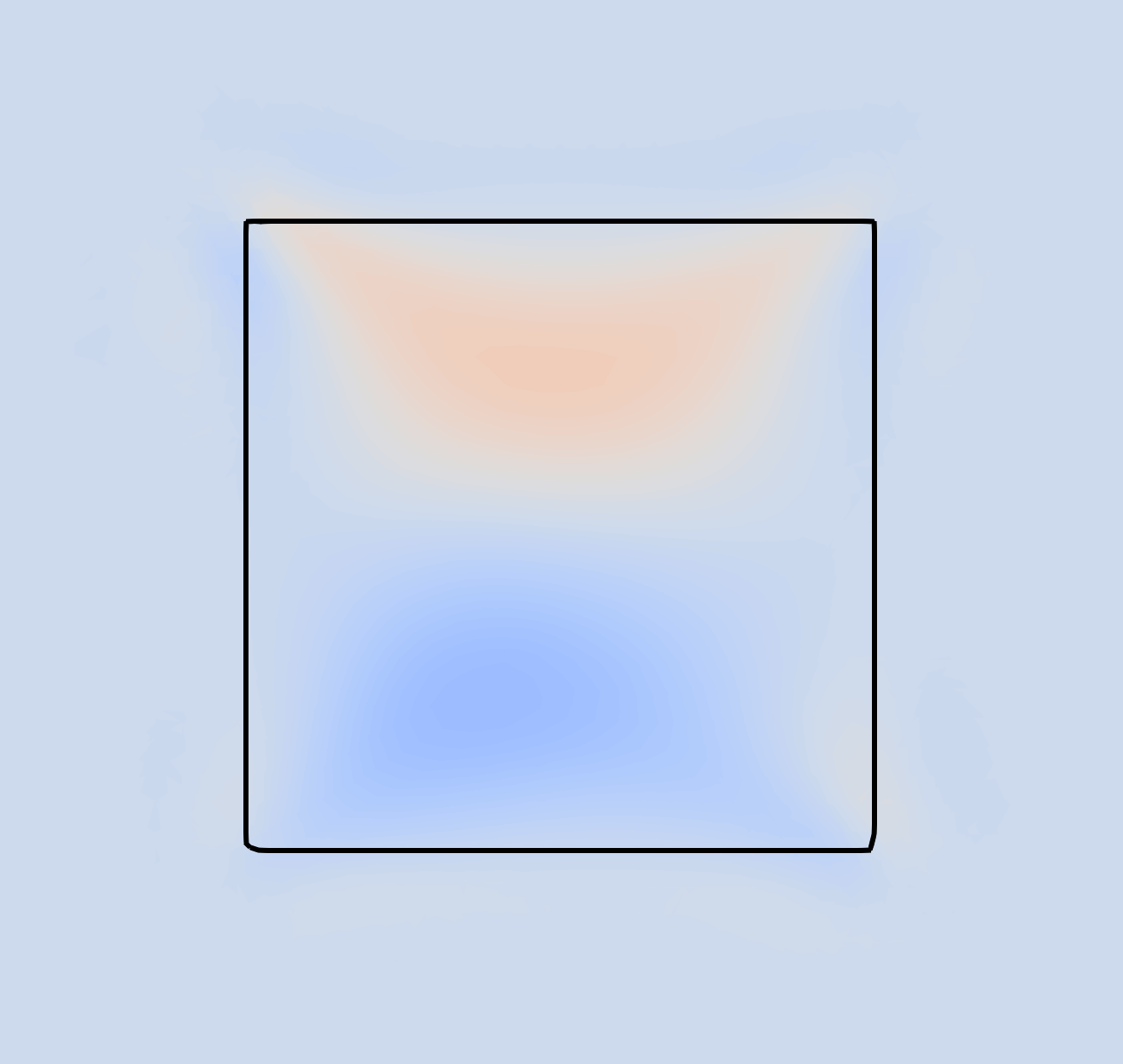}
  \end{minipage}
  \includegraphics[width=1.8cm]{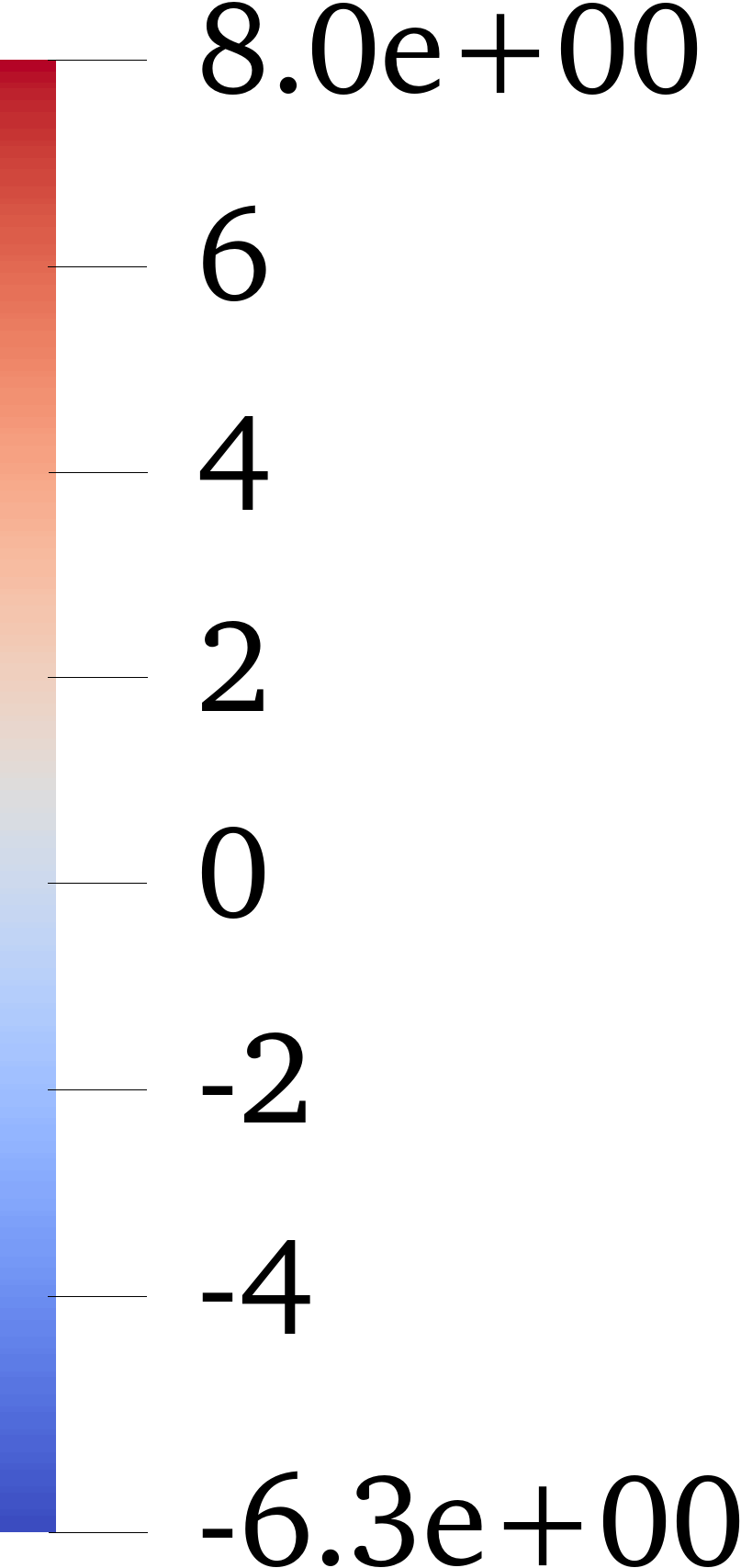}
  \caption{\hyperref[sec.numerics:subsec.ex3]{Example 3} -- Divergence of the velocity solution near the moving cylinder at times $t=0.0798, 0.0812, 0.0826$. Computed using Taylor-Hood elements with $k=2$, $h=0.1$, $\Delta t=0.0014$, $\nu=0.01$, $\gpv=\gpp=0.1$, $\ggd=0$. Top: Narrow-band stab.; Bottom: Global stab.}
  \label{fig.ex3.picture.div}
\end{figure}

\begin{figure}
  \centering
  \begin{minipage}[b]{14cm}
  \centering
   $t=0.3206$
   \hspace*{78pt}
   $t=0.3220$
   \hspace*{78pt}
   $t=0.3234$\\[2pt]
  \includegraphics[width=4.5cm]{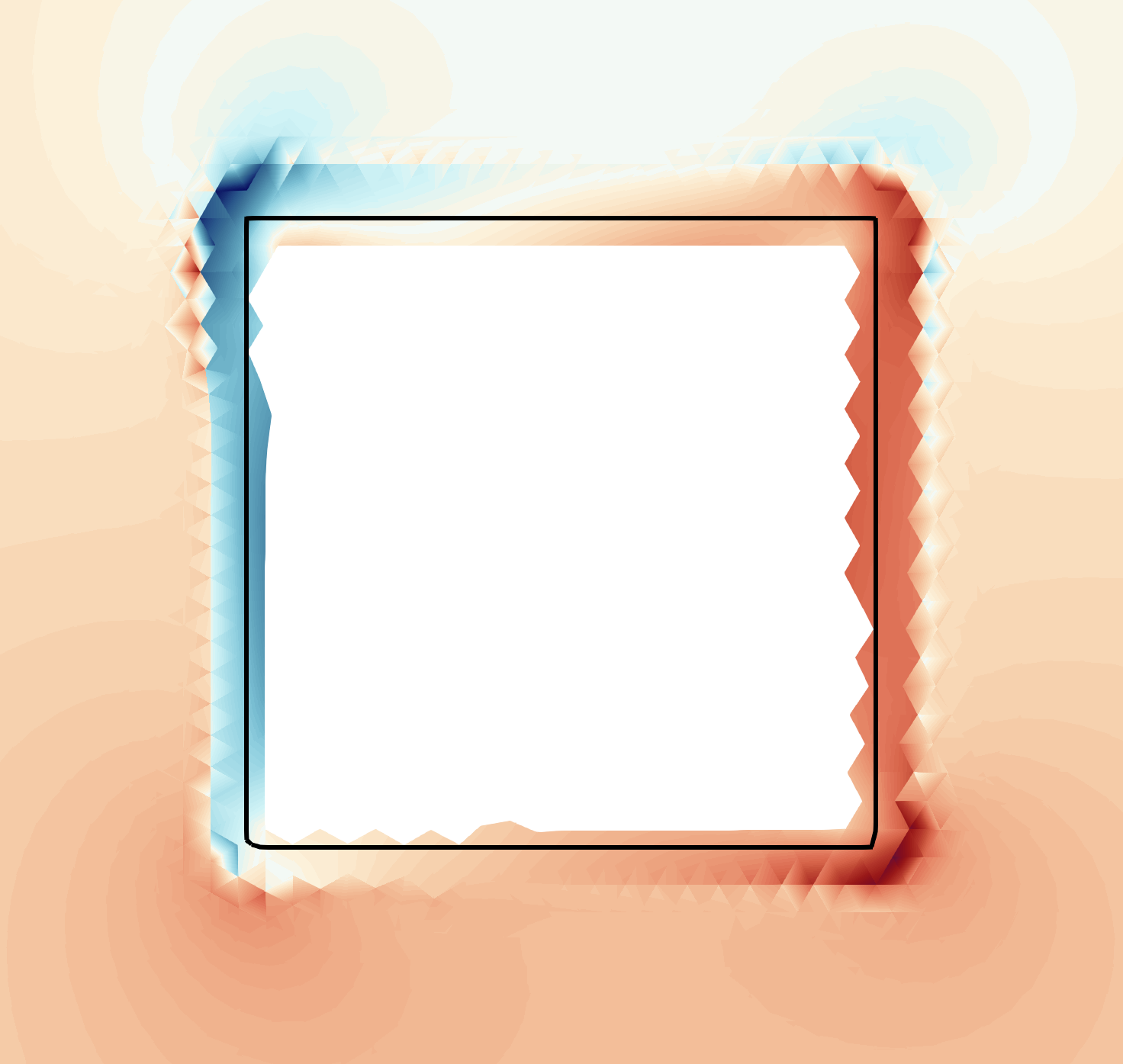}
  \includegraphics[width=4.5cm]{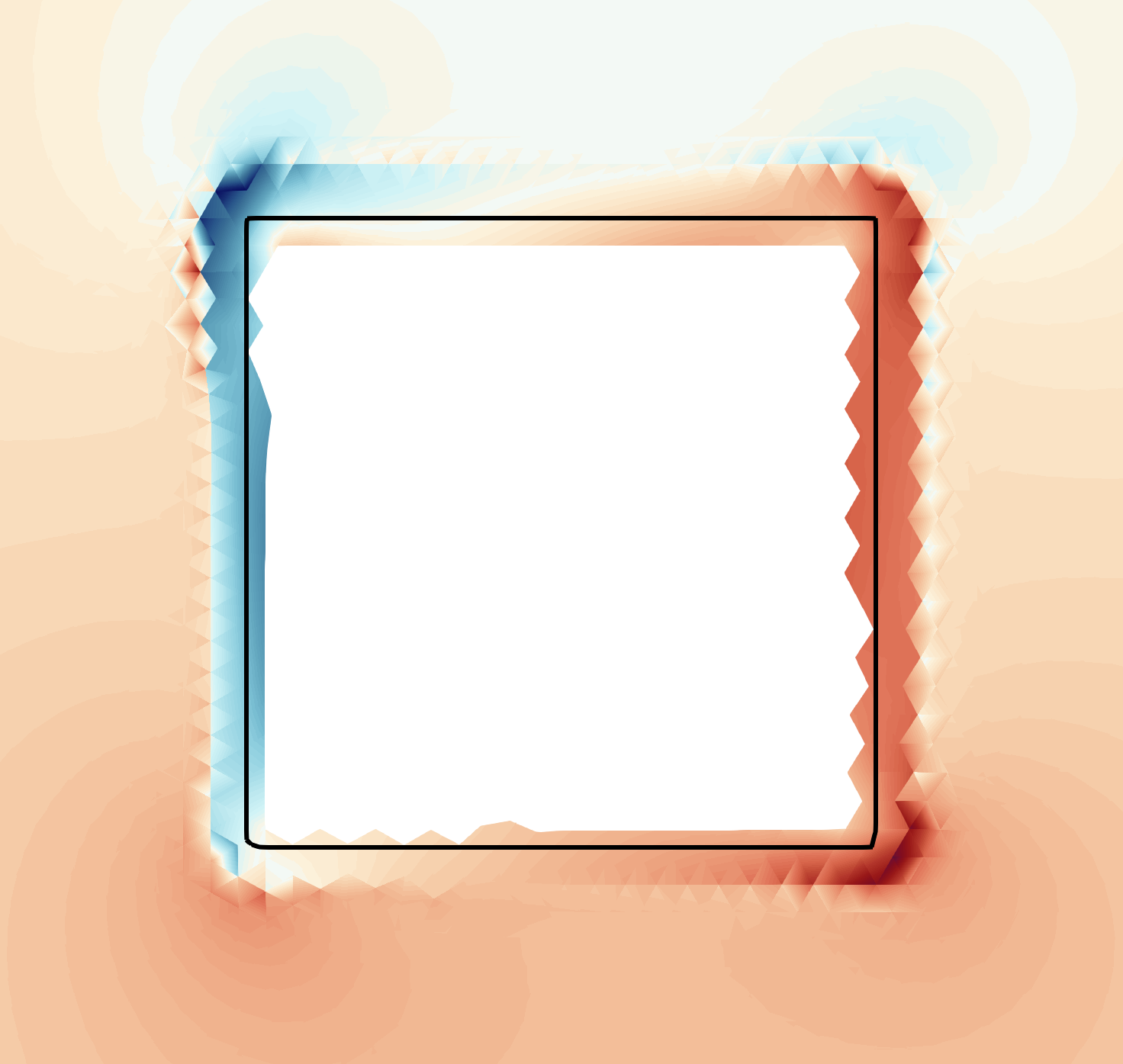}
  \includegraphics[width=4.5cm]{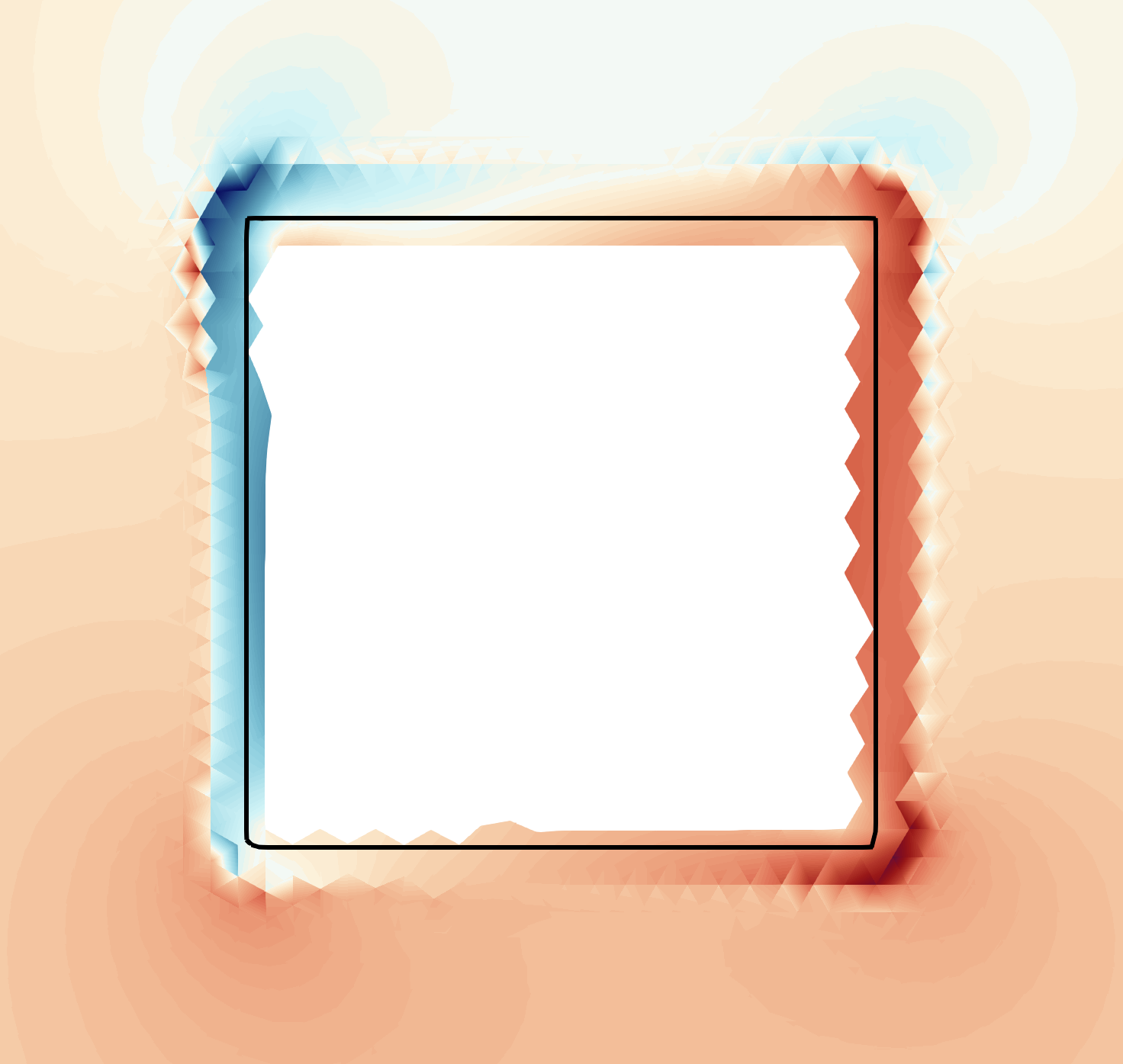}\\[2pt]
  \includegraphics[width=4.5cm]{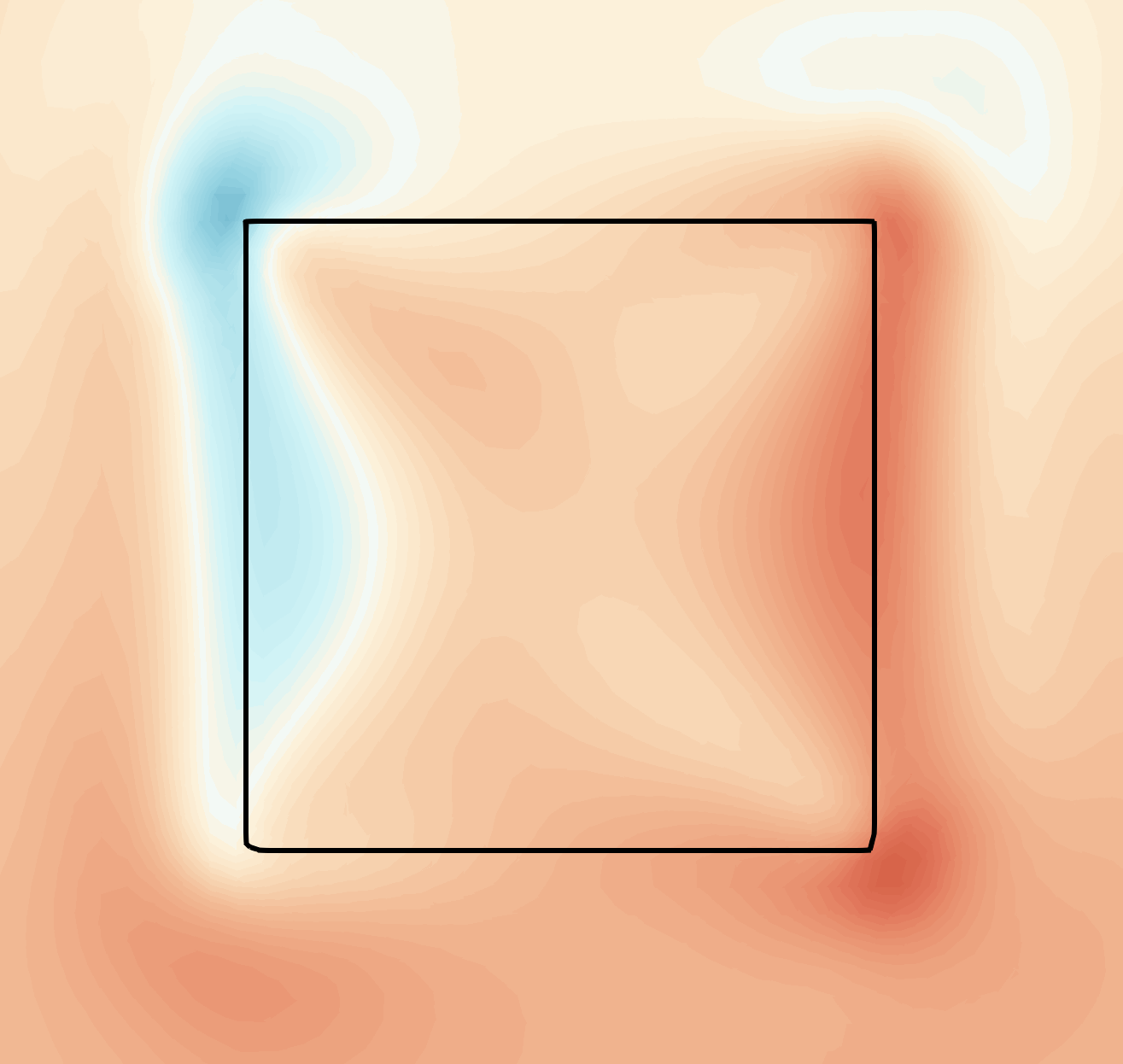}
  \includegraphics[width=4.5cm]{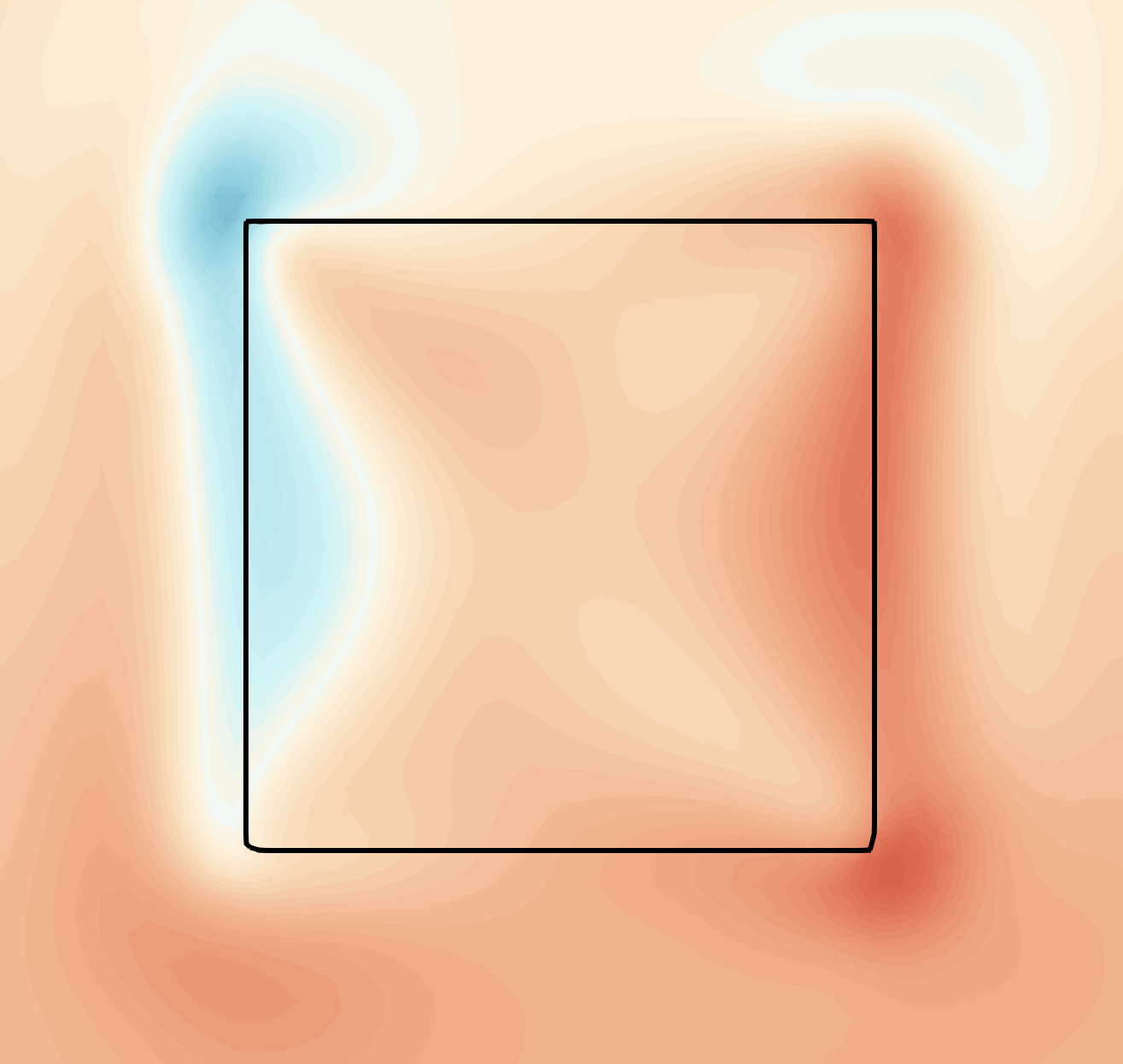}
  \includegraphics[width=4.5cm]{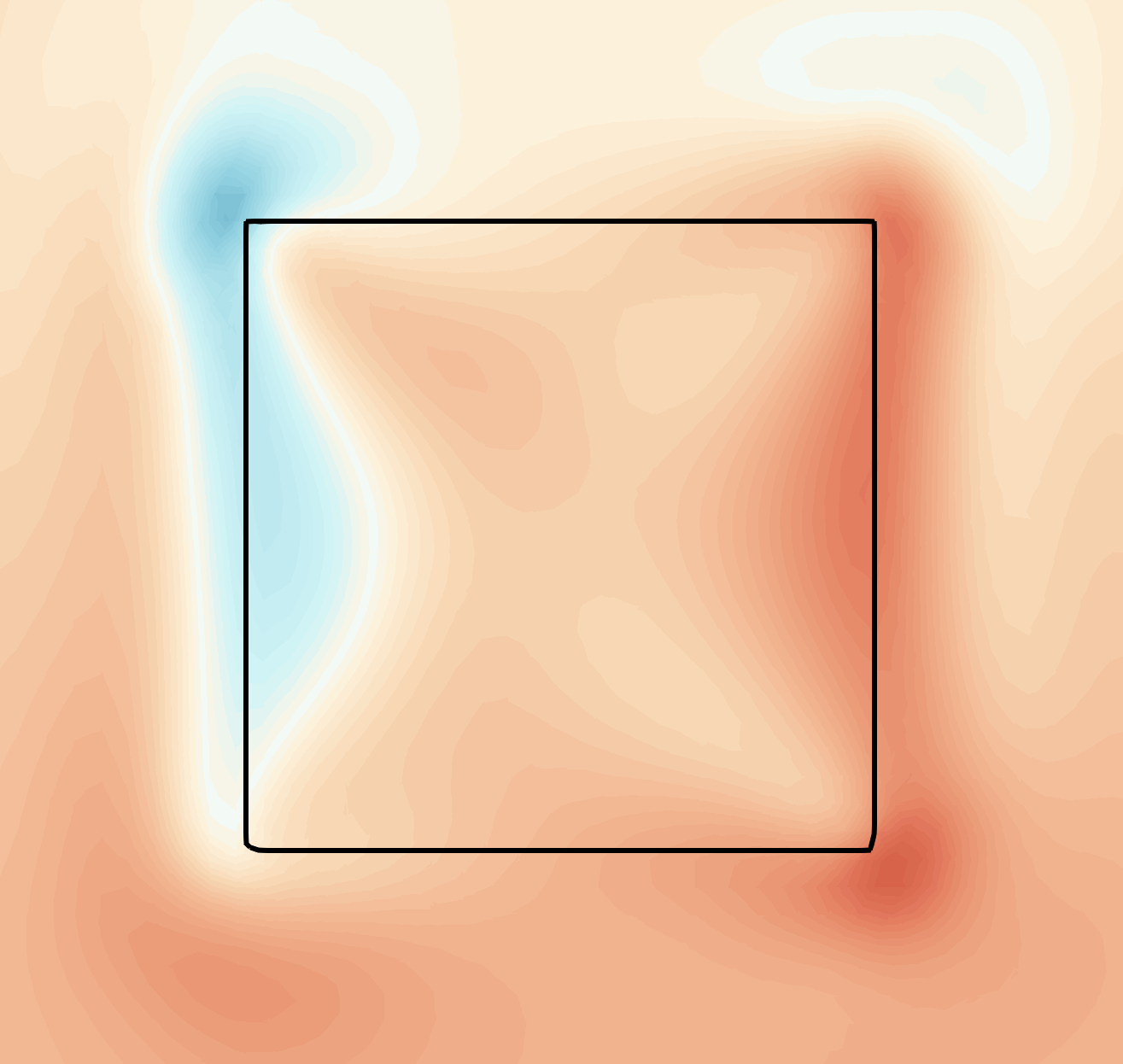}
  \end{minipage}
  \includegraphics[width=1.8cm]{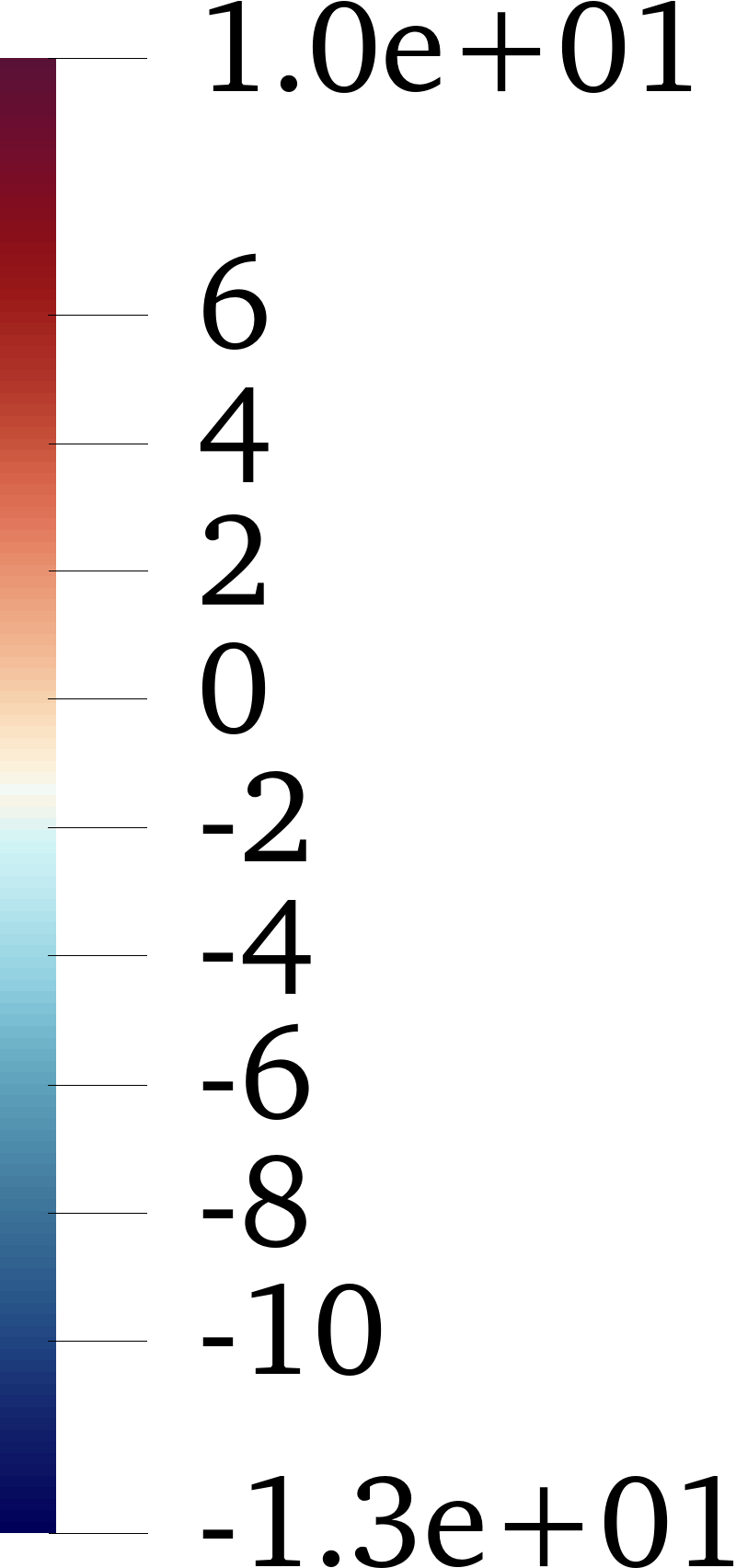}
  \caption{\hyperref[sec.numerics:subsec.ex3]{Example 3} -- Divergence of the velocity solution near the moving cylinder at times $t=0.0798, 0.0812, 0.0826$. Computed using Taylor-Hood elements with $k=2$, $h=0.1$, $\Delta t=0.0014$, $\nu=0.01$, $\gpv=\gpp=0.1$, $\ggd=0$. Top: Narrow-band stab.; Bottom: Global stab.}
  \label{fig.ex3.picture.vort}
\end{figure}

\Cref{fig.ex3.picture.pre,fig.ex3.picture.div,fig.ex3.picture.vort} show the solution near the moving boundary at three time instances, where spurious oscillations in the pressure drag coefficient were observed in the case of the narrow-band stabilization. Consistent with our quantitative results, the change in pressure here is much larger and more visible than in \Cref{fig.ex2.picture.pre}. The pressure change is, once again, most pronounced near regions where new pressure elements have recently become active.
For the global stabilization, we again do not see any spurious oscillations in the pressure and both the divergence and vorticity does not appear to be polluted in a strip around the moving interface.

\begin{figure}
  \centering
  \begin{minipage}[b]{0.325\textwidth}
    \centering
      Pressure\\[4pt]
      \includegraphics[width=4.8cm]{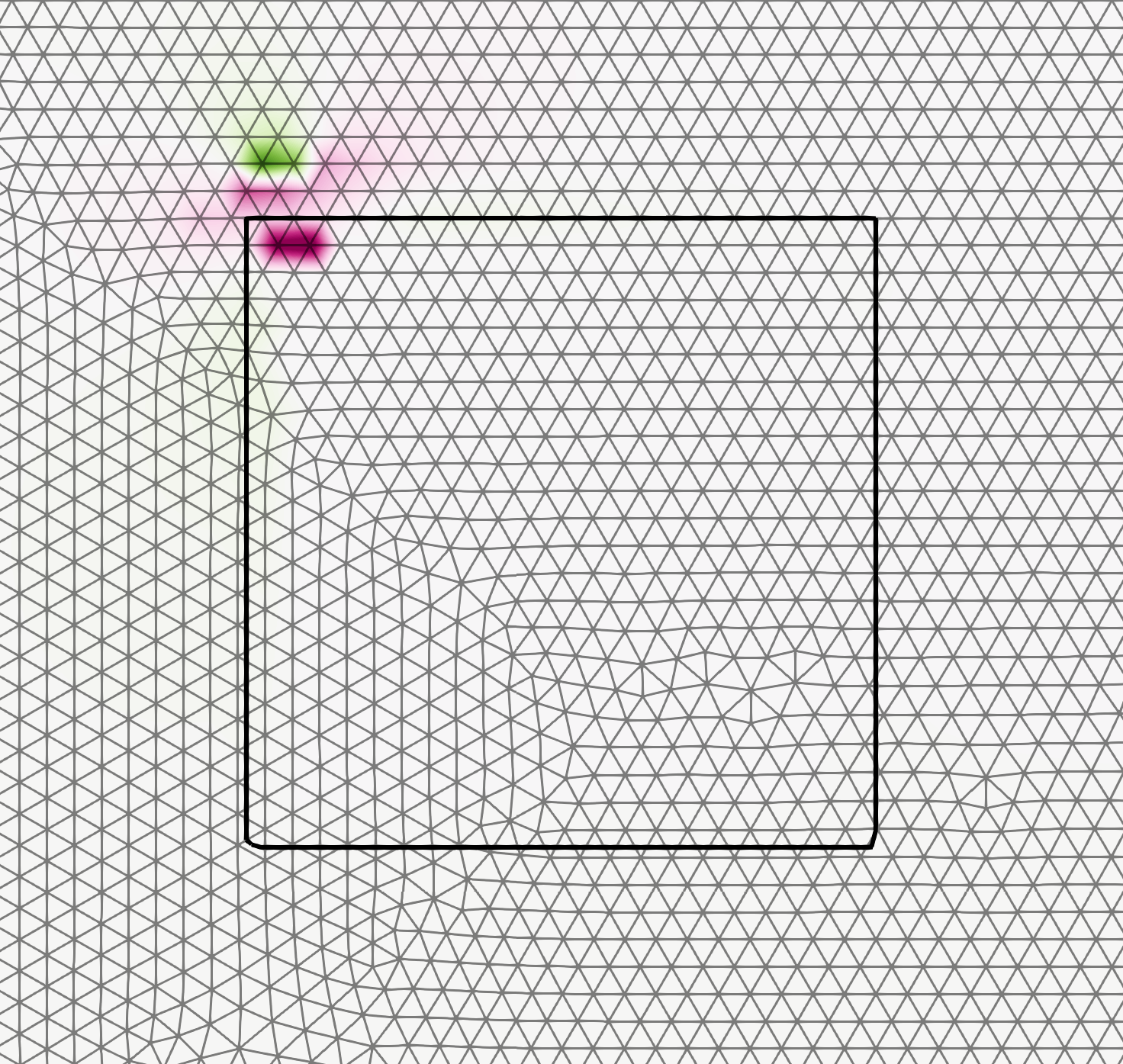}\\[8pt]
      \includegraphics[width=5cm]{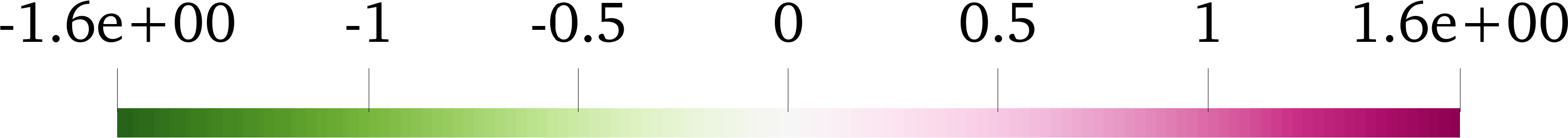}
  \end{minipage}
  \hfill
  \begin{minipage}[b]{0.325\textwidth}
    \centering
    Divergence\\[4pt]
    \includegraphics[width=4.8cm]{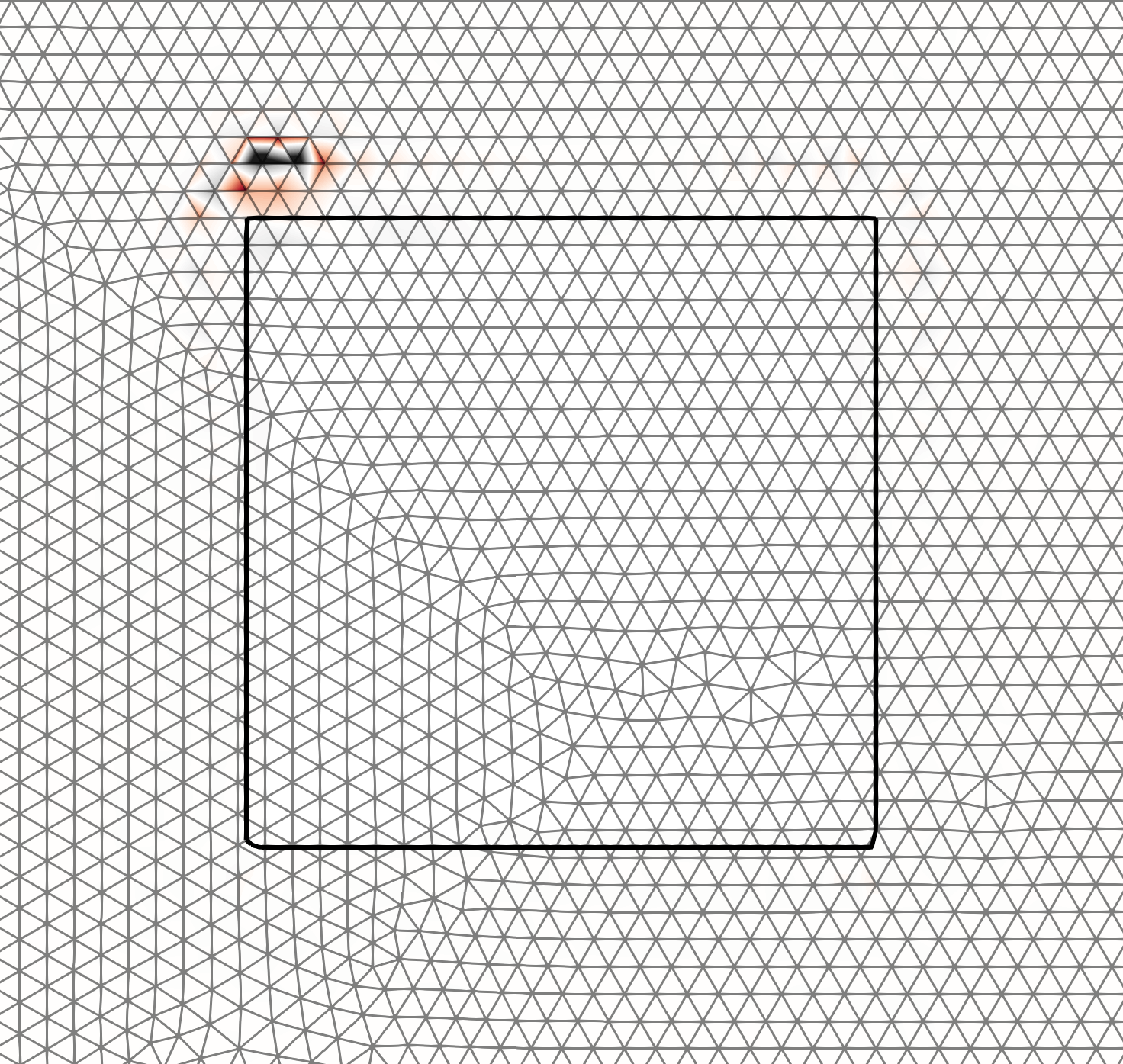}\\[8pt]
    \includegraphics[width=5cm]{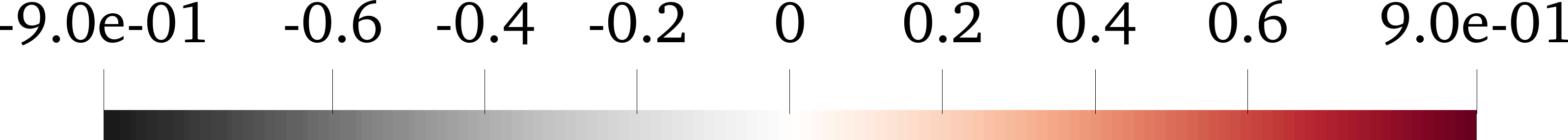}
  \end{minipage}
  \hfill
  \begin{minipage}[b]{0.325\textwidth}
    \centering
    Vorticity\\[4pt]
    \includegraphics[width=4.8cm]{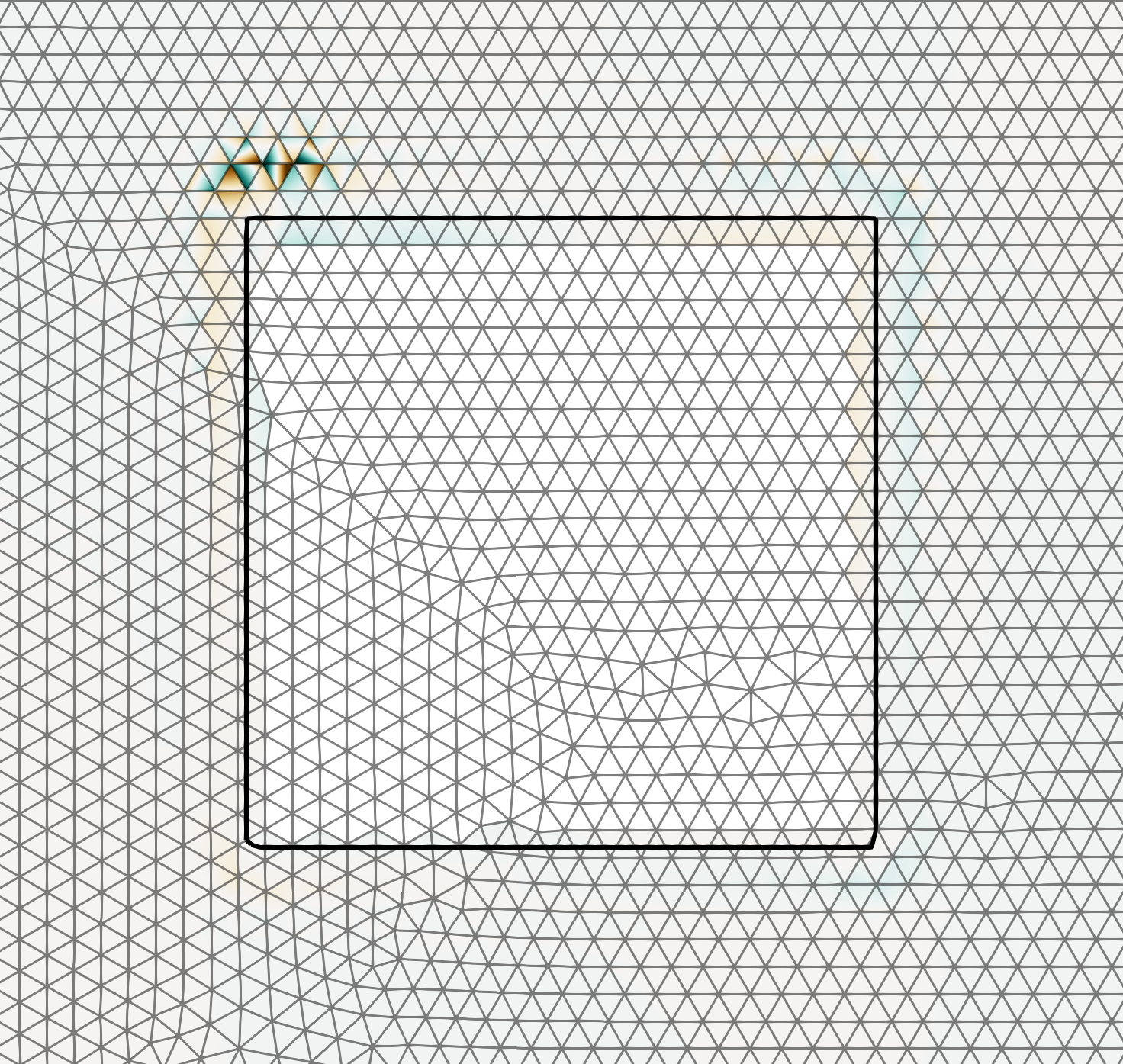}\\[8pt]
    \includegraphics[width=5cm]{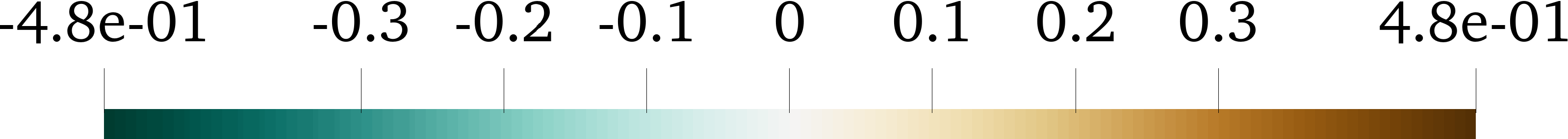}
  \end{minipage}
  \caption{\hyperref[sec.numerics:subsec.ex3]{Example 3} -- Difference between the solution at $t=0.3206$ and $0.322$. Computed with Taylor-Hood elements with $k=2$, $h=0.1$, $\Delta t=0.0014$, $\nu=0.01$, $\gpv=\gpp=0.1$, $\ggd=0$ and narrow-band ghost-penalty stabilisation.}
  \label{fig.ex3.picture_diff}
\end{figure}

While we observe no temporal oscillations in the divergence or vorticity within \Cref{fig.ex3.picture.div,fig.ex3.picture.vort}, these become noticeable when visualizing the difference between the first two times in \Cref{fig.ex3.picture_diff}. Additionally, more pronounced temporal oscillations appear in the vorticity, where pressure oscillations also occur.

\subsection{Other Variants of the Eulerian Unfitted FEM}
\label{sec.numerics:subsec.other-methods}

We also ran the same tests within the unfitted Eulerian time-stepping framework for several related methods all with the narrow band variant of the stabilization. The results were not significantly different, so instead of detailing all results, we provide a brief summary.

\subsubsection{BDF1 Time-Stepping}

The Eulerian unfitted finite element approach for moving domains has been analyzed in the literature using BDF1 time-stepping~\cite{LO19,vWRL21,BFM22,NO23}. We tested this less accurate variant, and although it produced slightly smaller spurious oscillations, the improvement was not substantial enough to justify the use of a first-order method.

\subsubsection{Equal-Order Discretization in Space}

In \cite{BFM22}, a variant of \eqref{eqn.discrete.method:TH} using equal-order unfitted elements for velocity and pressure was studied. This method incorporates continuous interior penalty (CIP) pressure stabilization in the bulk, which coincides with the "normal derivative jump" version of the ghost penalty operator for piecewise linear elements. Our experiments with this equal-order method revealed slightly larger spurious oscillations compared to the Taylor-Hood version for the oscillating circle case. This result aligns with expectations, as we have observed that spatial resolution is a major factor in the magnitude of spurious forces, and thus, using $P^1-P^1$ elements does not help. For the oscillating square case, where the solution is less regular, the equal-order choice showed slightly better results than Taylor-Hood without grad-div stabilization but performed worse than Taylor-Hood with grad-div stabilization.

\subsubsection{Scott-Vogelius Discretization in Space}

Divergence constraint and pressure variable are closely linked, and the lack of divergence conformity on the extended domain is currently a major challenge in proving FE pressure stability in stronger norms. This motivates the consideration of exactly divergence-free discretizations in space. The well-known Scott-Vogelius (SV) $P^k-P^{k-1}_{dc}$ element yields pointwise divergence-free finite element velocity in a steady domain setting. The SV element is inf-sup stable in the fitted setting for $k \geq 2d$ on meshes without singular vertices and for $k \geq d$ on barycentrically refined meshes.

In the unfitted setting, the SV element has been explored in \cite{LNO23} on Alfeld-split (barycentrically refined) meshes. However, pressure stabilization was achieved using the same narrow-band ghost-penalty operator:
\begin{equation*}
  j_h^n(p, q) = \sum_{F\in\FhnGAl}\int_{\omega_F}\jump{p}_\omega\jump{q}_\omega\dif \xb,
\end{equation*}
resulting in divergence errors on elements where the ghost-penalty operator is applied. To address this, a new mixed ghost-penalty type operator was introduced in \cite{FHNZ24} for an interface Darcy problem and later applied to the Stokes problem in \cite{FNZ23}. This mixed operator in $H^{div}$-conforming settings has the form:
\begin{equation}\label{aux1189}
  j_h^n(p, \diver(\vb)) + j_h^n(q, \diver(\ub)) =\sum_{F\in\FhndAl}\int_{\omega_F}\jump{p}_\omega\jump{\diver(\vb)}_\omega + \jump{q}_\omega\jump{\diver(\ub)}_\omega\dif \xb,
\end{equation}
where $\FhndAl$ is formed by selecting ghost-penalty elements based on the unrefined mesh and then considering all interior facets from the Alfeld-split elements. The fully discrete version of the scheme reads:
\begin{multline}\label{eqn.discrete.method:SV}
  \int_{\Ohn}\frac{3\uhn - 4\ub_h^{n-1} + \ub_h^{n-2}}{2 \Delta t} \vh \dif \xb + a_h^n(\uhn,\vh) + c_h^n(\uhn, \uhn, \vh) + b_{h, 0}^n(p_h^n, \vh) + b_h^n(q_h,\uhn)\\
   + s_{SV}^{gp, n}((\uhn, p_h^n), (\vh, q_h))= 0.
\end{multline}
Since the pressure ghost-penalty operator controls the solution in the $L^2$-norm, it ensures that the solution is exactly divergence-free over the entire active domain (see \cite[Theorem~1]{FNZ23}, which extends to the Scott-Vogelius setting). In this case, the velocity and pressure ghost-penalty facets must coincide. While an inf-sup stability result for this mixed ghost-penalty approach exists for the Darcy setting \cite{FHNZ24}, a similar result for the Stokes setting is currently unavailable.

Computational tests with this approach produced pointwise divergence-free solutions up to machine precision, but overall results were unsatisfactory. The convergence study for examples from \Cref{sec.numerics:subsec.convergence} revealed inconsistent performance for extension strips larger than $L=1$, and the direct linear solver\footnote{Pardiso via Intel MKL} occasionally failed due to poor conditioning of the matrix. Additionally, in the test cases from \Cref{sec.numerics:subsec.ex2,sec.numerics:subsec.ex3}, we observed spurious pressure oscillations one to two orders of magnitude larger than those in the Taylor-Hood results. Given the unresolved question of inf-sup stability for \eqref{eqn.discrete.method:SV}, further analysis in the stationary setting is necessary before advancing this method for moving boundaries.

Another well-understood variant of the Scott-Vogelius unfitted FEM, in the stationary setting, is the method in \cite{LNO23}. Extending this to the moving domain involves replacing the mixed ghost-penalty operator in \eqref{eqn.discrete.method:SV} with the pressure ghost-penalty operator $j_h(p,q)$, based on the facet set $\FhndAl$. This results in a solution that is exactly divergence-free in the bulk but not in the extension strip. To address this, a `heavy' grad-div stabilization is added, affecting only the extension strip solution. Computational tests with $\ggd=10h^{-1}$ (as used in \cite{LNO23}) showed no conditioning issues, unlike the mixed ghost-penalty term from \eqref{aux1189}. However, pressure boundary forces still exhibited spurious oscillations somewhat larger than those in the Taylor-Hood case.

\section{Conclusions and Outlook}
\label{sec.conclusion}

We investigated an Eulerian unfitted finite element method for the Navier-Stokes equations on moving domains, focusing on spurious temporal pressure oscillations --- a common issue in immersed boundary methods  often attributed to insufficient mass or volume conservation. The unfitted finite element framework, with its ability to precisely control these factors and its solid theoretical foundation, presents a promising approach to addressing the challenge of building an oscillation free IBM.

To gain analytical insights, we first examined the simpler stationary domain setting. Here, we demonstrated that controlling temporal pressure oscillations requires uniform control of the $L^\infty(L^2)$-norm of the pressure, which depends on the $L^\infty(H^1)$-norm of the velocity and the $L^\infty(L^2)$-norm of its discrete time derivative. Numerical experiments confirmed that control over the latter is sufficient, though not strictly necessary.

Identifying where the analysis fails in the case of to moving boundaries, we proposed a globally redefined ghost-penalty term to ensure unconditional stability of the instantaneous pressure. This modification addresses the limitations of the narrow-band variant, which restricts the penalty term to a small region near the interface.

Numerical experiments confirmed optimal convergence rates for both narrow-band and global ghost-penalty variants in the $L^2(L^2)$- and $L^2(H^1)$-norms of the velocity and the $L^2(L^2)$-norm of the pressure. However, the global variant exhibited unconditional stability in the $L^\infty(L^2)$-norm of the pressure, while the narrow-band variant showed $\mathcal{O}(\Delta t^{-1})$ growth in this norm for fixed mesh sizes.

The pressure drag coefficient, computed for two challenging test cases, revealed spurious oscillations in the narrow-band variant when new elements entered the active domain. These oscillations amplified as $\Delta t \to 0$ with fixed $h$, though they were less severe than in other state-of-the-art methods. Adjustments such as reducing geometric and mass-conservation errors, tuning stabilization parameters, or using higher-order elements mitigated the oscillations but did not fully eliminate them. Smaller viscosity coefficients increased oscillations, while spatial mesh refinement significantly reduced them.

While these conclusions pertain to the narrow-band variant, a key finding of this paper is that extending the ghost-penalty term globally ensures instantaneous pressure stability. This effectively suppresses spurious pressure oscillations and removes any dependence on $\Delta t$, while maintaining optimal convergence rates in the $L^2(L^2)$- and $L^2(H^1)$-norms for velocity and the $L^2(L^2)$-norm for pressure.

A detailed stability and convergence analysis of the global ghost-penalty variant in $L^\infty$ temporal norms remains an open problem. 

While this paper considers a prescribed interface motion, we expect these ideas to extend to cases where the interface motion is coupled to fluid flow, e.g., in fluid–structure interaction problems. In the latter case, reducing spurious hydrodynamic forces may be crucial for accurately predicting the structure’s motion. Extending and testing the unfitted FEM with global stabilization in such a setting is an interesting direction for future research.

\section*{Data Availability Statement}
The code used to realize the numerical examples is freely available on github \url{https://github.com/hvonwah/stable_inst_pressure_moving_domain_repro} and archived on zenodo \url{https://doi.org/10.5281/zenodo.14548166}.

\section*{Acknowledgment} 
M.O. was supported in part by the U.S. National Science Foundation under awards DMS-2309197 and DMS-2408978.
The authors would like to thank Guglielmo Scovazzi for bringing the issue of spurious pressure oscillations in immersed boundary methods to their attention.


\begin{thebibliography}{10}
\expandafter\ifx\csname url\endcsname\relax
  \def\url#1{\texttt{#1}}\fi
\expandafter\ifx\csname urlprefix\endcsname\relax\def\urlprefix{URL }\fi
\expandafter\ifx\csname href\endcsname\relax
  \def\href#1#2{#2} \def\path#1{#1}\fi

\bibitem{formaggia2010cardiovascular}
L.~Formaggia, A.~Quarteroni, A.~Veneziani, Cardiovascular {Mathematics}:
  Modeling and simulation of the circulatory system, Vol.~1, Springer Science
  \& Business Media, 2010.
\newblock \href {https://doi.org/10.1007/978-88-470-1152-6}
  {\path{doi:10.1007/978-88-470-1152-6}}.

\bibitem{bazilevs2013computational}
Y.~Bazilevs, K.~Takizawa, T.~E. Tezduyar, Computational fluid-structure
  interaction: Methods and applications, John Wiley \& Sons, 2013.
\newblock \href {https://doi.org/10.1002/9781118483565}
  {\path{doi:10.1002/9781118483565}}.

\bibitem{hirt1974arbitrary}
C.~W. Hirt, A.~A. Amsden, J.~L. Cook, An arbitrary {L}agrangian-{E}ulerian
  computing method for all flow speeds [reprint of {J}. {C}omput. {P}hys. {14}
  (1974), no. 3, 227--253], J. Comput. Phys. 135~(2) (1997) 198--216.
\newblock \href {https://doi.org/10.1006/jcph.1997.5702}
  {\path{doi:10.1006/jcph.1997.5702}}.

\bibitem{peskin1977numerical}
C.~S. Peskin, Numerical analysis of blood flow in the heart, J. Comput. Phys.
  25~(3) (1977) 220--252.
\newblock \href {https://doi.org/10.1016/0021-9991(77)90100-0}
  {\path{doi:10.1016/0021-9991(77)90100-0}}.

\bibitem{tezduyar1992new}
T.~E. Tezduyar, M.~Behr, S.~Mittal, J.~Liou, A new strategy for finite element
  computations involving moving boundaries and interfaces---the
  deforming-spatial-domain/space-time procedure. {II}. {C}omputation of
  free-surface flows, two-liquid flows, and flows with drifting cylinders,
  Comput. Methods Appl. Mech. Engrg. 94~(3) (1992) 353--371.
\newblock \href {https://doi.org/10.1016/0045-7825(92)90060-W}
  {\path{doi:10.1016/0045-7825(92)90060-W}}.

\bibitem{masud1997space}
A.~Masud, T.~J.~R. Hughes, A space-time {G}alerkin/least-squares finite element
  formulation of the {N}avier-{S}tokes equations for moving domain problems,
  Comput. Methods Appl. Mech. Engrg. 146~(1-2) (1997) 91--126.
\newblock \href {https://doi.org/10.1016/S0045-7825(96)01222-4}
  {\path{doi:10.1016/S0045-7825(96)01222-4}}.

\bibitem{glowinski1999distributed}
R.~Glowinski, T.-W. Pan, T.~I. Hesla, D.~D. Joseph, A distributed {L}agrange
  multiplier/fictitious domain method for particulate flows, Int. J. Multiph.
  Flow 25~(5) (1999) 755--794.
\newblock \href {https://doi.org/10.1016/S0301-9322(98)00048-2}
  {\path{doi:10.1016/S0301-9322(98)00048-2}}.

\bibitem{formaggia1999stability}
L.~Formaggia, F.~Nobile, \href{https://infoscience.epfl.ch/record/176278}{A
  stability analysis for the arbitrary {L}agrangian {E}ulerian formulation with
  finite elements}, East-West J. Numer. Math. 7~(2) (1999) 105--131.
\newline\urlprefix\url{https://infoscience.epfl.ch/record/176278}

\bibitem{duarte2004arbitrary}
F.~Duarte, R.~Gormaz, S.~Natesan, Arbitrary {L}agrangian-{E}ulerian method for
  {N}avier-{S}tokes equations with moving boundaries, Comput. Methods Appl.
  Mech. Engrg. 193~(45-47) (2004) 4819--4836.
\newblock \href {https://doi.org/10.1016/j.cma.2004.05.003}
  {\path{doi:10.1016/j.cma.2004.05.003}}.

\bibitem{gross2011numerical}
S.~Gross, A.~Reusken, Numerical methods for two-phase incompressible flows,
  Vol.~40, Springer Science \& Business Media, 2011.
\newblock \href {https://doi.org/10.1007/978-3-642-19686-7}
  {\path{doi:10.1007/978-3-642-19686-7}}.

\bibitem{uhlmann2005immersed}
M.~Uhlmann, An immersed boundary method with direct forcing for the simulation
  of particulate flows, J. Comput. Phys. 209~(2) (2005) 448--476.
\newblock \href {https://doi.org/10.1016/j.jcp.2005.03.017}
  {\path{doi:10.1016/j.jcp.2005.03.017}}.

\bibitem{seo2011sharp}
J.~H. Seo, R.~Mittal, A sharp-interface immersed boundary method with improved
  mass conservation and reduced spurious pressure oscillations, J. Comput.
  Phys. 230~(19) (2011) 7347--7363.
\newblock \href {https://doi.org/10.1016/j.jcp.2011.06.003}
  {\path{doi:10.1016/j.jcp.2011.06.003}}.

\bibitem{lee2011sources}
J.~Lee, J.~Kim, H.~Choi, K.-S. Yang, Sources of spurious force oscillations
  from an immersed boundary method for moving-body problems, J. Comput. Phys.
  230~(7) (2011) 2677--2695.
\newblock \href {https://doi.org/10.1016/j.jcp.2011.01.004}
  {\path{doi:10.1016/j.jcp.2011.01.004}}.

\bibitem{XCM24}
D.~Xu, O.~Colomés, A.~Main, K.~Li, N.~M. Atallah, N.~Abboud, G.~Scovazzi, A
  weighted shifted boundary method for immersed moving boundary simulations of
  stokes’ flow, J. Comput. Phys. 510 (2024) 113095.
\newblock \href {https://doi.org/10.1016/j.jcp.2024.113095}
  {\path{doi:10.1016/j.jcp.2024.113095}}.

\bibitem{mittal2008versatile}
R.~Mittal, H.~Dong, M.~Bozkurttas, F.~M. Najjar, A.~Vargas, A.~von Loebbecke, A
  versatile sharp interface immersed boundary method for incompressible flows
  with complex boundaries, J. Comput. Phys. 227~(10) (2008) 4825--4852.
\newblock \href {https://doi.org/10.1016/j.jcp.2008.01.028}
  {\path{doi:10.1016/j.jcp.2008.01.028}}.

\bibitem{luo20103d}
H.~Luo, B.~Yin, H.~Dai, J.~Doyle, A 3d computational study of the
  flow-structure interaction in flapping flight, in: 48th AIAA aerospace
  sciences meeting including the new horizons forum and aerospace exposition,
  2010, p. 556.
\newblock \href {https://doi.org/10.2514/6.2010-556}
  {\path{doi:10.2514/6.2010-556}}.

\bibitem{liao2010simulating}
C.-C. Liao, Y.-W. Chang, C.-A. Lin, J.~McDonough, Simulating flows with moving
  rigid boundary using immersed-boundary method, Comput Fluids 39~(1) (2010)
  152--167.
\newblock \href {https://doi.org/10.1016/j.compfluid.2009.07.011}
  {\path{doi:10.1016/j.compfluid.2009.07.011}}.

\bibitem{yang2009smoothing}
X.~Yang, X.~Zhang, Z.~Li, G.-W. He, A smoothing technique for discrete delta
  functions with application to immersed boundary method in moving boundary
  simulations, J. Comput. Phys. 228~(20) (2009) 7821--7836.
\newblock \href {https://doi.org/10.1016/j.jcp.2009.07.023}
  {\path{doi:10.1016/j.jcp.2009.07.023}}.

\bibitem{shirgaonkar2009new}
A.~A. Shirgaonkar, M.~A. MacIver, N.~A. Patankar, A new mathematical
  formulation and fast algorithm for fully resolved simulation of
  self-propulsion, J. Comput. Phys. 228~(7) (2009) 2366--2390.
\newblock \href {https://doi.org/10.1016/j.jcp.2008.12.006}
  {\path{doi:10.1016/j.jcp.2008.12.006}}.

\bibitem{lee2013implicit}
J.~Lee, D.~You, An implicit ghost-cell immersed boundary method for simulations
  of moving body problems with control of spurious force oscillations, J.
  Comput. Phys. 233 (2013) 295--314.
\newblock \href {https://doi.org/10.1016/j.jcp.2012.08.044}
  {\path{doi:10.1016/j.jcp.2012.08.044}}.

\bibitem{schneiders2013accurate}
L.~Schneiders, D.~Hartmann, M.~Meinke, W.~Schr{\"o}der, An accurate moving
  boundary formulation in cut-cell methods, J. Comput. Phys. 235 (2013)
  786--809.
\newblock \href {https://doi.org/10.1016/j.jcp.2012.09.038}
  {\path{doi:10.1016/j.jcp.2012.09.038}}.

\bibitem{ruberg2014fixed}
T.~Rüberg, F.~Cirak, A fixed‐grid b‐spline finite element technique for
  fluid–structure interaction, Internat. J. Numer. Methods Fluids 74~(9)
  (2013) 623--660.
\newblock \href {https://doi.org/10.1002/fld.3864}
  {\path{doi:10.1002/fld.3864}}.

\bibitem{goza2016accurate}
A.~Goza, S.~Liska, B.~Morley, T.~Colonius, Accurate computation of surface
  stresses and forces with immersed boundary methods, J. Comput. Phys. 321
  (2016) 860--873.
\newblock \href {https://doi.org/10.1016/j.jcp.2016.06.014}
  {\path{doi:10.1016/j.jcp.2016.06.014}}.

\bibitem{kim2019immersed}
W.~Kim, H.~Choi, Immersed boundary methods for fluid-structure interaction: A
  review, Int. J. Heat Fluid Fl. 75 (2019) 301--309.
\newblock \href {https://doi.org/10.1016/j.ijheatfluidflow.2019.01.010}
  {\path{doi:10.1016/j.ijheatfluidflow.2019.01.010}}.

\bibitem{BCH14}
E.~Burman, S.~Claus, P.~Hansbo, M.~G. Larson, A.~Massing, {CutFEM}:
  {Discretizing} geometry and partial differential equations, Internat. J.
  Numer. Methods Engrg. 104~(7) (2014) 472--501.
\newblock \href {https://doi.org/10.1002/nme.4823}
  {\path{doi:10.1002/nme.4823}}.

\bibitem{BH10}
E.~Burman, P.~Hansbo, Fictitious domain finite element methods using cut
  elements: {I}. {A} stabilized {Lagrange} multiplier method, Comput. Methods
  Appl. Mech. Engrg. 199~(41-44) (2010) 2680--2686.
\newblock \href {https://doi.org/10.1016/j.cma.2010.05.011}
  {\path{doi:10.1016/j.cma.2010.05.011}}.

\bibitem{BH12}
E.~Burman, P.~Hansbo, Fictitious domain finite element methods using cut
  elements: {II}. {A} stabilized {Nitsche} method, Appl. Numer. Math. 62~(4)
  (2012) 328--341.
\newblock \href {https://doi.org/10.1016/j.apnum.2011.01.008}
  {\path{doi:10.1016/j.apnum.2011.01.008}}.

\bibitem{Bur10}
E.~Burman, Ghost penalty, C.R. Math. 348~(21-22) (2010) 1217--1220.
\newblock \href {https://doi.org/10.1016/j.crma.2010.10.006}
  {\path{doi:10.1016/j.crma.2010.10.006}}.

\bibitem{olshanskii2002low}
M.~A. Olshanskii, A low order {Galerkin} finite element method for the
  {Navier}--{Stokes} equations of steady incompressible flow: a stabilization
  issue and iterative methods, Comput. Methods Appl. Mech. Engrg. 191~(47-48)
  (2002) 5515--5536.
\newblock \href {https://doi.org/10.1016/S0045-7825(02)00513-3}
  {\path{doi:10.1016/S0045-7825(02)00513-3}}.

\bibitem{LNO23}
H.~Liu, M.~Neilan, M.~Olshanskii, A {cutFEM} divergence{\textendash}free
  discretization for the {Stokes} problem, ESAIM Math. Model. Numer. Anal.
  57~(1) (2023) 143--165.
\newblock \href {https://doi.org/10.1051/m2an/2022072}
  {\path{doi:10.1051/m2an/2022072}}.

\bibitem{FNZ23}
T.~Frachon, E.~Nilsson, S.~Zahedi, Divergence-free cut finite element methods
  for {Stokes} flow, BIT Numerical Mathematics 64~(4) (Apr. 2023).
\newblock \href {https://doi.org/10.1007/s10543-024-01040-x}
  {\path{doi:10.1007/s10543-024-01040-x}}.

\bibitem{hansbo2016cut}
P.~Hansbo, M.~G. Larson, S.~Zahedi, A cut finite element method for coupled
  bulk-surface problems on time-dependent domains, Comput. Methods Appl. Mech.
  Engrg. 307 (2016) 96--116.
\newblock \href {https://doi.org/10.1016/j.cma.2016.04.012}
  {\path{doi:10.1016/j.cma.2016.04.012}}.

\bibitem{LO19}
C.~Lehrenfeld, M.~A. Olshanskii, An {Eulerian} finite element method for {PDE}s
  in time-dependent domains, ESAIM Math. Model. Numer. Anal. 53~(2) (2019)
  585--614.
\newblock \href {https://doi.org/10.1051/m2an/2018068}
  {\path{doi:10.1051/m2an/2018068}}.

\bibitem{claus2019cutfem}
S.~Claus, P.~Kerfriden, A {CutFEM} method for two-phase flow problems, Comput.
  Methods Appl. Mech. Engrg. 348 (2019) 185--206.
\newblock \href {https://doi.org/10.1016/j.cma.2019.01.009}
  {\path{doi:10.1016/j.cma.2019.01.009}}.

\bibitem{von2021falling}
H.~von Wahl, T.~Richter, S.~Frei, T.~Hagemeier, Falling balls in a viscous
  fluid with contact: Comparing numerical simulations with experimental data,
  Phys. Fluids 33~(3) (2021).
\newblock \href {https://doi.org/10.1063/5.0037971}
  {\path{doi:10.1063/5.0037971}}.

\bibitem{AB22}
M.~Anselmann, M.~Bause, Cut finite element methods and ghost stabilization
  techniques for space-time discretizations of the {Navier}-{Stokes} equations,
  Internat. J. Numer. Methods Fluids (Mar. 2022).
\newblock \href {https://doi.org/10.1002/fld.5074}
  {\path{doi:10.1002/fld.5074}}.

\bibitem{BFM22}
E.~Burman, S.~Frei, A.~Massing, Eulerian time-stepping schemes for the
  non-stationary {Stokes} equations on time-dependent domains, Numer. Math.
  150~(2) (2022) 423--478.
\newblock \href {https://doi.org/10.1007/s00211-021-01264-x}
  {\path{doi:10.1007/s00211-021-01264-x}}.

\bibitem{vWRL21}
H.~von Wahl, T.~Richter, C.~Lehrenfeld, An unfitted {Eulerian} finite element
  method for the time-dependent {Stokes} problem on moving domains, IMA J.
  Numer. Anal. 42~(3) (2021) 2505--2544.
\newblock \href {https://doi.org/10.1093/imanum/drab044}
  {\path{doi:10.1093/imanum/drab044}}.

\bibitem{NO23}
M.~Neilan, M.~Olshanskii, An {Eulerian} finite element method for the
  linearized {Navier}–{Stokes} problem in an evolving domain, IMA J. Numer.
  Anal. (Jan. 2023).
\newblock \href {https://doi.org/10.1093/imanum/drad105}
  {\path{doi:10.1093/imanum/drad105}}.

\bibitem{Leh16}
C.~Lehrenfeld, High order unfitted finite element methods on level set domains
  using isoparametric mappings, Comput. Methods Appl. Mech. Engrg. 300 (2016)
  716--733.
\newblock \href {https://doi.org/10.1016/j.cma.2015.12.005}
  {\path{doi:10.1016/j.cma.2015.12.005}}.

\bibitem{LL21}
Y.~Lou, C.~Lehrenfeld, Isoparametric unfitted {BDF}--finite element method for
  {PDEs} on evolving domains, SIAM J. Numer. Anal. 60~(4) (2021) 2069--2098.
\newblock \href {https://doi.org/10.1137/21m142126x}
  {\path{doi:10.1137/21m142126x}}.

\bibitem{BH14}
E.~Burman, P.~Hansbo, Fictitious domain methods using cut elements: {III}. {A}
  stabilized {Nitsche} method for {Stokes}' problem, ESAIM Math. Model. Numer.
  Anal. 48~(3) (2014) 859--874.
\newblock \href {https://doi.org/10.1051/m2an/2013123}
  {\path{doi:10.1051/m2an/2013123}}.

\bibitem{MLLR14}
A.~Massing, M.~G. Larson, A.~Logg, M.~E. Rognes, A stabilized {Nitsche}
  fictitious domain method for the {Stokes} problem, J. Sci. Comput. 61~(3)
  (2014) 604--628.
\newblock \href {https://doi.org/10.1007/s10915-014-9838-9}
  {\path{doi:10.1007/s10915-014-9838-9}}.

\bibitem{GO17}
J.~Guzm{\'{a}}n, M.~A. Olshanskii, Inf-sup stability of geometrically unfitted
  {Stokes} finite elements, Math. Comp. 87~(313) (2017) 2091--2112.
\newblock \href {https://doi.org/10.1090/mcom/3288}
  {\path{doi:10.1090/mcom/3288}}.

\bibitem{Pre18}
J.~Preu\ss, Higher order unfitted isoparametric space-time {FEM} on moving
  domains, Master's thesis, Georg-August-Universit\"at G\"ottingen (2018).
\newblock \href {https://doi.org/10.25625/UACWXS} {\path{doi:10.25625/UACWXS}}.

\bibitem{BNV22}
S.~Badia, E.~Neiva, F.~Verdugo, Linking ghost penalty and aggregated unfitted
  methods, Comput. Methods Appl. Mech. Engrg. 388 (2022) 114232.
\newblock \href {https://doi.org/10.1016/j.cma.2021.114232}
  {\path{doi:10.1016/j.cma.2021.114232}}.

\bibitem{BVM18}
S.~Badia, F.~Verdugo, A.~F. Mart{\'{\i}}n, The aggregated unfitted finite
  element method for elliptic problems, Comput. Methods Appl. Mech. Engrg. 336
  (2018) 533--553.
\newblock \href {https://doi.org/10.1016/j.cma.2018.03.022}
  {\path{doi:10.1016/j.cma.2018.03.022}}.

\bibitem{MKO13}
B.~Müller, F.~Kummer, M.~Oberlack, Highly accurate surface and volume
  integration on implicit domains by means of moment-fitting, Internat. J.
  Numer. Methods Engrg. 96~(8) (2013) 512--528.
\newblock \href {https://doi.org/10.1002/nme.4569}
  {\path{doi:10.1002/nme.4569}}.

\bibitem{Say15}
R.~I. Saye, High-order quadrature methods for implicitly defined surfaces and
  volumes in hyperrectangles, SIAM J. Sci. Comput. 37~(2) (2015) A993--A1019.
\newblock \href {https://doi.org/10.1137/140966290}
  {\path{doi:10.1137/140966290}}.

\bibitem{OS16}
M.~A. Olshanskii, D.~Safin, Numerical integration over implicitly defined
  domains for higher order unfitted finite element methods, Lobachevskii J.
  Math. 37~(5) (2016) 582--596.
\newblock \href {https://doi.org/10.1134/s1995080216050103}
  {\path{doi:10.1134/s1995080216050103}}.

\bibitem{FOSS17}
T.~P. Fries, S.~Omerovi{\'{c}}, D.~Schöllhammer, J.~Steidl, Higher-order
  meshing of implicit geometries{\textemdash}part {I}: {I}ntegration and
  interpolation in cut elements, Comput. Methods Appl. Mech. Engrg. 313 (2017)
  759--784.
\newblock \href {https://doi.org/10.1016/j.cma.2016.10.019}
  {\path{doi:10.1016/j.cma.2016.10.019}}.

\bibitem{HSK17}
S.~Hubrich, P.~D. Stolfo, L.~Kudela, S.~Kollmannsberger, E.~Rank, A.~Schröder,
  A.~Düster, Numerical integration of discontinuous functions: moment fitting
  and smart octree, Comput. Mech. 60~(5) (2017) 863--881.
\newblock \href {https://doi.org/10.1007/s00466-017-1441-0}
  {\path{doi:10.1007/s00466-017-1441-0}}.

\bibitem{Sch97}
J.~Sch\"oberl, {NETGEN} an advancing front {2D}/{3D}-mesh generator based on
  abstract rules, Comput. Vis. Sci. 1~(1) (1997) 41--52.
\newblock \href {https://doi.org/10.1007/s007910050004}
  {\path{doi:10.1007/s007910050004}}.

\bibitem{Sch14}
J.~Sch\"oberl, C++11 implementation of finite elements in {}{NGSolve}{}, Tech.
  rep. (Sep. 2014).

\bibitem{LHPvW21}
C.~Lehrenfeld, F.~Heimann, J.~Preu\ss, H.~von Wahl,
  \href{https://github.com/ngsxfem/ngsxfem}{ngsxfem: Add-on to {NGSolve} for
  geometrically unfitted finite element discretizations}, J. Open Source Softw.
  6~(64) (2021) 3237.
\newblock \href {https://doi.org/10.21105/joss.03237}
  {\path{doi:10.21105/joss.03237}}.
\newline\urlprefix\url{https://github.com/ngsxfem/ngsxfem}

\bibitem{FHNZ24}
T.~Frachon, P.~Hansbo, E.~Nilsson, S.~Zahedi, A divergence preserving cut
  finite element method for {Darcy} flow, SIAM J. Sci. Comput. 46~(3) (2024)
  A1793--A1820.
\newblock \href {https://doi.org/10.1137/22m149702x}
  {\path{doi:10.1137/22m149702x}}.

\end{thebibliography}
\end{document}